\documentclass[12pt,oneside]{book}
\usepackage{geometry}                		
\usepackage{graphicx}									
\usepackage{amssymb}\vspace{0.25cm}
\usepackage{mathtools}
\usepackage{setspace}
\usepackage{tikz-cd}
\usepackage[mathscr]{euscript}
\newcommand{\abs}[1]{\left\vert#1\right\vert}
\usepackage[normalem]{ulem}
\usepackage{manfnt}
\usepackage{pgfplots}
\usepackage{pifont}
\usepackage[safe]{tipa}
\usepackage{marvosym}
\usepackage{scalerel,stackengine}
\usepackage{titlesec}
\usepackage{makeidx}
\usepackage{centernot}
\usepackage{xspace}
\usepackage[noauto]{chappg}
\usepackage[OT2,T1]{fontenc}
\usepackage{hyperref}
\usepackage{amsmath}
\usepackage[amsmath,thmmarks]{ntheorem}
\usepackage{tabularx}
\usepackage{yfonts}
\usepackage{epstopdf}
\usepackage{pdfpages}
\newcolumntype{Y}{>{\raggedright\arraybackslash}X}

\newcommand\jA{%
\textgoth{A}
}
\newcommand\jAt{%
\textgoth{At}
}
\newcommand\jAlg{%
\textgoth{Alg}
}
\newcommand\jFil{%
\textgoth{F}
\hspace{0.02cm}
\textgoth{i}
\hspace{0.02cm}
\textgoth{l}
}

\newcommand\jB{%
\textgoth{B}
}
\newcommand\jBa{%
\textgoth{Ba}
}
\newcommand\jBo{%
\textgoth{Bo}
}
\newcommand\jD{%
\textgoth{Dyn}
}
\newcommand\jDyn{%
\textgoth{D}
}
\newcommand\jF{%
\textgoth{F}
}
\newcommand\jG{%
\textgoth{G}
}
\newcommand\jK{%
\textgoth{K}
}
\newcommand\jKol{%
\textgoth{Kol}
}
\newcommand\jL{%
\textgoth{L}
}
\newcommand\jLat{%
\textgoth{Lat}
}
\newcommand\jo{%
\textgoth{o}
}
\newcommand\jP{%
\textgoth{P}
}
\newcommand\jS{%
\textgoth{S}
}
\newcommand\jRin{%
\textgoth{Rin}
}
\newcommand\jT{%
\textgoth{T}
}
\newcommand\jTop{%
\textgoth{Top}
}
\newcommand\jX{%
\textgoth{X}
}
\newcommand\jZ{%
\textgoth{Z}
}

\newcommand\jstar{%
\star
}

\newcommand\fc{%
\mathfrak{c}
}
\newcommand\fF{%
\mathfrak{F}
}
\newcommand\fG{%
\mathfrak{G}
}

\newcommand\fm{%
\mathfrak{m}
}

\newcommand\fN{%
\mathfrak{N}
}

\newcommand\fp{%
\mathfrak{p}
}
\newcommand\fq{%
\mathfrak{q}
}
\newcommand\fS{%
\mathfrak{S}
}

\newcommand\gZ{%
\textgoth{Z}
}

\newcommand\A{%
\mathbb{A}
}

\newcommand\Cx{%
\mathbb{C}
}
\newcommand\E{%
\mathbb{E}
}

\newcommand\K{%
\mathbb{K}
}

\newcommand\N{%
\mathbb{N}
}

\newcommand\PP{%
\mathbb{P}
}
\newcommand\Q{%
\mathbb{Q}
}
\newcommand\R{%
\mathbb{R}
}

\newcommand\X{%
\mathbb{X}
}
\newcommand\XX{%
\mathbb{X}
}
\newcommand\Y{%
\mathbb{Y}
}
\newcommand\YY{%
\mathbb{Y}
}
\newcommand\Z{%
\mathbb{Z}
}

\newcommand\bBR{%
\textbf{BR}\xspace
}

\newcommand\bBS{%
\textbf{BS}\xspace
}

\newcommand\bRS{%
\textbf{RS}\xspace
}

\newcommand\bSR{%
\textbf{SR}\xspace
}

\newcommand\Top{%
\textbf{Top}\xspace
}
\newcommand\TopsubS{%
\textbf{Top$_\text{S}$}\xspace
}

\newcommand\Spec{%
\text{Spec}
}

\newcommand\Hom{%
\text{Hom}
}

\newcommand\ST{%
\text{ST}\xspace
}

\newcommand\sI{%
\mathcal{I}
}
\newcommand\sJ{%
\mathcal{J}
}

\newcommand\sP{%
\mathcal{P}
}

\newcommand\sS{%
\mathcal{S}
}

\newcommand\dash{%
\text{-}
}

\newcommand\Fnc{%
\text{Fnc}
}
\newcommand\Par{%
\text{Par}
}
\newcommand\Seq{%
\text{Seq}
}

\newcommand\Com{%
\text{Com}\hspace{0.03cm}
}

\newcommand\loc{%
\text{loc}
}

\newcommand\Pro{%
\text{Pro}
}

\newcommand\nth{%
\text{th}
}

\newcommand\td{%
\text{d}
}

\newcommand\tA{%
\text{A}
}
\newcommand\At{%
\text{At}
}
\newcommand\tB{%
\text{B}
}
\newcommand\tc{%
\text{c}
}
\newcommand\tC{%
\text{C}
}

\newcommand\tK{%
\text{K}
}

\newcommand\tM{%
\text{M}
}

\newcommand\subsetx{%
\ \subset \ 
}

\newcommand\tb{%
\text{b}
}

\newcommand\tS{%
\text{S}
}
\newcommand\ts{%
\text{s}
}
\newcommand\tsd{%
\text{sd}
}
\newcommand\tds{%
\text{ds}
}

\newcommand\tT{%
\text{T}
}
\newcommand\tu{%
\text{u}
}

\newcommand\perf{%
\text{perf}
}

\newcommand\scat{%
\text{scat}
}

\newcommand\isol{%
\text{isol}
}

\newcommand\dom{%
\text{dom}
}

\def\thereforex{\boldsymbol{\text{ }
\leavevmode
\lower0.4ex\hbox{\textbullet}
\kern-.9em\raise1.1ex\hbox{\textbullet}
\kern-0.9em\lower0.4ex\hbox{\textbullet}
\hspace{0.1cm}\thinspace\text{ }}}
\def\thereforez{\boldsymbol{\text{ }
\leavevmode
\lower0.4ex\hbox{$\circ$}
\kern-.9em\raise1.1ex\hbox{$\circ$}
\kern-0.9em\lower0.4ex\hbox{$\circ$}
\hspace{0.1cm}\thinspace\text{ }}}
\newcommand\ra{%
\rightarrow
}

\newcommand\ds{%
\displaystyle
}

\newcommand\card{%
\text{card}\hspace{0.05cm}
}

\newcommand\diam{%
\text{diam}\hspace{0.05cm}
}

\newcommand\sgn{%
\text{sgn}\hspace{0.05cm}
}

\newcommand\tr{%
\text{tr}\hspace{0.05cm}
}
\newcommand\trx{%
\text{r}
}
\newcommand\rs{%
\text{r\hspace{0.02cm}s}
}

\newcommand\un[1]{%
\underline{#1}\xspace
}

\newcommand\unz[1]{
\underline{#1}\xspace
}

\newcommand\hsx{%
\hspace{0.05cm}
}
\newcommand\hsy{%
\hspace{0.03cm}
}

\newcommand\boxtimesdmc{%
\raisebox{-.075cm}{$\boxtimes$}
}

\newcommand{\norm}[1]{\left\lVert #1 \right\rVert}

\newcommand\ov[1]{%
\overline{#1}
}

\newcommand\ovs[1]{
\mkern 1.5mu
\overline{\mkern-2.5mu\mbox{\hspace{0.03cm}$#1$}\raisebox{2.45mm}{}\mkern-1.5mu}
\mkern 1.5mu
}

\newcommand\chisubQ{\chi\raisebox{-.1cm}{$\scaleto{_{\Q}}{6.5pt}$}}
\newcommand\chisubR{\chi\raisebox{-.1cm}{$\scaleto{_{R}}{6.5pt}$}}
\newcommand\chisubS{\chi\raisebox{-.1cm}{$\scaleto{_{S}}{6.5pt}$}}
\newcommand\chisubT{\chi\raisebox{-.1cm}{$\scaleto{_{T}}{6.5pt}$}}

\newcommand\chisubsS{\chi\raisebox{-.1cm}{$\scaleto{_{\sS}}{6.5pt}$}}
\newcommand\chisubSi{\chi\raisebox{-.1cm}{$\scaleto{_{S_i}}{6.5pt}$}}
\newcommand\chisubSsubi{\chi\raisebox{-.1cm}{$\scaleto{_{S_i}}{6.5pt}$}}

\newcommand\chisublimsupSsubi{\chi\raisebox{-.1cm}{$\scaleto{_{\varlimsup \ S_i}}{8.5pt}$}}
\newcommand\chisubliminfSsubi{\chi\raisebox{-.1cm}{$\scaleto{_{\varliminf \ S_i}}{8.5pt}$}}

\newcommand\chisubSminusT{\chi\raisebox{-.1cm}{$\scaleto{_{S - T}}{8.5pt}$}}
\newcommand\chisubSDeltaT{\chi\raisebox{-.1cm}{$\scaleto{_{S \hsy \Delta \hsy  T}}{8.5pt}$}}

\newcommand\chisubScapT{\chi\raisebox{-.1cm}{$\scaleto{_{S \cap T}}{8.5pt}$}}
\newcommand\chisubScupT{\chi\raisebox{-.1cm}{$\scaleto{_{S \cup T}}{8.5pt}$}}

\newcommand\chisubXX{\chi\raisebox{-.1cm}{$\scaleto{_{\XX}}{6.5pt}$}}
\newcommand\chisubXsubi{\chi\raisebox{-.1cm}{$\scaleto{_{X_i}}{6.5pt}$}}

\definecolor{ultramarine}{RGB}{0, 32, 96}

\definecolor{darkcerulean}{rgb}{0.3, 0.27, 0.49}

\definecolor{forestgreen}{rgb}{0.0, 0.27, 0.13}
\definecolor{forestgreenweb}{rgb}{0.13, 0.55, 0.13}
\definecolor{deepjunglegreen}{rgb}{0.0, 0.29, 0.29}

\definecolor{midnightblue}{rgb}{0.1, 0.1, 0.44}
\definecolor{midnightgreen}{rgb}{0.0, 0.29, 0.33}
\definecolor{myrtle}{rgb}{0.13, 0.26, 0.12}
\definecolor{darkviolet}{rgb}{0.58, 0.0, 0.83}
\definecolor{darkgreen}{rgb}{0.0, 0.2, 0.13}
\definecolor{officegreen}{rgb}{0.0, 0.5, 0.0}

\definecolor{harvardcrimson}{rgb}{0.79, 0.0, 0.09}
\definecolor{hollywoodcerise}{rgb}{0.96, 0.0, 0.63}
\definecolor{debianred}{rgb}{0.84, 0.04, 0.33}
\definecolor{darkturquoise}{rgb}{0.0, 0.81, 0.82}

\definecolor{darktangernine}{rgb}{1.0, 0.66, 0.07}
\definecolor{aureolin}{rgb}{0.99, 0.93, 0.0}
\definecolor{canaryyellow}{rgb}{1.0, 0.94, 0.0}
\definecolor{amber}{rgb}{1.0, 0.75, 0.0}

\definecolor{urobilin}{rgb}{0.88, 0.68, 0.13}
\definecolor{uscgold}{rgb}{1.0, 0.8, 0.0}

\usepackage{scalerel,stackengine}
\stackMath
\newcommand\reallywidehat[1]{%
\savestack{\tmpbox}{\stretchto{%
  \scaleto{%
    \scalerel*[\widthof{\ensuremath{#1}}]{\kern-.6pt\bigwedge\kern-.6pt}%
    {\rule[-\textheight/2]{1ex}{\textheight}}
  }{\textheight}%
}{0.5ex}}%
\stackon[1pt]{#1}{\tmpbox}%
}
\makeatletter
\DeclareRobustCommand\widecheck[1]{{\mathpalette\@widecheck{#1}}}
\def\@widecheck#1#2{%
    \setbox\z@\hbox{\m@th$#1#2$}%
    \setbox\tw@\hbox{\m@th$#1%
       \widehat{%
          \vrule\@width\z@\@height\ht\z@
          \vrule\@height\z@\@width\wd\z@}$}%
    \dp\tw@-\ht\z@
    \@tempdima\ht\z@ \advance\@tempdima2\ht\tw@ \divide\@tempdima\thr@@
    \setbox\tw@\hbox{%
       \raise\@tempdima\hbox{\scalebox{1}[-1]{\lower\@tempdima\box
\tw@}}}%
    {\ooalign{\box\tw@ \cr \box\z@}}}
\makeatother

\newtheoremstyle{xx}
  {4pt}
  {0pt}
  {\upshape}
  {\bfseries}
  {}
  { }
  {}
  
\makeatletter 
 \newtheoremstyle{myu}%
  {\upshape\item[ \indent\indent\bf\underline{\theorem@headerfont ##2:}]}%
\makeatother
\makeatletter 
 \newtheoremstyle{myn}%
  {\item[\hskip\labelsep \ \bf ##1 \theorem@headerfont ##2.]}%
\makeatother
\theoremstyle{myn}
\theoremstyle{myu}
{\upshape}

 \newtheoremstyle{mr}%
  {\upshape\item[ \indent{\theorem@headerfont ##2. \hspace{.2cm}}]}%
\makeatother
\theoremstyle{mr}
{\upshape}

\title
{
\textbf
{
$\mathcal{SETS \  AND \  CLASSES :}$
\\
$\mathcal{OPERATIONAL \ \ THEORY}$
}
}
\author{Garth Warner\\
Department of Mathematics\\
University of Washington}
\date{}	

\titleformat{\chapter}[display]
{\normalfont\filcenter\huge\bfseries}{}{0pt}{\large}

\titleformat{\chapter}[display]
{\normalfont\filcenter\huge\bfseries}{}{0pt}{\large}

\usepackage[OT2,OT1]{fontenc} 
\newcommand\cyr
{
\renewcommand\rmdefault{wncyr} 
\renewcommand\sfdefault{wncyss} 
\renewcommand\encodingdefault{OT2} 
\normalfont
\selectfont
}

\DeclareTextFontCommand{\textcyr}{\cyr}

\makeindex 
\begin{document}

\maketitle                              

\titlespacing*{\chapter}{0pt}{-50pt}{40pt}
\setlength{\parskip}{0.1em}

\begingroup
\fontsize{11pt}{11pt}\selectfont







\[
\textbf{ABSTRACT}
\]
\\

The purpose of this book is to lay out certain aspects of descriptive set theory.  
After initially establishing notation and generalities we proceed to the following topics:
partitions, 
semirings, 
rings, $\sigma$-rings, $\delta$-rings, 
products and sums, 
extension and generation.  
Extensive references and historical comments are included at the end of each section, as are further examples in the form of exercises and problems.  
\\[2cm]

\[
\textbf{ACKNOWLEDGEMENT}
\]

Many thanks to David Clark for his rendering the original transcript into AMS-LaTeX.  
Both of us also thank Judith Clare for her meticulous proofreading.
\newpage

\[
\textbf{CONTENTS}
\]


\hspace{2.1cm} \ \S1. \quad GENERALITIES%
\\[-.25cm]

\hspace{3.7cm} Exercises
\\[-.25cm]

\hspace{3.7cm} Problems
\\

\hspace{2.1cm} \ \S2. \quad PARTITIONS%
\\[-.25cm]

\hspace{3.7cm} Exercises
\\[-.25cm]

\hspace{3.7cm} Problems
\\

\hspace{2.1cm} \ \S3. \quad SEMIRINGS%
\\[-.25cm]

\hspace{3.7cm} Exercises
\\[-.25cm]

\hspace{3.7cm} Problems
\\

\hspace{2.1cm} \ \S4. \quad RINGS, $\sigma$-RINGS, $\delta$-RINGS%
\\[-.25cm]

\hspace{3.7cm} Exercises
\\[-.25cm]

\hspace{3.7cm} Problems
\\

\hspace{2.1cm} \ \S5. \quad PRODUCTS AND SUMS%
\\[-.25cm]

\hspace{3.7cm} Exercises
\\[-.25cm]

\hspace{3.7cm} Problems
\\

\hspace{2.1cm} \ \S6. \quad EXTENSION AND GENERATION%
\\[-.25cm]

\hspace{3.7cm} Exercises
\\[-.25cm]

\hspace{3.7cm} Problems

\newpage

\[
\textbf{Prerequisites}
\]

It is assumed that the reader is familiar with the language and notation employed in elementary algebra, analysis, set theory, and general topology.
\\

As usual
\\

$\N$ 
\quad the set of positive integers
\\

$\Z$ 
\quad the set of integers
\\

$\Q$ 
\quad the set of rational numbers
\\

$\PP$ 
\quad the set of irrational numbers
\\

$\R$ 
\quad the set of real numbers
\\

$\Cx$ 
\quad the set of complex numbers.
\\


The symbols $\N^n$, $\Z^n$, $\Q^n$, $\PP^n$, $\R^n$, $\Cx^n$ ($n$ a positive integer) are then to be assigned their customary interpretations.
\\

Tacitly, we shall always operate within the strictures of ZFC 
(Zermelo-Fraenkel Axioms + Axiom of Choice), 
unless the contrary is explicity stated.


\endgroup 

\pagenumbering{bychapter}
\setcounter{chapter}{0}
\pagenumbering{bychapter}
\chapter{
$\boldsymbol{\S}$\textbf{1}.\quad  Generalities}
\setlength\parindent{2em}
\setcounter{theoremn}{0}
\renewcommand{\thepage}{\S1-\arabic{page}}



\qquad 
Throughout this book, whenever the word \unz{set} is used, it is always understood to mean a subset of a given set which, 
generically, is denoted by $\XX$; 
we shall use the word \unz{class} for a set of sets and the word \unz{collection} for a set of classes.  
If $S$ and $T$ are subsets of $\XX$, then the union, intersection, difference, and symmetric difference of $S$ and $T$ 
are denoted by $S \cup T$, $S \cap T$, $S - T$, and $S \hsx \Delta \hsy T$, respectively.  
$\jP(X)$ stands for the class of all subsets of $\XX$; 
$\emptyset$ stands for the empty set. 
\\

By $\card (X)$, we mean the cardinality of $\XX$.  
A set is said to be 
\unz{countable} if its cardinality is $\aleph_0$, 
\unz{finite} if its cardinality is $< \aleph_0$, 
\unz{uncountable} if its cardinality is $> \aleph_0$, 
\unz{infinite} if not finite, i.e., countable or uncountable.  
If $F$ is a finite set, then $\# (F)$ is the number of elements in $F$.  
\\[-.5cm]

As is customary, 
\[
\aleph_0
\ < \ 
\aleph_1
\ < \ 
\ldots
\ < \ 
\aleph_\alpha
\ < \ 
\ldots
\]
are the \unz{infinite cardinals} and 
\[
\omega_0
\ < \ 
\omega_1
\ < \ 
\ldots
\ < \ 
\omega_\alpha
\ < \ 
\ldots
\]
are the \unz{infinite initial ordinals}.   
In this connection, bear in mind that $\alpha$ is an arbitrary ordinal and $\omega_\alpha$ is the first ordinal such that 
\[
\card (\{\beta \hsy : \hsy \beta < \omega_\alpha\}) 
\ = \ 
\aleph_\alpha.  
\]
Traditionally, $\omega_0$ is denoted by $\omega$, while $\omega_1$ is denoted by $\Omega$.  
By $\fc$, we understand the cardinality of the \unz{continuum}, i.e., $\fc = 2^{\aleph_0}$.  
The \unz{continuum hypothesis} is the statement that $2^{\aleph_0} = \aleph_1$; 
the \unz{generalized continuum hypothesis} is the statement that $2^{\aleph_\alpha} = \aleph_{\alpha + 1}$ for all ordinals $\alpha$.  
Both of these statements are independent of ZFC.
\\

The \unz{characteristic function} of a subset $S$ of $\XX$ is the function 
$\chisubS : \X \ra \R$ defined by
\[
\chisubS (x) \ = \ 
\begin{cases}
\ 1 \hspace{0,25cm} \text{if} \  x \in S
\\[8pt]
\ 0 \hspace{0,25cm} \text{if} \  x \in X - S
\end{cases}
.
\]
There is a canonical identification between $\jP (X)$ and the set $\Fnc(X, \{0, 1\})$ of all functions from $\XX$ to $\{0, 1\}$, 
namely the rule $S \mapsto \chisubS$.
\\


If $S$ and $T$ are subsets of $\XX$, then
\[
\begin{cases}
\ \chisubScapT \ = \ \min (\chisubS, \chisubT)
\\[8pt]
\ \chisubScupT \ = \ \max (\chisubS, \chisubT)
\end{cases}
\]
with $\chisubS \leq \chisubT$ iff $S \subset T$.  
Furthermore
\[
\begin{cases}
\ \chisubSminusT \hspace{0.4cm} = \ \chisubS (1 - \chisubT)
\\[8pt]
\ \chisubSDeltaT \ = \ \abs{\chisubS - \chisubT}
\end{cases}
.
\]

Let $\{S_i\}$ be a sequence of subsets of $\XX$ $-$then 
the set of all those points of $\XX$ which belong to $S_i$ for infinitely many values of $i$ is called the 
\unz{upper limit} 
or 
\unz{limit superior}  
of the sequence and is denoted by 
$\varlimsup \ S_i$ or $\limsup \ S_i$, 
while 
the set of all those points of $\XX$ which belong to $S_i$ for all but a finite number of values of $i$ is called the 
\unz{lower limit} 
or 
\unz{limit inferior}  
of the sequence and is denoted by 
$\varliminf \ S_i$ or $\liminf \ S_i$.  
Evidently, 
\[
\begin{cases}
\ \varlimsup \chisubSi 
\ = \ 
\bigcap\limits_{i = 1}^\infty \ 
\left(
\bigcup\limits_{m = i}^\infty \ 
S_m
\right)
\\[18pt]
\ \varliminf \chisubSi
\ = \ 
\bigcup\limits_{i = 1}^\infty \ 
\left(
\bigcap\limits_{m= i}^\infty \ 
S_m
\right)
\end{cases}
.
\]
In the event that
\[
\begin{cases}
\ \varlimsup \chisubSi 
\\[8pt]
\ \varliminf \chisubSi
\end{cases}
= \  S, 
\qquad \text{say,} 
\]
then $S$ is said to be the \unz{limit} of the sequence $S_1, S_2, \ldots$ and we write 
$S = \lim \ S_i$.  
For instance, if $\{S_i\}$ is an 
\unz{increasing} (\unz{decreasing}) 
sequence in the sense that 
$S_i \subset S_{i + 1}$ $(S_i \supset S_{i + 1})$ $\forall \ i$, then 
$\lim \ S_i = \bigcup \ S_i$ $(\bigcap \ S_i)$.  
In general, it is always true that 
\[
\bigcap \ S_i 
\subsetx 
\varliminf \ S_i 
\subsetx
\varlimsup \ S_i 
\subsetx
\bigcup \ S_i.
\]

In terms of characteristic functions, 
\[
\begin{cases}
\ \chisublimsupSsubi \ = \ \varlimsup \chisubSi
\\[8pt]
\ \chisubliminfSsubi \ = \ \varliminf \chisubSi
\end{cases}
.
\]
\\[-1cm]

\textbf{Example} \ 
Suppose that $\{S_i\}$ is a sequence of pairwise disjoint subsets of $\XX$ $-$then $\lim \ S_i = \emptyset$.
\\

The preceding notions can be interpreted topologically.  
For this purpose, it will be convenient to consider first the elements of a useful abstract construction. 
\\[-.5cm]

Thus let $(\XX, \jT)$ be a topological space $-$then by the 
\unz{sequential modification} 
of $(\XX, \jT)$  
we mean the topological space whose underlying set is still $\XX$ itself but whose topology $\jT_\tS$ consists, 
by definition, of the complements of those subsets $S$ of $\XX$ which are closed under pointwise convergence of sequences, 
i.e., a subset 
$S$ of $\XX$ is $\jT_\tS$-closed iff for every sequence $\{x_i\} \subset S$, $x_i \ra x$ $\implies x \in S$.  
It is easy to check that the class of closed subsets thereby singled out does in fact satisfy the usual axioms involved 
in defining a topology by closed sets.  
The canonical map 
$(\XX, \jT_\tS) \ra (\XX, \jT)$ is continuous, or, what amounts to the same, 
the $\jT_\tS$-topology on $\XX$ is finer than the $\jT$-topology.  
In addition, it is clear that a sequence $\{x_i\}$ in $\XX$ is $\jT$-convergent to a point $x$ iff it is $\jT_\tS$-convergent to $x$.  
These remarks enable one to characterize the sequential modification of $(\XX, \jT)$ in a simple way.  
Indeed, $\jT_\tS$ is the finest topology of all topologies $\jT_0 \supset \jT$ on $\XX$ which have the following property: \ 
A sequence in $\XX$ is $\jT$-convergent iff it is $\jT_0$-convergent.
\\

The essential significance of the sequential modification is contained in : 
\\

\textbf{Lemma 1} \quad 
Let $f : \XX \ra \Y$ be a map from $\XX$ into a topological space $\Y$ $-$then $f$ is continuous per $\jT_\tS$ iff 
$f$ is sequentially continuous per $\jT$.
\\


[We omit the elementary verification.]
\\

In connection with the preceding developments, a modicum of caution must be exercised, viz.: \ 
The $\jT_\tS$-closure of a subset $S$ of $\XX$ need not consist just of the sequential limit from $S$ but, in general, 
will be much larger, as can be seen by simple examples (cf. Exer. 8).  
This can easily be made precise.  
Given $S$, let $\tu S$ be the set of all $\jT$-limits of sequences in $S$.  
Putting $\tu_0 S = S$, define by transfinite recursion
\[
\tu_\alpha S 
\ = \ 
\tu 
\left(\hsx
\bigcup\limits_{\beta < \alpha} \ \tu_\beta S
\hsx\right)
\qquad (\alpha < \Omega).
\]
\\[-1cm]

\noindent
Then the $\jT_\tS$-closure of $S$ is $\bigcup\limits_{\alpha < \Omega} \ \tu_\alpha S$.  
Another way to look at it is to let $S_0$ run through those subsets of $S$ having cardinality $\leq \aleph_0$ 
$-$then the union of the $\jT_\tS$-closure of $S_0$ is the $\jT_\tS$-closure of $S$.  
In any event, the moral is that sequences do not ordinarily suffice; 
nets (or filters) will usually be needed.
\\[-.5cm]

[Note: \ 
Let \Top be the category whose objects are topological spaces and whose morphisms are continuous maps; 
let \TopsubS be the category whose objects are the sequential topological spaces, i.e., 
those topological spaces in which every sequentially closed subset is closed, 
and whose morphisms are continuous maps $-$then there is a canonically defined coreflective functor
\[
\Top \ra \TopsubS, 
\]
viz. the rule
\[
(\XX, \jT) \ra (\XX, \jT_\tS)
\]
together with the obvious assignment of morphisms.  
\TopsubS thus appears as a coreflective subcategory of \Top which, in fact, is monocoreflective, 
hence, on the basis of standard categorical generalities, is closed under the formation of quotients and coproducts in \Top.]

Suppose now that $\XX$ is again merely an abstract set but that $\Y$ is a topological space.  
Let $\Fnc(X, \Y)$ be the set of all functions from $\XX$ to $\Y$ equipped with the topology of pointwise convergence 
$-$then by $\Fnc(X, \Y)_\tS$ we understand the sequential modification of $\Fnc(X, \Y)$.  
The class of closed sets for the associated topology is thus comprised of those subsets of $\Fnc(X, \Y)$
which are closed under pointwise convergence of sequences. 
\\

\textbf{Example} \  
If $\XX$ and $\Y$ are both topological spaces, then the closure in $\Fnc(X, \Y)_\tS$ of the subset 
of all continuous maps is known as the class of \unz{Baire functions} (from $\XX$ to $\Y$).
\\

The identification 
$\jP(\XX) = \Fnc(X, \{0, 1\})$ enables one to topologize $\jP (\X)$ in a canonical way.  
Indeed, equipping $\{0, 1\}$ with the discrete topology, place on $\Fnc(X, \{0, 1\})$ the topology of pointwise convergence 
$-$then this topology may be pulled back to $\jP(\X)$, the upshot being that $\jP(\X)$ thus topologized is a compact Hausdorff space which, 
moreover, is totally disconnected.  
Write $\jP (\X)_\tS$ for the corresponding sequential modification $-$then $\jP (\X)_\tS$ is still Hausdorff and totally disconnected but, 
in general, need not be compact (cf. Exer. 12).  
Given a sequence $\{S_i\} \subset \jP (\XX)$, the relations 
\[
\begin{cases}
\ \chisublimsupSsubi \ = \ \varlimsup \chisubSi
\\[8pt]
\ \chisubliminfSsubi \ = \ \varliminf \chisubSi
\end{cases}
\]
then make it clear that $\lim \ S_i$ exists topologically, i.e., per $\jP(\XX)_\tS$, iff 
$\lim \ S_i$ exists in the sense that 
$\varlimsup \ S_i = \varliminf \ S_i$.  
\\

We shall terminate this \S \ with some definitions and related notation.  
\\[-.25cm]

Let $\jS$ be a nonempty subset of $\jP(X)$.  
Write $\jS_\ts$, $\jS_\sigma$, $\jS_\td$, $\jS_\delta$ for the class of subsets of $\XX$ comprised of all nonempty 
finite unions, 
countable unions, 
nonempty finite intersections, 
countable intersections, 
of sets in $\jS$ (repetitions being permissible); 
write $\jS_\Upsilon$ for the class of subsets of $\XX$ comprised of all sets in $\jS$ and all differences of sets in $\jS$; 
write $\jS_\tc$ for the class of subsets of $\XX$ comprised of all complements of sets in $\jS$.  
Successive applications of these operations is represented by juxtaposition of the symbols, e.g., 
$\jS_{\sigma \delta} \equiv (\jS_\sigma)_\delta$,  
the class of all countable intersections of countable unions of sets belonging to $\jS$.  
Obviously, 
\[
\begin{cases}
\ \jS \subsetx \jS_\ts \ = \ \jS_{\ts\ts} \subsetx \jS_\sigma \ = \ \jS_{\sigma\sigma}
\\[8pt]
\ \jS \subsetx \jS_\td \ = \ \jS_{\td\td} \subsetx \jS_\delta \ = \ \jS_{\delta\delta}
\end{cases}
. \qquad 
\begin{cases}
\ \jS 
\hspace{0.4cm} = \ 
\jS_{\tc\tc}, 
\\[8pt]
\jS_{\tc\sigma}
\hspace{0.2cm} = \ 
\jS_{\delta\tc}
\\[8pt]
\ \jS_{\tc\delta}
\ = \ 
\jS_{\sigma\tc}
\end{cases}
.
\]

The class $\jS$ is termed 
\unz{additive} 
(\unz{$\sigma$-additive}) 
if it is nonempty and closed under the formation of nonempty finite (countable) unions, 
i.e., provided $\jS = \jS_\ts$ ($\jS_\sigma$).  
The class $\jS$ is termed 
\unz{multiplicative} 
(\unz{$\delta$-multiplicative}) 
if it is nonempty and closed under the formation of nonempty finite (countable) intersections, 
i.e., provided $\jS = \jS_\td$ ($\jS_\delta$).   
If $\emptyset \in \jS$ and if $\jS$ is both additive and multiplicative ($\sigma$-additive or $\delta$-multiplicative), 
then $\jS$ is called a 
\unz{lattice} 
(\unz{$\sigma$-lattice}  or \unz{$\delta$-lattice}).  
Every $\sigma$-lattice or $\delta$-lattice is a lattice but, of course, not conversely.  
Naturally, a lattice of sets is an abstract lattice.
\\

\textbf{Example} \ 
Let $\XX$ be a topological space $-$then the class of all open (closed) subsets of $\XX$ is a $\sigma$-lattice ($\delta$-lattice).
\\

If $\jS$ is a nonempty subset of $\jP (\XX)$ and if $X_0$ is an arbitrary subset of $\XX$, 
then the \unz{trace} of $\jS$ on $X_0$ is the class
\[
\tr_{X_0} (\jS) 
\ = \ 
\{ S \cap X_0 \hsy : \hsy S \in \jS\}.
\]
The trace operation will preserve certain structures, e.g., the trace of a lattice is again a lattice.
\\[-.75cm]

\[
\textbf{\un{Notes and Remarks}}
\]

The notion of characteristic function is due to 
Ch. de la Vall\'ee Poussin\footnote[1]{\vspace{.11 cm}\textit{Trans. Amer. Math. Soc.}, \textbf{16} (1915), pp. 435-501.}.
Its use was, however, first anticipated by 
E$\acute .$ Borel\footnote[2]{\vspace{.11 cm}\un{Lecons sur la Th\'eorie des Fonctions}, Gauthier-Villars, Paris  (1898).}.  
E$\acute .$ 
Borel also introduced the upper limit and lower limit of a sequence of sets; cf.  
E$\acute .$ Borel
\footnote[3]{\vspace{.11 cm}\un{Lecons sur les Fonctions de Variables R\'eelles}, Gauthier-Villars, Paris  (1905).}.
(see p. 18).
Here
\[
\begin{cases}
\ \text{upper limit = limite compl\`ete}
\\[4pt]
\ \text{lower limit = limite restreinte}
\end{cases}
.
\]

\noindent
Strangely enough, the limit of a sequence of sets was formalized only later, viz. by 
Ch. de la Vall\'ee Poussin (op. cit.), the term being limite unique, the notation being lim, and also, 
independently, by 
F. Hausdorff\footnote[4]{\vspace{.11 cm}\un{Gr\"undzuge der Mengenlehre}, Veit \& Comp., Leipzig,  (1914).}
 in his classic
\un{Gr\"undzuge der Mengenlehre}, 
where also will be found the limit superior, limit inferior terminology.  
The notation $\varlimsup$ and $\varliminf$ was codified by 
Ch. de la Vall\'ee Poussin\footnote[5]{\vspace{.11 cm}\un{Int\' egrales de Lebesgue, Fonctions d'Ensemble, Classes de Baire}, Gauthier-Villars,Paris(1916),  p. 8}.
For an exhaustive study of the closure operations and their modifications, consult 
E. \v Cech\footnote[6]{\vspace{.11 cm}\un{Topological Spaces}, Academia, Prague,  (1966).}.
The topologization of $\jP(\XX)$ is the subject of a paper by 
R. Bagley\footnote[7]{\vspace{.11 cm}\textit{Michigan Math. J.}, \textbf{3} (1955-56), pp. 105-108.}
see also 
L. Savel\'ev\footnote[8]{{\fontencoding{OT2}\selectfont
[L. Savelev]}
\vspace{.11 cm}
\textit{Sibirsk. Math. \u Z.}, \textbf{6} (1965), pp. 1357-1364.}.
An elementary but useful survey (with extensive references) on the various operations
$\jS_\ts$, $\jS_\sigma$, $\jS_\td$, $\jS_\delta$, $\jS_\trx$ ,$\jS_\tc$ (and much more) has been given by 
W. Sierpin\'ski\footnote[9]{\vspace{.11 cm}\textit{Proc. Benares Math. Soc.}, N. S. \textbf{9} (1947), pp. 1-24.}.
The origin of the various subscripts used therein is this
\[
\begin{cases}
\ \ts, \hsy \sigma: \  \text{Summe}
\\[4pt]
\ \td, \hsy \delta: \  \text{Durchschnitt}
\end{cases}
,
\]

\noindent
r: relative \ (complement), \ c: complement.  \ 
Sierpin\'ski's\footnote[1]{\vspace{.11 cm}\un{Hypoth\` ese du Continu}, Chelsea, New York,  (1956).}
\un{Hypoth\` ese du Continu} 
is highly recommended as a source for additional information about the continuum hypothesis and its consequences.  
Many of the statements in this book have subsequently been approached from the point of view of Martin's axiom; 
cf. 
D. Martin and 
R. Solovay\footnote[2]{\vspace{.11 cm}\textit{Ann. Math. Logic}, \textbf{2} (1970), pp. 143-178.}.


\chapter{
$\boldsymbol{\S}$\textbf{1}.\quad  Exercises}
\setlength\parindent{2em}
\setcounter{theoremn}{0}
\renewcommand{\thepage}{\S1-\text{E-}\arabic{page}}



\qquad 
(1) \quad 
Let $S_i = [0, 1]$ for odd values of $i$ and $S_i = [-1, 0]$ for even values of $i$ $-$then 
$\varlimsup \ S_i = [-1, 1]$ and $\varliminf \ S_i = \{0\}$.
\\[-.25cm]

(2) \ 
Let $\{x_i\}$ be a sequence of real numbers; let $S_i = ]-\infty, x_i[$ $-$then 

\[
\begin{cases}
\ ]-\infty, \varlimsup \ x_i [ \subsetx  \varlimsup \ S_i \subsetx  ]-\infty, \varlimsup \ x_i ]
\\[11pt]
\ ]-\infty, \varliminf \ x_i [ \subsetx  \varliminf \ S_i \subsetx  ]-\infty, \varliminf \ x_i ]
\end{cases}
.
\]
\\[-1cm]

(3) \quad
Let $\{S_i\}$ , $\{S_i^\prime\}$, $\{S_i^{\prime\prime}\}$ be sequences of sets with 
$S_i \subset S_i^\prime \subset S_i^{\prime\prime}$ for all $i$.  
Suppose that 
$\lim S_i ^\prime= \lim S_i^{\prime\prime} = S$, say, $-$then $\lim S_i$ exists and is equal to $S$.
\\[-.25cm]

(4) \quad 
The union (intersection) of a sequence of sets $\{S_i\}$ can always be represented as the limit of an increasing (decreasing) sequence of sets.
\\[-.5cm]

[In fact
\[
\begin{cases}
\ \bigcup \ S_i \ = \ \lim (S_1, \cup \ldots \cup S_i)
\\[11pt]
\ \bigcap \ S_i \ = \ \lim (S_1, \cap \ldots \cap S_i)
\end{cases}
.]
\]
\\[-1cm]

(5)  \quad
Let $\{S_i\}$ be a sequence of sets $-$then 
$\lim (S_1 \ \Delta \ \ldots \ \Delta \ S_i)$ exists iff $\lim \ S_i = \emptyset$.
\\[-.25cm]

(6) \quad 
If $\{S_i\}$ is a sequence of sets, then
\[
X - \varlimsup S_i 
\ = \ 
\varliminf (X - S_i), 
\quad 
 X - \varliminf S_i 
\ = \ 
\varlimsup (X - S_i).
\]
\\[-1cm]

(7) \quad 
If $\{S_i\}$ and $\{T_i\}$ are sequences of sets, then
\[
\begin{cases}
\ \varlimsup \ (S_i \cup T_i) \ = \ \varlimsup \ S_i \cup \varlimsup \ T_i
\\[11pt]
\ \varliminf \ (S_i \cup T_i) \ \supset \ \varliminf \ S_i \cup \varliminf \ T_i
\end{cases}
,
\]

\[
\begin{cases}
\ \varlimsup \ (S_i \cap T_i) \subsetx \varlimsup \ S_i \cap \varlimsup \ T_i
\\[11pt]
\ \varliminf \ (S_i \cap T_i) \ = \ \varliminf \ S_i \cap \varliminf \ T_i
\end{cases}
,
\]

\[
\begin{cases}
\ \varlimsup \ (S_i - T_i) \subsetx \varlimsup \ S_i - \varlimsup \ T_i
\\[11pt]
\ \varliminf \ (S_i - T_i) \ \supset \ \varliminf \ S_i - \varliminf \ T_i
\end{cases}
.
\]
Consequently, if $\lim S_i$ and $\lim T_i$ exist, then so do $\lim (S_i \cup T_i)$, $\lim (S_i \cap T_i)$, and $\lim (S_i - T_i)$ with 

\[
\begin{cases}
\ \lim \ (S_i \cup T_i) \ = \ \lim \ S_i \cup \lim \ T_i
\\[11pt]
\ \lim \ (S_i \cap T_i) \ = \  \lim \ S_i \cap \lim \ T_i
\\[11pt]
\ \lim \ (S_i - T_i) \ = \  \lim \ S_i - \lim \ T_i
\end{cases}
.
\]
\\[-1cm]

(8) \quad 
Let $\chisubQ$ be the characteristic function of the rationals $-$then $\chisubQ$ is the pointwise limit of no sequence of continuous 
real valued functions on $\R$.  
However, $\chisubQ$ is a Baire function on $\R$ since

\[
\chisubQ (x) 
\ = \ 
\lim\limits_{m \ra \infty} \ 
[
\lim\limits_{n \ra \infty}
\{\cos^2 (m! \hsy \pi \hsy x)\}^{2 n} 
]
\qquad (x \in \R).
\]
In addition, 
\[
1 - \chisubQ (x) 
\ = \ 
\lim\limits_{m \ra \infty} \ \sgn \{\sin^2 (m! \hsy \pi \hsy x)\} 
\qquad (x \in \R).
\]
\\[-1.25cm]

[Note: \ 
This example shows that sequences do not suffice to describe a closure in the sequential modification of a space.]
\\[-.25cm]

(9) \quad 
Let $(\XX, \jT)$ be a topological space, $(\XX, \jT_\ts)$ its sequential modification.  
Let $\YY$ be a subset of $\XX$; 
let $\jT(\YY)$ and $\jT_\ts(\YY)$ be the corresponding relative topologies $-$then 
$\jT_\ts(\YY) \supset \jT(\YY)$, i.e., the sequential modification of the relative topology on $\YY$ is finer than the 
relativization to $\YY$ of the sequential modification of the topology on $\XX$, there being strict containment in general, 
but equality if $\YY$ is in addtion $\jT_\ts$-closed.
\\[-.5cm]

[To illustrate this phenomenon, take for $\XX$ the following subset of the upper half-plane + the origin:

\[
\{a_{mn} \ = \ \big(\frac{1}{m}, \frac{1}{n}\big) : m \hsy, n = 1, 2, \ldots \} 
\cup 
\{b_n \ = \ \big(0, \frac{1}{n}\big) : n = 1, 2, \ldots \} 
\cup 
\{c = (0, 0)\}.
\]

Topologize $\XX$ by specifying local open neighborhoods: 
The open neighborhoods of $a_{mn}$ and $b_n$ are to be the relativized usual open neighborhoods but the open neighborhoods of $\fc$ 
are to be the relativization of the usual open neighborhoods of $\{0\} \times ]0, \varepsilon[$ $(\varepsilon > 0)$ with $\fc$ added in.  
Consider $Y = \{a_{m n}\} \cup \{\fc\}$.]
\\[-.25cm]

(10) \quad 
Let $(\XX^\prime, \jT^\prime)$, $(\XX^{\prime\prime}, \jT^{\prime\prime})$, be topological spaces $-$then
\[
(\XX^\prime \times \XX^{\prime\prime}, \jT^\prime \times  \jT^{\prime\prime})_\ts) 
\ = \ 
(\XX^\prime \times \XX^{\prime\prime}, (\jT_\ts^\prime \times \jT_\ts^{\prime\prime})_\ts) 
\]

[To illustrate this phenomenon, take $\XX^\prime = \Q$ in the relative topology, $\jT^\prime$ and take 
$\XX^{\prime\prime} = \Q$ in the relative topology, $\jT^{\prime\prime}$ obtained by specifying that the open neighborhoods at the 
nonzero poitns are to be the relativized usual open neighborhoods but the open neighborhoods at zero itself are to be the 
relativization of the usual open neighborhoods of the sequence 
$\ds \big\{\frac{\sqrt{2}}{n} \hsy : \hsy n \in \N\big\}$ with 0 added in.  
Consider the diagonal $D$, as well as $D - \{(0,0)\}$.]
\\[-.25cm]

(11) \quad 
Suppose that $\XX$ is finite or countable $-$then the sequential modification $\jP(\XX)_\ts$ of $\jP(\XX)$ leaves $\jP(\XX)$ unchanged.
\\[-.5cm]

[Observe that if $\XX$ is finite or countable, then the topology of pointwise convergence on $\jP(\XX)$ is metrizable.]
\\[-.25cm]

(12) \quad 
Suppose that $\XX$ is uncountable $-$then the sequential modification $\jP(\XX)_\ts$ of $\jP(\XX)$ is never the same as $\jP(\XX)$.
\\[-.5cm]

[In the topology of pointwise convergence, $\jP(\XX)$ is, of course, compact.  
Show, therefore, that the uncountability of $\XX$ necessarily forces $\jP(\XX)_\ts$ to be noncompact.]
\\[-.25cm]

(13) \quad 
Let $\jS$ be a nonempty class of subsets of $\XX$ $-$then 
$\jS_{\ts\td} = \jS_{\td\ts}$ but, in general, 
$\jS_{\sigma\delta} \neq \jS_{\delta\sigma}$.
\\[-.25cm]

[Note: \ 
The second point can be seen by taking for $\jS$ the class of all bounded closed intervals of the line which have 
positive length $-$then, by a category argument, $\Q \in \jS_{\delta\sigma} - \jS_{\sigma\delta}$.]
\\[-.25cm]

(14) \quad 
There exist classes $\jS$ such that 

\[
\jS \neq \jS_\ts = \jS_\td, 
\quad 
\begin{cases}
\ \jS \neq \jS_\ts = \jS_{\ts \td}
\\[11pt]
\ \jS \neq \jS_\td = \jS_{\td \ts}
\end{cases}
.
\]
There exist classes $\jS$ such that 
\[
\begin{cases}
\ \jS \neq \jS_\sigma = \jS_{\sigma \delta}
\\[11pt]
\ \jS \neq \jS_\sigma \neq \jS_{\sigma \delta} = \jS_{\sigma \delta \sigma}
\end{cases}
.
\]
Admitting the continuum hypothesis, there exists a class $\jS$ of subsets of the line such that 
\[
\jS 
\ \neq \ 
\jS_\sigma
\ \neq \ 
\jS_{\sigma \delta}
\ \neq \ 
\jS_{\sigma \delta \sigma}
\ = \ 
\jS_{\sigma \delta \sigma\delta}.
\]
One can go much further (to any $\alpha < \Omega$!); cf. \S6.
\\[-.5cm]

[Note: \ 
The last assertion is tied up with an old problem of 
A. Kolmogoroff\footnote[4]{\vspace{.11 cm}\textit{Fund. Math.},\textbf{25} (1935), 578.}
For the details on the line, see
W. Sierpi\'nski\footnote[4]{\vspace{.11 cm}\textit{Mat. Sb.}, N. S. \textbf{43} (1936), 303-306.}
.]
\\[-.25cm]

(15)  \quad 
There exist classes $\jS$ for which $\jS$, $\jS_\trx$, $\jS_{\trx \trx}$, $\jS_{\trx \trx \trx}$, $\ldots$ are all distinct.  
If $\jS = \jS_\trx$, then $\jS = \jS_\td$ but, in general, if $\jS = \jS_\tr$, then $\jS \neq \jS_\ts$.
\\[-.25cm]

(16) \quad 
If $\jS$ is a lattice, then $\jS_{\trx \trx}$ is the class consisting of all unions of two sets from the class $\jS_\trx$.  
\\[-.5cm]

[Use the identities

\[
\begin{cases}
\ (S_1 - S_2) - (S_3 - S_4) \ = \ [S_1 - (S_2 \cup S_3)] \cup [(S_1 \cap S_4) - S_2]
\\[11pt]
\ (S_1 - S_2) - (S_3 - S_4) 
\ = \ 
[(S_21 \cup S_3) - (S_2 \cap S_4)]
- 
[(S_2 \cup S_4) - ((S_1 \cap S_4) \cup (S_2 \cap S_3))].]
\end{cases}
.
\]


\chapter{
$\boldsymbol{\S}$\textbf{6}.\quad  Problems}
\setlength\parindent{2em}
\setcounter{theoremn}{0}
\renewcommand{\thepage}{\S1-\text{P-}\arabic{page}}



\textbf{I} \ LIMITS OF LATTICES
\\[-.25cm]

Let $\jS$ be a lattice in $\XX$; 
let $\varlimsup \ \jS$ $(\varliminf \ \jS)$ stand for the subsets of $\XX$ which are the upper limit (lower limit) of a 
sequence of sets from $\jS$ $-$then
\[
\tu \hsx  \jS 
\ = \ 
\varlimsup \ \jS \cap \varliminf \ \jS.
\]
\\[-1cm]

[It suffices to prove that $\tu \jS = \jS_{\sigma \delta} \hsx \cap \hsx \jS_{\delta \sigma}$.  
For this purpose, establish the following generality.  
Let $\{S_{i, j}^\prime\}$, $\{S_{i, j}^{\prime\prime}\}$ be two double sequences of sets in $\XX$ such that 
\[
S_{i, j}^\prime \supset S_{i, j + 1}^\prime
\quad 
S_{i, j}^{\prime\prime} \subset S_{i, j+1}^{\prime\prime}
\]
with 
\[
\begin{tikzcd}
\bigcup\limits_i \bigcap\limits_j S_{i, j}^\prime
\arrow[rr,shift right=0.45,dash] \arrow[rr,shift right=-0.45,dash] 
\arrow[dr,shift right=0.5,dash] \arrow[dr,shift right=-0.5,dash]
&&
\bigcap\limits_i \bigcup\limits_j S_{i, j}^{\prime\prime}
\arrow[dl,shift right=0.5,dash] \arrow[dl,shift right=-0.5,dash]
\\
&
S
\end{tikzcd}
.
\]
Then
\[
S 
\ = \ 
\lim(
(S_{1, j}^\prime \cap S_{1, j}^{\prime\prime})
\cup 
(S_{2, j}^\prime \cap S_{1, j}^{\prime\prime} \cap S_{2, j}^{\prime\prime})
\cup \ldots \cup
(S_{j, j}^\prime \cap S_{1, j}^{\prime\prime} \cap \ldots \cap S_{j, j}^{\prime\prime}
)).]
\]

\unz{Ref} \
W. Sierpi\'nski\footnote[2]{\vspace{.11 cm}\textit{C. R. Acad. Sci. Paris}, \textbf{192} (1931), 1625-1627.}.

\newpage
\noindent
\textbf{II} \ A THEOREM OF INSERTION
\\[-.25cm]

Let $\jS$ be a lattice in $\XX$; 
let $S_\sigma \in \jS_\sigma$, $S_\delta \in \jS_\delta$ with $S_\sigma \supset S_\delta$ $-$then 
there exists an $S \in \jS_\sigma \cap \jS_\delta$ such that $S_\sigma \supset S \supset S_\delta$.
\\[-.25cm]


[Use the following generality.  
Let $\{S_i^\prime\}$, $\{S_i^{\prime\prime}\}$ be two sequences of sets in $\XX$ such that 
\[
S_i^\prime \supset S_{i +1}^\prime
\quad 
S_i^{\prime\prime} \subset S_{i+1}^{\prime\prime}
\]
with
\[
\bigcup \ S_i^\prime
\supset
\bigcap \ S_i^{\prime\prime}.
\]
Then
\[
\bigcup \ 
(S_i^\prime \cap  S_i^{\prime\prime})
\ = \ 
S_1^{\prime\prime}
\cap
(S_1^\prime\cup S_2^{\prime\prime}) 
\cap
(S_2^\prime\cup S_3^{\prime\prime}) 
\cap \ldots \hsx .]
\]

\unz{Ref} \
W. Sierpi\'nski\footnote[2]{\vspace{.11 cm}\textit{Fund. Math.}, \textbf{6} (1924), 1-5.}.

\newpage
\noindent
\textbf{III} \ UPPER LIMIT OF A SEQUENCE OF SETS
\\[-.25cm]

Let $\sI$ be the class of all infinite subsets of $\N$ $-$then, given any sequence $\{S_i\}$ of subsets of $\XX$, 
\[
\varlimsup \ S_i 
\ = \ 
\bigcup\limits_{I \in \sI} \ 
\bigcap\limits_{i \in I} \ 
S_i.
\]
Supposing that $I \in \sI$, say $I = \{i_j \hsy : \hsy j = 1, 2, \ldots\}$, let us agree to write 
$\varlimsup\limits_I \ S_i$ for $\varlimsup \ S{i_j}$.
\\[-.5cm]

It is easy to give examples where $\card(\bigcap\limits_{i \in I} \ S_i) \leq 1$ $\forall \ I \in \sI$ and yet, 
e.g., $\card(\varlimsup S_i) = \fc$.  
Accordingly, one asks instead: 
How does the cardinality of $\bigcap\limits_{i \in I} \ S_i$ influence the cardinality of $\varlimsup\limits_I S_i$?
\\[-.5cm]

(1) \quad True or False?
\\[-.5cm]

\qquad (a) \quad 
$\exists$ a sequence $\{S_i\}$ such that $\forall \ I \in \sI$,
$\bigcap\limits_{i \in I} \ S_i$ is finite but
$\varlimsup\limits_I \  S_i$ is infinite.
\\[-.25cm]

\qquad (b) \quad 
$\exists$ a sequence $\{S_i\}$ such that $\forall \ I \in \sI$,
$\card(\bigcap\limits_{i \in I} \ S_i)$ $\leq$ $\aleph_0$ 
but 
$\card(\varlimsup\limits_I \  S_i)$ $\geq \aleph_1$.
\\[-.25cm]

(2) True or False?
\\[-.5cm]

\qquad (a) \quad 
If $\forall \ I \in \sI$, $\card(\bigcap\limits_{i \in I} \ S_i) \leq N$, 
then $\exists$ an $I_0 \in \sI$ such that 
$\card(\varlimsup\limits_{I_0} \  S_i) \leq N$ $(N < \aleph_0)$.
\\[-.25cm]

\qquad (b) \quad 
If $\forall \ I \in \sI$, $\card(\bigcap\limits_{i \in I} \ S_i) \leq \aleph_0$, 
then there exists an $I_0 \in \sI$ such that 
$\card(\varlimsup\limits_{I_0} \  S_i) \leq \aleph_0$.
\\

In conclusion, let $\{S_i\}$ be a sequence of subsets of $\XX$ such that 
$\card(\bigcap\limits_{i \in I} \ S_i) \leq \aleph_0$ $\forall \ I \in \sI$ $-$then $\{S_i\}$ admits a convergent subsequence 
iff $\exists$ an $I_0 \in \sI$ such that 
$\card(\varlimsup\limits_{I_0} \  S_i) \leq \aleph_0$.
\\[-.25cm]

\unz{Ref} \
M. Laczkovich\footnote[2]{\vspace{.11 cm}\textit{Ann. of Math.},\textbf{3} (1977), 199-206.}.


\newpage
\noindent
\textbf{IV} \ THE CHARACTERISTIC FUNCTION OF A SEQUENCE OF SETS
\\[-.25cm]

Denote by $\Seq (\jP(\XX))$ the class of all sequences of subsets of $\XX$ $-$then 
by the \unz{characteristic function} of an element $\sS = \{S_i\}$ of $\Seq (\jP(\XX))$ we shall understand the function 
$\chisubsS : \XX \ra \R$ defined by the series
\[
\chisubsS (x) 
\ = \ 
2 \cdot \sum\limits_{i = 1}^\infty \ \frac{\chisubSi (x)}{3^i} 
\qquad (x \in \XX).
\]
The range of $\chisubsS$ is evidently a subset of C, the classical Cantor set.  
In fact, the map $\sS \ra \chisubsS$ implements an identification between $\Seq (\jP(\XX))$ and $\Fnc(\XX, \tC)$.
\\[-.25cm]

Here are some elementary remarks.
\\[-.25cm]

\qquad (1) \quad 
The sets in the sequence $\sS$ are all one and the same iff $\chisubsS$ assumes only the values 0 and 1.
\\[-.25cm]

\qquad (2) \quad 
The sets in the sequence $\sS$ are pairwise disjoint iff $\chisubsS$ assumes only the value 0 and values of the form 
$\ds\frac{2}{3^n}$.
\\[-.25cm]

\qquad (3) \quad 
A sequence $\sS$ of sets is increasing iff $\chisubsS$ assumes only the values 0, 1, and values of the form $\ds\frac{1}{3^n}$.
\\[-.25cm]

\qquad (4) \quad 
A sequence $\sS$ of sets is decreasing iff $\chisubsS$ assumes only the values 0, 1, and values of the form 
$\ds 1 - \Big(\frac{1}{3^n}\Big)$.
\\[-.25cm]

\qquad (5) \quad 
A sequence $\sS$ of sets is convergent iff $\chisubsS$ assumes only the values 0, 1, 
and values of the form $\ds\frac{m}{3^n}$.
\\

Characterize those $\sS$ for which: 
\\[-.25cm]

\qquad (a) \quad 
$\ov{\chisubsS (X)} = \tC$; 
\quad 
(b) \quad 
$\chisubsS (X)  = \tC$.
\\[-.25cm]

Show that if $\XX$ is in addition a topological space, then $\chisubsS$ is continuous iff all the sets in $\sS$ are open and closed. 
\\[-.25cm]

[Note: \ 
Suppose that $\XX$ is a metric space with weight $\aleph_0$ $-$then, upon consideration of the characteristic function 
of a base of cardinality $\aleph_0$, one can readily establish the following well-known results: 
\\[-.25cm]

$\XX$ is the continuous image of a subset of C.  
Furthermore, if
\[
\begin{cases}
\ \XX \ \text{is compact}
\\[11pt]
\ \XX \ \text{is compact and totally disconnected}
\\[11pt]
\ \XX \ \text{is compact, totally disconnected, and perfect}
\end{cases}
,
\]
then
\[
\begin{cases}
\ \XX \ \text{is a continuous image of C}
\\[11pt]
\ \XX \ \text{is a homeomorphic image of a closed subset of C}
\\[11pt]
\ \XX \ \text{is a homeomorphic image of C}
\end{cases}
.]
\]
\\[-1cm]

\unz{Ref} \
E. Szpilrajn\footnote[2]{\vspace{.11 cm}\textit{Fund. Math.},\textbf{31} (1938), 207-223.}
\\[-.25cm]

[A transfinite generalization can be found in 
M. Stone\footnote[3]{\vspace{.11 cm}\textit{Fund. Math.},\textbf{33} (1945), 27-33.}
.]

\newpage
\noindent
\textbf{V} \ THE EQUALITY (INEQUALITY) OF $\jS_{\sigma \hsy \delta}$ AND  $\jS_{\delta \hsy \sigma}$
\\[-.25cm]

Let $\XX$ be a set of cardinality $\aleph_0$; 
let $\jS \subset \jP(\XX)$ be nonempty $-$then, of necessity, $\jS_{\sigma \hsy \delta} = \jS_{\delta \hsy \sigma}$.
\\[-.5cm]

[This is easy, the point being that the complement of a countable subset of $\XX$ is either countable or finite.]
\\[-.5cm]

Let $\XX$ be a set of cardinality $> \aleph_0$ $-$then there necessarily exists an $\jS \subset \jP(\XX)$ for which 
$\jS_{\sigma \hsy \delta} \neq \jS_{\delta \hsy \sigma}$.
\\[-.5cm]

[There is no loss of generality in supposing that $\XX$ is a subset of $\R$ of cardinality $\aleph_1$.  
Let $\jS$ be the class of all sets of the form $X \cap I_{k, n}$, 
$\ds I_{k, n} = \Big] \frac{k}{2^n}, \frac{k+1}{2^n}  \Big[$ 
a generic dyadic open interval. 
We claim that 
$\jS_{\sigma \hsy \delta} \neq \jS_{\delta \hsy \sigma}$.
To prove this, select in each nonempty $X \cap I_{k, n}$ some point $x_{k , n}$, say. 
Denote by $X_0$ the totality of all such $-$then 
$X - X_0 \in \jS_{\sigma \hsy \delta}$ but $X - X_0 \notin \jS_{\delta \hsy \sigma}$.
\\[-.25cm]

\unz{Ref} \
W. Sierpi\'nski\footnote[2]{\vspace{.11 cm}\textit{Spis. Bulgar. Akad. Nauk}, \textbf{53} (1936), 181-195.}
\\[-.25cm]


[Note: \ 
Let $\XX$ be any set; let $\jS$ be a nonempty subset of $\jP(\XX)$.  
Write $\jS_\Sigma$, $\jS_\Delta$ for the class of subsets of $\XX$ comprised of all nonempty unions, 
nonempty intersections of sets in $\jS$ (repetitions being permissible) $-$then always 
$\jS_{\Sigma \Delta} = \jS_{\Delta \Sigma}$.]

\newpage
\noindent
\textbf{VI} \ DIFFERENCES, UNIONS, INTERSECTIONS
\\[-.25cm]

Let $\XX$ be a set of cardinality $\aleph_0$; 
Let $\jS \subset \jP(\XX)$ be nonempty $-$then
\[
\jS_{\trx \sigma \trx \sigma \trx} 
\ = \ 
\jS_{\trx \sigma \trx \sigma}, 
\quad 
\jS_{\trx \delta \trx \delta \trx} 
\ = \ 
\jS_{\trx \delta \trx \delta} 
\]
but, in general, 
\[
\jS_{\trx \sigma \trx \sigma} 
\ \neq \ 
\jS_{\trx \sigma \trx}, 
\quad 
\jS_{\trx \delta \trx \delta} 
\ \neq \ 
\jS_{\trx \sigma \trx}.
\] 
Discuss the effects of permuting the roles of r and $\sigma$ or of r and $\delta$.
\\[-.5cm]

What happens if $\XX$ has cardinality $> \aleph_0$?
\\[-.5cm]

\unz{Ref} \ 
S. Picard\footnote[2]{\vspace{.11 cm}\textit{Fund. Math.},\textbf{26} (1936), 262-266.}
\\[-.5cm]

[See also the paper of Sierpi\'nski's referred to in Prob. V.]
\\

\newpage
\noindent
\textbf{VII} \ FILTERS AND ULTRAFILTERS
\\[-.25cm]

Let $\jS$ be a nonempty subset of $\jP(\XX)$ $-$then $\jS$ is said to be a \unz{filter} on $\XX$ if:
\\[-.25cm]

\qquad (i) \quad 
$\emptyset \notin \jS$; 
\\[-.25cm]

\qquad (ii) \quad 
$\jS = \jS_\td$; 
\\[-.25cm]

\qquad (iii) \quad 
$S \in \jS, \hsx S \subset T \implies T \in \jS$.
\\[-.25cm]

The collection $\jFil(\XX)$ of all filters on $\XX$ is ordered by the inclusion relation 
(induced from that on 
$\jP(\jP(\XX)$).
\\[-.25cm]

[Note: \ 
Occasionally, condition (i) is dropped, $\jP(\XX)$ itself being regarded as a filter (cf., e.g., Exer. 9) \S4)).]
\\[-.25cm]


An \unz{ultrafilter} on $\XX$ is a filter which is properly contained in no other filter on $\XX$.  
I.e.: The ultrafilters on $\XX$ are the maximal elements in the ordered set $\jFil(\XX)$.
\\[-.5cm]

A filter $\jS$ on $\XX$ is an ultrafilter iff for each $S \subset \XX$, either $S \in \jS$ or $\XX - S \in \jS$.  
If $S_1, \ldots, S_m$ are subsets of $\XX$ whose union 
$\ds \bigcup\limits_{i = 1}^m \ S_i$ is an element of an ultrafilter $\jS$ on $\XX$, then at least one of the $S_i$ 
belongs to $\jS$.
\\[-.5cm]

Every filter is contained in an ultrafilter; 
moreover, every filter is the intersection of the ultrafilters containing it.  
\\[-.5cm]

A \unz{filterbase} on $\XX$ is a class $\jS$ of nonempty subsets of $\XX$ with the property that
\[
\forall \ S_1, \hsx S_2 \in \jS, \ \exists \ S_3 \in \jS \text{ st } S_3 \subset S_1 \cap S_2.
\]

A class $\jS$ is contained in a filter iff it is a filterbase.  
If $\jS$ is a filterbase, then
\[
\jFil(\jS) 
\ = \ 
\{T \subset \XX \hsy : \hsy \exists  \ S \in \jS \text{ st } S \subset T\}
\]
is the smallest filter containing $\jS$ or still, the filter generated by $\jS$.
\\[-.5cm]

A class $\jS$ is said to have the \unz{finite intersection property} if the intersection of the members of any finite 
subclass of $\jS$ is nonempty.  
Suppose that $\jS$ has the finite intersection property $-$then $\jS_\td$ is a filterbase, thus $\jS$ is contained in $\jFil(\jS_\td)$, 
thence in an ultrafilter.  
Every filter has the finite intersection property.
\\[-.25cm]

\unz{Ref} \ 
H. Cartan\footnote[2]{\vspace{.11 cm}\textit{C. R. Acad. Sci. Paris},\textbf{205} (1937), 595-598 and 777-779.}
\\[-.5cm]

[Note: \ 
The purpose of this problem is merely to fix the terminology and recall some basic facts.]


\newpage
\noindent
\textbf{VIII} \ COMPACT AND COUNTABLY COMPACT CLASSES
\\[-.25cm]

Let $\jS$ be a nonempty subset of $\jP(\XX)$ $-$then $\jS$ is said to be \unz{compact} (\unz{countably} \un{compact}) 
if every subclass (countable subclass) of $\jS$ with the finite
intersection property has a nonempty intersection.
\\[-.5cm]

There is no a priori connection between the compactness (countable compactness) of a class and the topology 
of pointwise convergence on $\jP(\XX)$ (or of its sequential modification).
\\[-.5cm]

\textbf{Example} \ 
Let $\XX$ be a compact (countably compact) Hausdorff space $-$then the class of all closed subsets of $\XX$ is 
compact (countably compact).
\\[-.5cm]

There are countably compact classes which are not compact. 
\\[-.5cm]

The main stabilization result is this: 
Suppose that 
\[
\jS \ \text{is} \ 
\begin{cases}
\ \text{compact}
\\[11pt]
\ \text{countably compact}
\end{cases}
.
\]
Then
\[
\begin{cases}
\ \jS_{\ts \Delta} \ \text{is compact}
\\[11pt]
\ \jS_{ts \delta} \ \text{is countably compact}
\end{cases}
.
\]

[Since compactness (countable compactness) is evidently preserved by operation $\Delta (\delta)$, 
it suffices in either case to deal with just $\jS_\ts$.  
Consider therefore, a class (countable class) $\{S_i : i \in \sI\}$ of elements of $\jS_\ts$ with the finite intersection property.  
Fix an ultrafilter $\jS_0$ on $\XX$ such that the $S_i \in \jS_0$ $\forall \ i$.  
Write $S_i = \bigcup\limits_i \ S_{i \hsy j}$ $(j \in J_i)$, $J_i$ a finite set $(S_{i \hsy j} \in \jS \ \forall \ j)$.  
Choose, as is possible (cf. Prob. VII), an index $j_i \in J_i$ for which $(S_{i \hsy j} \in \jS_0$ $-$then 
the class consisting of the $S_{i \hsy j_i}$ $(i \in \sI)$ has the finite intersection property, so 
$\bigcap \ S_{i \hsy j_i} \neq \emptyset \implies \bigcap \ S_i \neq \emptyset$.]
\\[-.25cm]

\unz{Ref} \ 
E. Marczewski\footnote[2]{\vspace{.11 cm}\textit{Fund. Math.},\textbf{40} (1953), 113-124.}
[The notion of a countably compact class of sets is due to 
W. Sierpin\'ski\footnote[2]{\vspace{.11 cm}\textit{Fund. Math.},\textbf{21} (1933), 250-275.}


\chapter{
$\boldsymbol{\S}$\textbf{2}.\quad  Partitions}
\setlength\parindent{2em}
\setcounter{theoremn}{0}
\renewcommand{\thepage}{\S2-\arabic{page}}



\qquad 
Let $\XX$ be a nonempry set $-$then by a \unz{partition} of $\XX$ we understand a class 
$P(\XX) = \{X_i \hsy : \hsy i \in I\}$ of nonempty, pairwise disjoint subsets $X_i$ of $\XX$ such that 
$\XX = \bigcup \ X_i$, i.e., such that $\chisubXX = \sum \ \chisubXsubi$.  
The $X_i$ are called the \un{components} of $P(\XX)$.  
Associated with the partition $P(\XX)$ is a surjective map $f : \XX \ra I$, viz. the map assigning to 
$x \in \X_i \hookrightarrow \XX$ the index $i$; 
conversely, associated with a surjective map $f : \XX \ra I$ is  a partition $P(\XX)$, viz. the partition whose $i^\nth$-component 
$X_i$ is the fiber $f^{-1} (i)$.  
\\

{\small\bf Example} \ 
The equivalence classes determined by an equivalence relation on $\XX$ form a partition of $\XX$.  
\\

A partition $Q(\XX)$ is said to be a \un{refinement} of the partition $P(\XX)$, written
\[
Q(\XX) \ \succeq \ P(\XX) 
\quad \text{or} \quad 
P(\XX) \  \preceq\ Q(\XX),
\]
if every component of $Q(\XX)$ is contained in a component of $P(\XX)$.  
Evidently, $Q(\XX)$ refines $P(\XX)$ iff every component of $P(\XX)$ is a union of components of $Q(\XX)$.  
The partition whose components are the elements of $\XX$ refines every partition of $\XX$; 
every partition
of $\XX$ refines the partition whose sole component is $\XX$ itself.  
\\[-.25cm]

Let

\[
\begin{cases}
\ P^\prime(\XX) \ = \ \{X_{i^\prime} \hsy : \hsy i^\prime \in I^\prime\}
\\[8pt]
\ P^{\prime\prime}(\XX) \ = \ \{X_{i^{\prime\prime}} \hsy : \hsy i^{\prime\prime} \in I^{\prime\prime} \}
\end {cases}
\]

\noindent
be two partitions of $\XX$ $-$then by the \un{meet} of $P^\prime(\XX)$ and $P^{\prime\prime}(\XX)$ we mean that partition 
$P^\prime(\XX) \wedge P^{\prime\prime}(\XX)$ of $\XX$ 
whose components are the nonempty $X_{i^\prime} \cap X_{i^{\prime\prime}}$.  
It is clear that 
$P^\prime(\XX) \wedge P^{\prime\prime}(\XX)$
is a simultaneous refinement of both $P^\prime(\XX)$ and $P^{\prime\prime}(\XX)$; 
moreover, 
$P^\prime(\XX) \wedge P^{\prime\prime}(\XX)$
is refined by every partition with this property.  
Since the relation of refinement is reflexive and transitive, it follows that the collection of all partitions of $\XX$ is in fact a 
directed set. 
\\[-.25cm]

[Note: \ 
The collection of all partitions of $\XX$ carries the structure of a lattice possessing certain supplementary characteristics 
(cf. Exer. 3).]
\\

{\small\bf Example} \ 
Suppose that $f: \XX \ra \XX$ is a map without fixed points $-$then there exists a disjoint decomposition 
\[
\XX
\ = \ 
X_1 \cup X_2 \cup X_3
\]
of $\XX$ such that $X_i \cap f(X_i) = \emptyset$ $(i = 1, 2, 3)$.  
\\[-.25cm]

[Note: \ 
Strictly speaking, this decomposition need not be a partition of $\XX$ since a given $X_i$ may be empty.  
For the easy details, see 
M. Kat\u etov\footnote[3]{\vspace{.11 cm}\textit{Comment. Math. Univ. Carolin.}, \textbf{8} (1967), 431-433.}
.]
\\

In many of the applications, the emphasis is not so much on partitioning $\XX$ by certain of its subsets as it is on 
partitioning the elements of a given class of subsets of $\XX$ by elements from that class.
\\[-.25cm]

Let, then, $\jS$ be a nonempty class of subsets of $\XX$; 
it is not required but is not excluded that $\XX$ itself belongs to $\jS$.  
Let $S \in \jS$ $-$then by an \un{$\jS$-partition} of $S$, we understand a class 
$P(S) = \{S_i \hsy : \hsy i \in I\}$ of nonempty, pairwise disjoint subsets $S_i \in \jS$ such that 
$S = \bigcup \ S_i$, i.e., such that $\chisubS = \sum \ \chisubSsubi$.  
The $S_i$ are called the \un{components} of $P(S)$.
\\[-.25cm]

The collection of all $\jS$-partitions of $S$ need not be directed by the relation of refinement, 
the point being that there is no reason to expect that the meet of two $\jS$-partitions is again an $\jS$-partition.  
However, there is a simple condition on $\jS$ which will guarantee this, namely that $\jS$ be a multiplicative class.  
The multiplicativity of $\jS$, an essentially minimal requirement, also ensures that it is permissible to take the trace of an 
$\jS$-partition. 
Thus let $S, \hsy T \in \jS$ with $S \supset T \neq \emptyset$.  
Suppose that 
$P(S) = \{S_i \hsy : \hsy i \in I\}$ 
is an $\jS$-partition of $S$ $-$then by the \un{trace} of $P(S)$ on $T$ we mean that $\jS$-partition 
$\tr_T (P(S))$ of $T$ whose components are tne nonempty $S_i \cap T$.  
To within the empty set, this notation agrees with that introduced in \S1.
\\[-.25cm]

Partitions of restricted cardinality (viz. $\leq \aleph_0$) figure prominently in the theories of the integral and derivative.  
To stress this, let us agree that an $\jS$-partition of $S \in \jS$ is \un{finite} (\un{countable}) if this is so of the corresponding 
index set.  
The class of all components arising  from all possible finite (countable) $\jS$-partitions of $S$ will be denoted by 
$\Com_\jS (S)$ ($\sigma\dash\Com_\jS (S)$) while the collection of all possible finite (countable) $\jS$-partitions of $S$ 
will be denoted by $\Par_\jS (S)$ ($\sigma\dash\Par_\jS (S)$ ).  
If $\jS$ is multiplicative, then, per the relation of refinement, both $\Par_\jS (S)$ and $\sigma\dash\Par_\jS (S)$  are directed sets.  
Conventionally, $S$ admits \un{infinite $\jS$-partitions} if $\sigma\dash\Par_\jS (S)$ is nonempty; 
of course, for this to be the case, $\jS$ itself must be at least countable.
\\


{\small\bf Example} \ 
Take $\XX = [0, 1]$ $-$then the traditional notion of a partition of $\XX$ consists in the specification of points 
$0 = x_0 < x_1 < \ldots < x_{n-1} < x_n = 1$.  
Observe, however, that the intervals 
$[x_0, x_1], \ldots, [x_{n-1}, x_n]$ do not partition $\XX$.  
The way out is to use instead the intervals 
$[x_0, x_1[, [x_1, x_2[, \ldots, [x_{n-1}, x_n[$ or the intervals 
$]x_0, x_1], ]x_1, x_2], \ldots, ]x_{n-1}, x_n]$.  
Note too that while the intervals 
$]x_0, x_1[, \ldots, ]x_{n-1}, x_n[$ do not partition $\XX$, they do constitute a topological partition of $\XX$; cf. infra.  
In passing, we remark that it is easy to exhibit countable partitions of $\XX$, e.g., $\{0\}$, and the 
$\left]\frac{1}{n+1}, \frac{1}{n}\right[$ $(n = 1, 2, \ldots)$.  
Consider now the class $\jS$ of all closed subintervals of $\XX$; 
$\jS$ is multiplicative, singletons (as well as the empty set) belonging to $\jS$.  
Given $[a,b]$ in $\jS$, it is clear that
\[
\Par_\jS ([a,b]) 
\ = \ 
\{[a,b]\}, 
\quad 
\sigma\dash\Par_\jS ([a,b]) 
\ = \ 
\emptyset.
\]
Therefore, in so far as it is a question of finite or countable partitions, $\jS$ is inutile.  
Trivially, of course, $[a, b] \ = \ \bigcup\limits_{a \leq x \leq b} \ \{x\}$, 
an uncountable union (if $b > a)$.  
\\

The preceding example, its essential simplicity notwithstanding, already contains a degree of unpleasantness.  
Our strictly set theoretic definition of partition allows for no overlap in the components.  
In certain situations, however, this turns out to be an unduly restrictive condition, particularly in the presence of other structures, 
for instance, a topology.  
Though this will not be a
point of conscern at present, nevertheless an illustration may prove helpful.
\\

Let $\XX$ be a topological space $-$then by a \un{topological partition} of $\XX$ we understand a class 
$P(\XX) = \{X_i \hsy : \hsy i \in I\}$ of nonempty, pairwise disjoint, open and connected subsets $X_i$ of $\XX$ 
such that $\bigcup \ X_i$ is dense in $\XX$.  
The $X_i$ are called the \un{components} of $P(\XX)$.  
A topological partition $Q(\XX)$ is said to be a \un{refinement} of the topological partition $P(\XX)$, written 
\[
Q(\XX) \ \underset{\tT}{\succeq} \ P(\XX) 
\quad \text{or} \quad 
P(\XX) \  \underset{\tT}{\preceq} \ Q(\XX), 
\]
if every component of $Q(\XX)$ is contained in a component of $P(\XX)$.  
\\[-.25cm]

Specialize now and suppose that $\XX$ is actually a metric space with metric $\td$.  
Let $\varepsilon > 0$ $-$then an \un{$\varepsilon$-partition} of $\XX$ is a topological partition with the property that 
each of its components has diameter $< \varepsilon$.  
$\XX$ is called \un{d-partitionable} if for every $\varepsilon > 0$, there exists an $\varepsilon$-partition of $\XX$.
\\

{\small\bf Example} \ 
The metric space $(\XX, \td)$ is \un{strongly $\td$-partitionable} if for every $\varepsilon > 0$, 
there exists a finite $\varepsilon$-partition of $\XX$.  
We then ask: \ 
What metric spaces are strongly $\td$-partitionable? 
It turns out that there is a very simple answer.  
to give it, recall that $\XX$ has \un{property S} if for every $\varepsilon > 0$, $\XX$ can be written as 
the union of a finite number of connected subsets each of diameter less than $\varepsilon$.  
In terms of this notion, the sought for characterization then reads:  \
$\XX$ is strongly $\td$-partitionable iff $\XX$ has property S.  
Consequently, if $\XX$ is strongly $\td$-partitionable, then for every $\varepsilon > 0$, it is possible to find a finite 
$\varepsilon$-partition of $\XX$ all components of which have property S, 
hence there is a partition $P_1(\XX),  P_2(\XX), \ldots$ such that $P(\XX)$ is a finite $1/i$-partition of $\XX$ and 
$P_{i+1} (\XX)$ is a refinement of $P_i (\XX)$.  
Assume in addition that $\XX$ is compact and connected, i.e., that $\XX$ is a continuum $-$then, as is well known, 
$\XX$ is locally connected iff $\XX$ has property S.  
By definition, a \un{continuous curve} is a locally connected continuum.  
In view of what has been said, therefore, every continuous curve is strongly $\td$-partitionable, a theorem of R. Bing.
\\

[Note: \ 
For a complete discussion of these and related results, see 
R. Bing\footnote[2]{\vspace{.11 cm}\textit{Bull. Amer. Math. Soc.}, \textbf{55} (1949), 1101-1110.}
\textit{Bull. Amer. Math. Soc.}, \textbf{55} (1949), 1101-1110, 
and 
R. Bing\footnote[3]{\vspace{.11 cm}\textit{Bull. Amer. Math. Soc.}, \textbf{58} (1952), 536-556.}
\textit{Bull. Amer. Math. Soc.}, \textbf{58} (1952), 536-556.]
\\

\[
\textbf{\unz{Notes and Remarks}}
\]

Partitions, in one guise or another, have been around from the beginning.  
They will play a central role in the sequel.  
Incidentally, it should be noted that partitions and equivalence relations are coextensive notions, 
both being descriptions of the same mathematical reality.  
Observe too that the axiom of choice is entirely equivalent to the statement that every partition of every set has a 
\unz{set of representatives}, i.e., if 
$\sP(\XX) = \{X_i : i \in \sI\}$ 
is a partition of $\XX$, 
then there exists a subset $\tC_{\sP(\XX)}$ of $\XX$ such that 
$\tC_{\sP(\XX)} \cap X_i = \{x_i\}$ $\forall \ X_i$.  
The discovery that continuous curves could be topologically partitioned was one of the most important 
combinatorial developments of the 1950's.  
The term continuous
curve arises, of course, from the famous theorem of Hahn-Mazurkiewicz which states that a metric space 
is a continuous curve iff it is the continuous image of $[0,1]$.  
For this reason, continuous curves are sometimes referred to as 
\unz{Peano spaces}.  
A systematic treatment of these matters can be found in 
G.T. Whyburn\footnote[2]{\vspace{.11 cm}\un{Analytic Topology},  \textit{Amer. Math. Soc. Colloquium Publications}, \textbf{28}, New York,  1942} 
and 
T. Rad\'o\footnote[3]{\vspace{.11 cm}\un{Length and Area}, 
\textit{Amer. Math. Soc. Colloquium Publications}, \textbf{28}, New York,  1948} 
Finally, for much additional information on the general theory of partitions, the reader can consult with profit 
O. Ore\footnote[4]{\vspace{.11 cm}\textit{Duke. Math. J.}, \textbf{9} (1942), 573-627.}


\chapter{
$\boldsymbol{\S}$\textbf{2}.\quad  Exercises}
\setlength\parindent{2em}
\setcounter{theoremn}{0}
\renewcommand{\thepage}{\S2-\text{E-}\arabic{page}}


\qquad
(1) \quad 
For $n = 1, 2, \ldots,$ let $p_n$ be the number of partitions of a set of $n$ elements $-$then 
the $p_n$ satisfy the recursion relation
\[
p_{n + 1} 
\ = \ 
1 + \sum\limits_{k = 1}^n \ \binom{n}{k} \hsy p_k.
\]
What is the relationship between the $p_n$ and $\exp (\exp x - 1)$?
\\

(2)   
Let $\jX = \{X_i : i \in \sI\}$ be a class of nonempty subsets of a set $\XX$ $-$then
$\jX$ determines a partition $P_\jX$ of $\XX$ which partitions each of the $X_i$ and is refined by any partition of $\XX$ with this property.

[Given a subset $E$ of $\sI$, put
\[
X_E 
\ = \ 
\bigcap\limits_{i \in E} \ X_i 
\hsx \cap \hsx 
\bigcap\limits_{i \in \hsy \sI - E} \ (\XX - X_i).
\]

Consider the nonempty $X_E$.]
\\[-.25cm]

(3) \quad 
Let
\[
\begin{cases}
\ P^\prime (\XX) = \{X_{i^\prime} : i \in \sI^\prime\} \\[4pt]
\ P^{\prime\prime} (\XX) = \{X_{i^{\prime\prime}} : i^{\prime\prime} \in \sI^{\prime\prime}\}
\end{cases}
\]
be two partitions of $\XX$ $-$then by the \unz{join} of 
$P^\prime (\XX)$
and 
$P^{\prime\prime} (\XX)$
we mean that partition 
$P^\prime (\XX) \vee P^{\prime\prime} (\XX)$
of $\XX$ whose components are the minimal nonempty 
$\bigcup \ X_{i^\prime} = \bigcup \ X_{i^{\prime\prime}}$.  
It is clear that 
$P^\prime (\XX) \vee P^{\prime\prime} (\XX)$
is refined simultaneously by both 
$P^\prime (\XX)$
and 
$P^{\prime\prime} (\XX)$; 
moreover 
$P^\prime (\XX) \vee P^{\prime\prime} (\XX)$
refines every partition with this property.
\\[-.5cm]

[Note: \ 
In the technical language of the trade, the collection of all partitions of $\XX$ is a relatively complemented, semimodular, complete lattice 
with largest and smallest elements.  
It is called the \unz{partition lattice} attached to $\XX$.  
Up to isomorphism, every abstract lattice appears as a sublattice of some such partition lattice; 
cf. 
P. Whitman\footnote[2]{\vspace{.11 cm}\textit{Bull. Amer. Math. Soc.}, \textbf{52} (1946), 507-522.}. 
\\[-.25cm]

(4) \quad 
Suppose that $\XX = \bigcup\limits_{i = 1}^m \ X_i$ is the union of $m = 2^m$ nonempty, distinct subsets $X_i$ $-$then 
there exist $n + 1$ nonempty pairwise disjoint subsets $Y_j$ of $\XX$ such that 
$\XX = \bigcup\limits_{j = 1}^{n + 1} \ Y_j$.
\\[-.5cm]

[There are two ways to look at this.  
The first method consists in remarking that $\XX$ must have at least $n + 1$ distinct elements, say $x_1, \ldots, x_{n + 1}$, so
\[
\XX 
\ = \ 
\{x_1\} \cup \ldots \cup \{x_n\} 
\ \bigcup \ 
(\XX - \{x_1, \ldots, x_n\})
\]
which is certainly a partition of $\XX$ with the desired property.  
However, while the axiom of choice has not been used, the construction can hardly be considered effective.  
The second (effective) method consists in considering $M = \{1, \ldots, m\}$, the $2^m - 1$ nonempty subsets of which 
$\{i_1, \ldots, i_s\}$, can be arranged into a finite sequence according to the size of the number $2^{i_1} + \cdots + 2^{i_s}$.
Denoting by $\{M_k\}$ the sequence thereby obtained, put
\[
Z_k 
\ = \ 
\bigcap\limits_{i \in M_k} \ X_i \ - \ \bigcup_{i \in M - M_k} \ X_i. 
\]
The $Z_k$ may be used to determine the $Y_j$.]

(5)  \quad
Suppose that $\XX = \bigcup\limits_{i = 1}^\infty \ X_i$ is the union of countably many nonempty, distinct subsets $X_i$ $-$then 
there exist countably many nonempty, pairwise disjoint subsets $Y_j$ of $\XX$ such that 
$\XX = \bigcup\limits_{j = 1}^\infty \ Y_j$.
\\[-.5cm]

[The axiom of choice is not needed here (Kuratowski); 
cf. 
A. Tarski\footnote[3]{\vspace{.11 cm}\textit{Fund. Math.}, \textbf{6} (1924), 94-95.}.
.]
\\[-.25cm]

(6) \quad 
Let $\XX$ be a set; let $f : \XX \ra \XX$ be a map.  
Suppose that $f$ is injective $-$then $\XX$ can be uniquely decomposed as a countable union of pairwise disjoint sets 
$X_0, X_1, \ldots$ (possibly $\emptyset$) such that 
\[
f(X_0) = X_0, 
\quad 
f(X_i) = f(X_{i + 1}) 
\quad (i \geq 1).
\]

[Take 
\[
X_0 
\ = \ 
\bigcap\limits_{i = 1}^\infty \ f^i(\XX), 
\quad 
X_i = f^{i - 1} (\XX) - f^i (\XX) 
\qquad (i \geq 1),
\]
where $f^0 (\XX) = \XX$.]
\\[-.25cm]

(7) \quad 
Let $\XX$ and $\YY$ be sets; let $f: \sP(\XX) \ra \sP(\YY)$ and $g: \sP(\YY) \ra \sP(\XX)$ be maps.  
Suppose that 
\[
\begin{cases}
\forall \ S, \hsx T \in \sP(\XX), \ S \subset T \implies f(S) \subset f(T)
\\[4pt]
\forall \ S, \hsx T \in \sP(\YY), \ S \subset T \implies g(S) \subset g(T)
\end{cases}
.
\]
Then there exist disjoint decompositions $X = X_1 \cup X_2$, $Y = Y_1 \cup Y_2$ such that $f(X_1) = Y_1$, $g(Y_2) = X_2$.   
Must these decompositions be partitions of $\XX$ or $\YY$?
\\[-.5cm]

[First prove that if $M$ is a set, $\Phi : \sP(M) \ra \sP(M)$ a map such that 
\[
\forall \ A, \hsx B \in \sP(M), \ A \subset B \implies \Phi(A) \subset \Phi(B), 
\]
then for some subset $M_0$ of $M$, $\Phi(M_0) = M_0$.  
This done, specialize and for $S \subset \XX$, put

\[
\Phi(S) 
\ = \ 
X - g(\YY - f(S)).
\]
The preceding remark implies that $\Phi$ has a fixed point $X_1$, say.  
Take, then, $X_2 = \XX - X_1$, $Y_1 = f(X_1)$, $Y_2 = \YY - Y_1$.]
\\[-.25cm]

(8) \quad 
There exists a nonempty set $\XX$ and a nonempty class $\jS$ of subsets of $\XX$ with the following property: 
Every nonempty $S \in \jS$ admits a partition by three elements of $\jS$ but no nonempty $S \in \jS$ admits a partition by two elements of $\jS$. 
Can $\jS$ be taken multiplicative?
\\[-.25cm]

(9) \quad 
Let $\jS$  be a nonempty class of subsets of $\XX$ with the property that every nonempty element of $\jS$ can be written as the union of 
three distinct elements of $\jS$ $-$then every nonempty element of $\jS$ can be written as the union of two distinct elements of $\jS$.
\\[-.25cm]

(10) \quad 
There exist a nonempty set $\XX$ and a nonempty class $\jS$ of subsets of $\XX$ with the following property: 
Every nonempty $S \in \jS$ admits a partition by two elements of $\jS$ but no nonempty $S \in \jS$ admits a partition by countably many 
elements of $\jS$.  
Can $\jS$ be taken multiplicative?
\\[-.5cm]

[Note: \ 
Suppose that $\XX = \N$ $-$then in this case, 
if every nonempty $S \in \jS$ can be partitioned by two elements of $\jS$, 
it must actually be the case that every nonempty $S \in \jS$ can be partitioned by countably many elements of $\jS$.] 
\\[-.25cm]

(11) \quad 
Exhibit an explicit countable partition of $\N$, each component of which is countable.
\\[-.25cm]

(12) \quad 
Exhibit an explicit countable partition of $\R$, each component of which consists of two elements.
\\[-.25cm]

(13) \quad 
Exhibit an explicit countable partition of $[0,1]$, each component of which consists of two elements.
\\[-.25cm]

(14) \quad
Take $\XX = \R$ $-$then there exists a subset $S$ of $\XX$ and a countable set of real numbers $\{s_i\}$ such that
\[
X 
\ = \ 
\bigcup\limits_{i = 1}^\infty \ (s_i + S),
\]
where
\[
i \neq j 
\implies 
(s_i + S) \cap (s_j + S) = \emptyset.
\]

[This is easy: 
Put $S = [0,1[$ and choose the $s_i$ in the obvious way.]
\\[-.25cm]

(15)  \quad
Take $\XX = [0,1]$ $-$then there exists a subset $S$ of $\XX$ and a countable set of real numbers $\{s_i\}$ such that 
\[
X 
\ = \ 
\bigcup\limits_{i = 1}^\infty \ (s_i + S),
\]
where
\[
i \neq j 
\implies 
(s_i + S) \cap (s_j + S) = \emptyset.
\]

[This is difficult; 
cf. 
J. v. Neumann\footnote[3]{\vspace{.11 cm}\textit{Fund. Math.}, \textbf{6} (1924), 94-95.}.
We remark that the axiom of choice is needed here; naturally, neither $S$, nor any of its translates is Lebesgue measurable.]
\\[-.25cm]

(16) \quad 
The continuum hypothesis is equivalent to the statement that the real line - the origin can be partitioned into countably many rationally independent sets.  
\\[-.5cm]

[This result is due to 
P. Erd\"os and 
S. Kakutani\footnote[2]{\vspace{.11 cm}\textit{Bull. Amer. Math. Soc.}, \textbf{49} (1943), 457-461.}. 
In brief, the argument runs as follows.
\\[-.5cm]

Admit the continuum hypothesis.  
Let $\{x_\beta : \beta < \omega_1\}$ be a Hamel basis for $\R$.  
Given nonzero rational numbers $r_1, \ldots, r_n$, write $\R(r_1, \ldots, r_n)$ for the set of all $x \in \R$ such that 
$x = r_1 x_{\beta_1} + \cdots + r_n x_{\beta_n} $ $(\beta_1 < \ldots < \beta_n)$ $-$then, in an obvious notation
\[
\R 
\ = \ 
\{0\} \cup 
\bigcup _{(r_1, \ldots, r_n)} \ \R(r_1, \ldots, r_n)
\qquad \text{(disjoint union).}
\]

Decompose each $\R(r_1, \ldots, r_n)$ by considering $\beta < \omega_1$ the subset comprised of those $x$ for which $\beta_n = \beta$.
\\[-.5cm]

Deny the continuum hypothesis.  
Let $\{x_\beta : \beta < \omega_\alpha\}$ be a Hamel basis for $\R$ $-$then $\alpha \geq 2$. 
Let $X_i$ be any countable partition of $\R - \{0\}$ $-$then there exists an index $i$ for which

\[
\card(\{\omega_1 \leq \beta < \omega_\alpha : i(\beta) = i\}) 
\ \geq \ \aleph_2,
\]
where $i(\beta)$ is defined by requiring that there be ordinals $\beta_{\beta}^\prime$,  $\beta_{\beta}^{\prime\prime}$ with

\[
\begin{cases}
\ \beta_\beta^\prime < \omega_1, \quad \beta_\beta^{\prime\prime} < \omega_1
\\[4pt]
\ \beta_\beta^\prime < \beta_\beta^{\prime\prime}
\end{cases}
\qquad
\begin{cases}
\ x_{\beta_\beta^\prime} + x_\beta \in X_{i(\beta)}
\\[4pt]
x_{\beta_\beta^{\prime\prime}} + x_\beta \in X_{i(\beta)}
\end{cases}
\]

\noindent
Conclude from this that there exist ordinals

\[
\begin{cases}
\ \beta^\prime < \omega_1, \quad \beta^{\prime\prime} < \omega_1
\\[4pt]
\ \beta^\prime < \beta^{\prime\prime}
\end{cases}
\qquad 
\begin{cases}
\ \omega_1 \leq \beta^\star \quad \beta^{\star\star} < \omega_\alpha
\\[4pt]
\ \beta^\star < \beta^{\star\star}
\end{cases}
\qquad 
i = i(\beta^\star) = i(\beta^{\star\star})
\]

\noindent
such that
\[
x_{\beta^\prime} + x_{\beta^\star}, 
\quad 
x_{\beta^{\prime\prime}} + x_{\beta^\star}, 
\qquad 
x_{\beta^\prime}+ x_{\beta^{\star\star}}, 
\quad 
x_{\beta^{\prime\prime}} + x_{\beta^{\star\star}}
\]
all belong to $X_i$.]


\chapter{
$\boldsymbol{\S}$\textbf{2}.\quad  Problems}
\setlength\parindent{2em}
\setcounter{theoremn}{0}
\renewcommand{\thepage}{\S2-\text{P-}\arabic{page}}

\textbf{I} \ MESH FUNCTIONS 
\\[-.25cm]

Let $\jS$ be a nonempty class of subsets of $\XX$; 
let $S \in \jS$ $-$then by a \unz{mesh function} $\delta$ on $\Par_\jS (S)$ we understand a rule which assigns to each 
$P(S)$ in 
$\Par_\jS (S)$  
a positive real number $\delta ( P(S))$ subject to the following rule: 
$\forall \ \varepsilon > 0$, $\exists \ P_\varepsilon (S) \in \Par_\jS (S)$ such that 
$\delta (P_\varepsilon (S)) < \varepsilon$.  
If $\Par_\jS (S)$  admits a mesh function $\delta$, then $\delta$ can be used to direct $\Par_\jS (S)$: 
$Q(S) \underset{\delta}{\leq} P(S)$ iff 
$\delta(Q(S)  \leq \delta(P(S))$.  
It is to be stressed that if $Q(S)$ is a refinement of $P(S)$, then there may be no relation between 
$\delta(Q(S)$ and $\delta(P(S))$; in fact, $\delta$ need not decrease upon refinement.
\\[-.25cm]

[Take $\XX = [0, 1[$ and let $\jS$ be the class of all left closed and right open subintervals $[a,b[$ of $\XX$.  
Fix $S = [a,b[$ in $\jS$ $-$then an element $P(S)$ in $\Par_\jS (S)$  has the form 
$\{[a_i, b_i[ \hsy : i = 1, \ldots, n\}$, 
where, say, $a_1 = a$, $b_n = b$ and $a_1 < b_1 = a_2 < b_2 \cdots$.  
Put $\delta (P(S)) = \max (b_i - a_i)$ $-$then $\delta$ is a mesh function on $\Par_\jS (S)$ which, 
moreover, does not decrease upon refinement.  
Define now a function $\sigma$ on $\XX$ via the following stipulation: \ 
$\sigma(x) = 0$ if $x$ is irrational, 
$\sigma(x) = 1/q$ if $x = p/q$ is rational $(0 \leq p \leq q, \hsx q \min.)$.  
Put 
$
\ds
\delta (P(S)) = 
\sum \ 
(\sigma (a_i) + \sigma (b_i)) + \max (b_i - a_i) - (\sigma (a) + \sigma (b))
$ 
$-$then 
$\delta$ is a mesh function on $\Par_\jS (S)$ which, this time, need not decrease upon refinement.]
\\[-.5cm]

\noindent
\unz{Ref} \quad
L. Cesari\footnote[4]{\vspace{.11 cm}\textit{Trans. Amer. Math. Soc.}, \textbf{102} (1962), 94-113.}. 
\\[-.25cm]

[Note: \ 
Suppose that $\XX$ is a continuous curve.  
Let TOP-Par($\XX$) be the collection of all finite $\varepsilon$-partitions of $\XX$ $-$then the rule which assigns to each 
$P(\XX)$ in TOP-Par($\XX$) the maximum diameter of its components can be viewed, 
in the obvious way, as a mesh function on TOP-Par($\XX$) which moreover, decreases upon refinement.]

\newpage
\noindent
\textbf{II} \ THEOREMS OF RAMSEY AND SIERPI\'NSKI
\\[-.25cm]

Given a set $\XX$ and a natural number $n$, let us agree to write $\langle \XX \rangle_n$ for the class of all subsets 
of $\XX$ of cardinality $n$.
\\

\unz{\textbf{Theorem}} \ (Ramsey) \ 
Let $\XX$ be a set of cardinality $\aleph_0$; 
let $\{\jX_1, \ldots, \jX_m\}$ be a finite partition of $\langle \XX \rangle_n$ $-$then 
there exists an infinite subset $S$ of $\XX$ and an index $i$ such that 
$\langle S \rangle_n \subset \jX_i$.
\\[-.5cm]

[There is no loss of generality in taking $\XX = \N$.  
This being so, it will then be enough to prove that for any map 
$f : \langle \N \rangle_n \ra \{1, \ldots, m\}$, there exists an infinite subset $S$ of $\N$ such that $f$ is constant on 
$\langle S \rangle_n$.  
If $n = 1$, the result is clear so assume that it holds for $n \geq 1$.  
Let 
$f : \langle \N \rangle_{n + 1} \ra \{1, \ldots, m\}$
be a map.  
Given $x \in \N$, write $f_x$ for the function on 
$\langle \N - \{x\} \rangle_n$
defined by the rule
\[
f_x (?) 
\ = \ 
f(\{x\} \cup ?).
\]
Apply the induction hypothesis in an appropriate way to $f_x$.]
\\[-.25cm]

\unz{Ref} \quad 
F. Ramsey\footnote[2]{\vspace{.11 cm}\textit{Proc. London Math. Soc.} (2),\textbf{30} (1930), 264-286.}.
\\[-.25cm]

One possible generalization of Ramsey's theorem might read: 
Let $\XX$ be a set of cardinality $\aleph_1$; 
let 
$\{\jX_1, \ldots, \jX_m\}$ 
be a finite partition of 
$\langle \XX \rangle_n$ 
$-$then there exists a subset $S$ of $\XX$ of cardinality $\aleph_1$ and an index $i$ such that 
$\langle S \rangle_n \subset \jX_i$.  
This statement is, however, false.  
In fact, even more can be said:
\\[-.25cm]

\unz{\textbf{Theorem}} \ (Sierpi\'nski) \ 
Let $\XX$ be a set of cardinality $2^{\aleph_0}$ $-$then 
there exists a finite partition 
$\{\jX_1, \ldots, \jX_m\}$ of $\langle \XX \rangle_n$ 
with the following property: 
For every subset $S$ of $\XX$ of cardinality $\aleph_1$, 
$\langle S \rangle_n \not\subset \jX_i$ 
$(i = 1, \ldots, m)$.  
\\[-.5cm]

[There is no loss of generality in taking $\XX = \R$.  
Furthermore, it can be supposed that $m = 2$, $n = 2$, the general case being a consequence of this one.  
Let $<$ be the usual ordering of $\R$; 
let $<_w$ be some well-ordering of $\R$ $-$then we define a map 
$f : \langle \R \rangle_2 \ra \{0, 1\}$ 
by requiring that $f(\{x, y\}) = 0$ if $<$ and $<_w$ order the pair $\{x, y\}$ in the same way 
and 
$f(\{x, y\}) = 1$ 
if $<$ and $<_w$ order the pair $\{x, y\}$ in the opposite way.  
If now $S$ were a subset of $\R$ of cardinality 
$\aleph_1$ 
such that either 
$f(\langle S \rangle_2) = 0$ 
or 
$f(\langle S \rangle_2) = 1$, 
then of necessity either the natural order or its inverse would well-order $S$, an impossibility.  
The partition of 
$\langle \R \rangle_2$ 
canonically associated with $f$ thus has the desired properties.]
\\[-.25cm]

\unz{Ref} \quad 
W. Sierpi\'nski\footnote[3]{\vspace{.11 cm}\textit{Ann. Scoula Norm. Sup. Pisa Cl. Sci.} (2), \textbf{2} (1933), 285-287.}.
\\[-.25cm]

[A useful survey on this interesting subject was given by 
P. Erd\"os and 
R. Rado\footnote[4]{\vspace{.11 cm}\textit{Bull. Amer. Math. Soc.}, \textbf{62} (1956), 427-489.}. 
See also 
P. Erd\"os,
A Hajnal, and 
R. Rado\footnote[5]{\vspace{.11 cm}\textit{Acta Math. Acad. Sci. Hungar.}, \textbf{16} (1965), 93-196.}. 
P. Erd\"os and
A Hajnal\footnote[6]{\vspace{.11 cm}\textit{Proc. Symp. Pure Math.} ,\textbf{13} (1971), 17-48.}. 
For an account of recent developments (and additional references), cf. 
R. Graham, 
B. Rothschild, and 
J. Spencer\footnote[7]{\vspace{.11 cm}\un{Ramsey Theory}, Wiley New York, 1980}.]
\\[.5cm]

\noindent
\textbf{III} \ DISJOINT AND NONDISJOINT CLASSES
\\

Let $\jS$ be an infinite class of sets $-$then there necessarily exists an infinite subclass $\jS_0$ of $\jS$ such that 
\[
\forall \ S^\prime, \hsx S^{\prime\prime} \in \jS_0 \hsx : \hsx S^\prime \neq S^{\prime\prime}
\implies 
S^\prime \cap S^{\prime\prime}
= \emptyset
\]
or
\[
\forall \ S^\prime, \hsx S^{\prime\prime} \in \jS_0 \hsx : \hsx S^\prime \neq S^{\prime\prime}
\implies 
S^\prime \cap S^{\prime\prime}
\neq \emptyset.
\]
On the other hand, there exists an uncountable class $\jS$ of sets such that $\jS$ contains no uncountable subclass 
having one or the other of the preceding properties.
\\

\unz{Ref} \quad 
W. Sierpi\'nski\footnote[2]{\vspace{.11 cm}\textit{Fund. Math.}, \textbf{35} (1948), 165-174.}.
\\

\newpage
\noindent
\textbf{IV} \ PARTITIONS OF THE PLANE
\\

The continuum hypothesis is equivalent to the statement that the plane can be partitioned into two sets $X$ and $Y$, 
where $X$ ($Y)$ intersects every line parallel to the $x$ ($y$)-axis in a finite or countable set.
\\[-.25cm]

\unz{Ref} \quad 
W. Sierpi\'nski\footnote[2]{\vspace{.11 cm}\textit{Bull. Acad. Sci Cracovie} , \textbf{A} (1919), 1-3.}, 
W. Sierpi\'nski\footnote[3]{\vspace{.11 cm}\textit{Fund. Math.}, \textbf{5} (1924), 177-187.}.
\\[-.25cm]

The plane cannot be partitioned into two sets $X$ and $Y$, where $X$ intersects every line parallel to the $x$-axis 
in a finite set and $Y$ intersects every line parallel to the $y$-axis in a finite or countable set.
\\[-.5cm]

\unz{Ref} \quad 
H. Tietze\footnote[4]{\vspace{.11 cm}\textit{Math. Ann.}, \textbf{88} (1923), 290-312.}.
\\[-.5cm]


[Note: \ By comparison, the continuum hypothesis is equivalent to the statement that 
space 
can be partitioned into three sets $X$, $Y$, and $Z$, where $X$ ($Y, \hsx Z$) intersects every line parallel to the 
$x$ $(y, \hsx z)$-axis in a finite set; 
see 
W. Sierpi\'nski\footnote[5]{\vspace{.11 cm}\textit{Rend. Circ. Mat. Palermo}, (2) \textbf{1} (1952), 7-10.}.]
\\[-.25cm]

The axiom of choice implies that the plane can be partitioned into two sets $X$ and $Y$, where $X$ ($Y$) 
intersects every line parallel to the $x$ ($y$)-axis in a set of cardinality $< 2^{\aleph_0}$.
\\[-.5cm]

\unz{Ref} \quad 
W. Sierpi\'nski\footnote[6]{\vspace{.11 cm}\textit{Soc. Sci. Lett. Varsovie C. R. Cl. III Sci. Math. Phys.}, 
\textbf{25} (1932), 9-12.}
\\[-.25cm]

The continuum hypothesis is equivalent to the statement that there exist in the plane three straight lines 
$L_1$, $L_2$, $L_3$, with the property that the plane is the union of three sets 
$S_1$, $S_2$, $S_3$ such that $S_i$ intersects every line parallel to $L_i$ $(i = 1, 2, 3)$ in a finite set.
\\[-.5cm]

\unz{Ref} \quad 
F. Bagemihl\footnote[7]{\vspace{.11 cm}\textit{Rend. Circ. Mat. Palermo},  \textbf{7} (1961), 77-79.}.
R. Davies\footnote[8]{\vspace{.11 cm}\textit{Z.Math. Logic Grundlag. Math.},  \textbf{8} (1962), 109-111.}.
\\[-.25cm]

The axiom of choice implies that the plane can be partitioned into countably many sets, 
none of which contains the vertices of an equilateral triangle.

\unz{Ref} \quad 
J. Ceder\footnote[9]{\vspace{.11 cm}\textit{Rend. Circ. Mat. Palermo},  \textbf{14} (1969), 1247-1251.}.
\\[-.25cm]

The continuum hypothesis imples that the plane can be partitioned into countably many sets, 
none of which contains the vertices of an isosceles triangle.

\unz{Ref} \quad 
R. Davies\footnote[1]{\vspace{.11 cm}\textit{Proc. Canbridge Philos. Soc.},  \textbf{72} (1972), 179-183.}.
\\[-.25cm]


[Note: \ 
There is an extensive literature on these and related themes.  
For additional results, together with a variety of conjectures, see 
P. Erdos\footnote[2]{\vspace{.11 cm}\textit{Real Anal. Exchange},  \textbf{4} (1978-79), 113-138.}
.]


\chapter{
$\boldsymbol{\S}$\textbf{3}.\quad  Semirings}
\setlength\parindent{2em}
\setcounter{theoremn}{0}
\renewcommand{\thepage}{\S3-\arabic{page}}



\qquad
Let $\XX$ be a nonempty set; 
let $\jS$ be a subset of $\jP (\XX)$ containing the empty set $-$then 
$\jS$ is said to be a 
\un{semiring} (\un{$\sigma$-semiring}) 
if $\jS$ is multiplicative and if for all nonempty $S$, $T \in \jS$, wth $S \supset T$, 
there exists a finite (finite or countable) $\jS$-partition of $S$ having $T$ has a component.  
A 
\un{semialgebra} (\un{$\sigma$-semialgebra})  is a 
semiring ($\sigma$-semiring) containing $X$.  
It is clear that every semiring is a $\sigma$-semiring but the converse is not true.  
Conventionally, $\{\emptyset\}$ is both a semiring and a $\sigma$-semiring.
\\

{\small\bf Examples}
\\[-.5cm]

\qquad (1) \quad 
Take $\XX = \R$.  
Let $\jS$ be the class consisting of all bounded open intervals, and all singletons $-$then $\jS$ is a semiring.
\\

\qquad (2) \quad 
Take $\XX = \R$.  
Let $\jS$ be the class consisting of all bounded, left closed and right open intervals and all singletons $-$then 
$\jS$ is a $\sigma$-semiring but not a semiring.
\\

Partition theory leads at once to the consideration of 
semirings ($\sigma$-semirings).  
Indeed, let $\jS$ be a multiplicative class; let $S \in \jS$ $-$then the class consisting of the empty set and the 
elements of $\Com_\jS (S)$ ($\sigma\dash\Com_\jS (S)$) is a 
semiring ($\sigma$-semiring).
\\[-.25cm]

[Note: \ 
Tacitly, of course, $S \neq \emptyset$.  
Accordingly, $S \in \Com_\jS (S)$, hence $\Com_\jS (S)$ is not empty.  
On the other hand, 
$\sigma\dash\Com_\jS (S)$ 
may very well be empty (cf. \S2).]
\\

Semirings (or $\sigma$-semirings) also arise naturally in the presence of certain chain conditions.  
Thus let $\jS$ be a multiplicative class containing the empty set $-$then $\jS$ is said to satisfy the 
\un{finite} (\un{countable}) \un{chain condition} 
if for all $S$, $T \in \jS$ with $S \supset T$, there exists a finite (countable) class $\{S_i\} \subset \jS$ such that 
\[
T 
\ = \ 
S_1 \subset S_2 \subset \ldots \subset 
\bigcup\limits_i \ S_i 
\ = \ 
S, 
\]
where $S_i - S_{i - 1} \in \jS$ for each $i > 1$.  
Here, of course, repetitions are allowed.  
Any multiplicative class containing the empty set for which the finite (countable) chain condition holds is evidently a 
semiring ($\sigma$-semiring).
\\

{\small\bf Example} \ 
Let $\jS \subset \jP (\XX)$ be a lattice $-$then the class of all sets of the form $S - T$, where $S, \hsx T \in \jS$ and $S \supset T$,  
is a semiring.  
Indeed, the condition as regards the empty set is trivial (take $S = T$).  
Let now $S_1 - T_1$ and $S_2 - T_2$ be in our class.  
Multiplicativity is then a consequence of the identity 

\[
(S_1 - T_1) \cap (S_2 - T_2) 
\ = \ 
(S_1 \cap S_2) - (S_1 \cap S_2) \cap (T_1 \cup T_2).  
\]
If in addition, $S_1 - T_1$ is contained in $S_2 - T_2$, then 
\[
S_1 - T_1
\subset 
(S_1 \cap S_2) - (T_1 \cap T_2) 
\subset
S_2 - T_2,
\]
from which it follows that the finite chain condition is in force, as can be seen by a direct set-theoretic calculation.
\\

{\small\bf Lemma 1} \ 
Let $\jS$ be a semiring; 
let $S_1, \ldots, S_m$ be nonempty, pairwise disjoint elements of $\jS$, contained in some fixed element $S$ of $\jS$ $-$then 
there exists a finite $\jS$-partition $P(S)$ of $S$ of the form 
\[
\{S_1, \ldots, S_m, S_{m+1}, \ldots, S_n\}.
\]
\\[-1cm]

\un{Proof} \ 
The proof is by induction on the integer $m$.  
If $m = 1$, then the assertion is true by the very definition of semiring.  
Assuming now the validity for $m \geq 1$, suppose that $\emptyset \neq T \subset S$ and intersects none of the $S_1, \ldots, S_m$ $-$then 
\[
T 
\ = \ 
T \cap S_{m+1} \cup \ldots \cup T \cap S_n
\qquad \text{(disjoint union)}.   
\]
In turn, making the obvious conventions, write
\[
\begin{matrix*}[l]
&S_{m+1} \ 
&= \ 
T \cap S_{m+1} \cup  S_{m+1} (1) 
\cup \ldots \cup 
S_{m+1} (\gamma_{m + 1}) 
\\[11pt]
&& \  \vdots
\\[11pt]
&S_n \ 
&= \ 
T \cap S_n \cup  S_n (1) 
\cup \ldots \cup 
S_n (\gamma_n) 
\end{matrix*}
\qquad \text{(disjoint union)}.
\]
Then
\[
\{S_1, \ldots, S_m, T, S_{m + i} (j)\}
\]
is an $\jS$-partition of $\jS$, thereby completing the proof.
\\

{\small\bf Lemma 2} \ 
Let $\jS$ be a semiring; 
let $S_1, \ldots, S_m$ be nonempty distinct elements of $\jS$ $-$then the union $S_1 \cup \ldots \cup S_m$ can be represented in the form 
\[
S_1 (1) 
\cup \ldots \cup 
S_1 (\gamma_1) 
\cup \ldots \cup 
S_m (1)
\cup \ldots \cup 
S_m (\gamma_m),
\]
where the $S_i (j)$ are nonempty, pairwise disjoint, belong to $\jS$, and 
\[
S_i \supset S_i (1), \ldots, S_i (\gamma_i) 
\qquad (i = 1, \ldots, m).
\]
\\[-1cm]

\un{Proof} \ 
The proof is by induction on the integer $m$.  
As the assertion is trivially true when $m = 1$, let us assume that it is valid for $m \geq 1$.  
Given $S_{m + 1}$, consider the $S_{m+1} \cap S_i (j)$.  
If each of these intersections is empty, then our contention is evident.  
Suppose, herefore, that 
$S_{m+1} \cap S_i (j) \neq \emptyset$ for certain $i$ and $j$ $-$then there are two possibilities:

\[
\begin{cases}
\ S_{m+1} \cap S_i (j)  \ = \ S_{m+1}
\\[8pt]
\ S_{m+1} \cap S_i (j)  \ \neq \ S_{m+1}
\end{cases}
.
\]

\noindent
If the first possibility obtains, then $i$ and $j$ are unique.  
Accordingly, in view of the definition of semiring, the difference $S_i (j) - S_{m+1}$, if not empty, can be written as a finite sum of nonempty, 
pairwise disjoint elements of $\jS$, leading, thereby, to the desired decomposition.  
If the second possibility obtains, then the $S_{m+1} \cap S_i (j)$ are proper, pairwise disjoint subsets of $S_{m+1}$.  
The proof can then be completed by an appeal to Lemma 1.
\\

We shall leave it up to the reader to decide if Lemmas 1 and 2 admit meaningful formulations in terms of $\sigma$-semirings, 
the issue being, of course, countable versus finite (cf. Exer. 5).  
\\[-.5cm]

In passing, it should be noted that the trace of a semiring ($\sigma$-semiring) is again a semiring ($\sigma$-semiring).
\\

{\small\bf Example} \ 
Take for $\XX$ a bounded, closed interval in $\R^n$, say: 
\[
\XX 
\ = \ 
\{x \hsy : \hsy a_1 \leq x_1 \leq b_1, \ldots, a_n \leq x_n \leq b_n\}.
\]
Let $\jS$ be the class consisting of the empty set and all intervals
\[
\{x \hsy : \hsy \alpha_1 \leq x_1 \leq \beta_1, \ldots, \alpha_n \leq x_n \leq \beta_n\}
\qquad 
(a_i \leq \alpha_i < \beta_i \leq b_i)
\]
if $\beta_i < b_i$ for every $i$, but if $\beta_i = b_i$ for some $i$, then the inequality $x_i < \beta_i$ is to be replaced by $x_i \leq \beta_i$.  
With this agreement, $\jS$ is a semialgebra.
By comparison, note that the class of all closed subintervals of $\XX$, while multiplicative, is not a semiring, 
although the class of all finite unions of such is a lattice. 
\\[-.25cm]

[Note: \ 
There are, of course, numerous simple variants on this theme.]
\\

\newpage 
\[
\textbf{\un{Notes and Remarks}}
\]

The notion of a semiring is frequently attributed to 
J. v. Neumann\footnote[1]{\vspace{.11 cm}\un{Functional Operators}, Annals of Mathematics Studies,\textbf{21} Princeton, (1950)}
This, however, is inaccurate, the priority belonging to 
A. Kolmogoroff\footnote[2]{\vspace{.11 cm}\textit{Math. Ann.},\textbf{103} (1930), 654-696.}
There one will find the term \unz{zerlegbarer Bereich} employed in context for what we have called semiring or $\sigma$-semiring.  
Actually, v. Neumann (op. cit.) did not work with semirings per se but rather with multiplicative classes satisfying the finite chain condition; 
the were called by him \unz{halfrings} (see p. 85 of that work).
The term semiring appears in 
Halmos\footnote[3]{\vspace{.11 cm}\un{Measure Theory}, D. Van Nostrand, New York, (1950)} (see p. 22), 
but still only in reference to the finite chain condition.  
Semirings were used early on by 
V. Glivenko\footnote[4]{\vspace{.11 cm}
{\fontencoding{OT2}\selectfont
V. Glivenko, 
Intgral Stieltjes, 
ONTI, 
Moskva-Leningrad,
}
(1936)
(see pp. 175-207).
}
in his book 
\un{The Stieltjes Integral}.
That semirings and  $\sigma$-semirings might be made the basis for measure theory was suggested by 
N. d. Bruijn and 
A. Zaanen\footnote[5]{\vspace{.11 cm}\textit{Indag. Math.},\textbf{16} (1954), 456-466.}
their perspective is quite different from that of Kolmogoroff's (op. cit.), being didactic rather than innovative.


\chapter{
$\boldsymbol{\S}$\textbf{3}.\quad  Exercises}
\setlength\parindent{2em}
\setcounter{theoremn}{0}
\renewcommand{\thepage}{\S3-\text{E-}\arabic{page}}



\qquad
(1) \quad 
Give an example of a semiring of finite cardinality which does not satisfy the finite chain condition.
\\

(2) \quad 
Give an example of a semiring of infinite cardinality which satisfies neither the finite chain condition nor the countable chain condition.
\\

(3) \quad 
Give an example of a semiring of infinite cardinality which does not satisfy the finite chain condition but does satisfy the countable chain condition.
\\

(4) \quad 
Give an example of a $\sigma$-semiring of finite cardinality which is not a semiring and which does not satisfy the countable chain condition.
\\

(5)  \quad 
Take $\XX = [0, 1[$ and consider the semiring $\jS$ consisting of all left closed and right open subintervals of $\XX$ $-$then
every $\jS$-partition of $\XX$ is finite or countable.  
Does there exist an $\jS$-partition $P(\XX)$ of $\XX$ such that each $[a,b[ \hsx \subsetx \XX$ $(a < b)$ with rational endpoints 
is partitioned by the components of $P(\XX)$ lying therein?
\\[-.5cm]

[What is the relevance of this exercise to Lemmas 1 and 2?]
\\

(6) \quad 
Let $\XX$ be a set of cardinality $n$, say.  In terms of $n$, how many semirings does $\jP (\XX)$ contain?
\\

(7) \quad 
By definition, a nonempty, bounded subset of $\R^n$ is called a \unz{convex polyhedron} provided that it can be written as a finite 
intersection of open or closed halfspaces.  
Show that the class consisting of the empty set and all convex polyhedra is a semiring satisfying the finite chain condition.
\\

(8) \quad 
Take for $\XX$ the Banach space $(c_0)$ of all real sequences $x = \{x_i\}$ which converge to zero, 
the norm being given by $\norm{x} = \sup \abs{x_i}$.  
Let $\{r_i(+)\}$ be a sequence in $\ov{\R}$ such that $0 < r_i(+) \leq +\infty$, 
$\varliminf \ r_i(+) > 0$; 
let $\{r_i(-)\}$ be a sequence in $\ov{\R}$ such that $0 > r_i(-) \geq -\infty$, 
$\varlimsup \ r_i(-) < 0$ $-$then 
by $S(\{r_i(-)\}, \{r_i(+)\})$ we understand the set of all $x \in \XX$ such that 
$r_i(-) \leq x_i < r_i(+)$ $\forall \ i$.  
Explain why the class consisting of the empty set and all possible 
$S(\{r_i(-)\}, \{r_i(+)\})$ 
is neither a semiring nor a $\sigma$-semiring.  
\\[-.5cm]

[Note: \ 
It was claimed to the contrary by 
P. Maserick\footnote[2]{\vspace{.11 cm}\textit{Pacific J. Math.}, \textbf{17} (1966), 137-148.}, 
that the class in question was a $\sigma$-semiring satisfying the countable chain condition.]


\chapter{
$\boldsymbol{\S}$\textbf{3}.\quad  Problem}
\setlength\parindent{2em}
\setcounter{theoremn}{0}
\renewcommand{\thepage}{\S3-\text{P-}\arabic{page}}
\vspace{-1.25cm}


\[
\textbf{NORMAL CLASSES} \ 
\]
\\[-1.5cm]

Let $\jS$ be a multiplicative class $-$then $\jS$ is said to be \unz{normal} if for any $S \in \jS$ admitting infinite $\jS$-partitions, each element 
$P(S) = \{S_1, \ldots, S_m, \ldots\}$ in $\sigma$-$\Par_\jS (S)$ has the property that $\forall \ m$, there exists a finite $\jS$-partition
\\[-1.cm]

\[
\{S_1, \ldots, S_m, T_1, \ldots, T_{r_m}\}
\]
\\[-1.5cm]

\noindent
of $S$.  
Every semiring is a normal class (cf. Lemma 1). 
\\[-.5cm]

(1) \quad
There exist multiplicative classes which are not normal.
\\[-.5cm]

[Take for $\XX$ a countable set $\{x_1, x_2, \ldots, \}$.  
Put $\jS = \{\emptyset, \XX, \{x_1\}, \{x_2\}, \ldots \}$ $-$then $\jS$ is multiplicative but not normal.]
\\[-.5cm]

(2) \quad
There exist $\sigma$-rings which are not normal classes.
\\[-.5cm]

[Let $\XX$ be an infinite set.  
Let $P(\XX) = \{X_1, \ldots, X_m\}$ be a finite partition of $\XX$ by subsets $X_i$, each of which we suppose in turn can be countably 
partitioned by subsets $X_{i \hsy j}$ $-$then the class $\jS$ consisting of $\emptyset$, $\XX$, the $X_i$ and the $X_{i \hsy j}$  is a 
$\sigma$-semiring but is not normal.]
\\[-.5cm]

(3) \quad
There exist normal classes which are not $\sigma$-rings.
\\[-.5cm]

[Take for $\XX$ a countable set $\{x_1, x_2, \ldots, \}$. 
Put $\jS = \{\emptyset, \XX, \{x_2\}, \{x_3\}, \ldots, \{x_i, x_{i + 1} (i = 2, 3, \ldots) \}$ 
$-$then $\jS$ is normal but is not a $\sigma$-ring.]
\\[-.25cm]

Let $\jS$ be a multiplicative class $-$then $\jS$ is normal iff for any $S \in \jS$ admitting infinite $\jS$-partitions and for any 
$P(S) \in \sigma \dash \Par_\jS (S)$ each element 
$P(S) = \{S_1, \ldots, S_m, \ldots\}$ in $\sigma$-$\Par_\jS (S)$ has the property that $\forall \ m$, 
there exists a finite $\jS$-partition
\[
\{S_1, \ldots, S_m, T_1, \ldots, T_{r_m}\}
\]
refining $P(S)$.  
\\[-.5cm]

\unz{Ref} \quad 
D. Procenko\footnote[3]{{\fontencoding{OT2}\selectfont
D. Procenko}, \vspace{.11 cm}\un{Soob\v s\v c. Acad. Nauk Gruzin. SSR}, \textbf{40} (1965), 271-278.}.


\chapter{
$\boldsymbol{\S}$\textbf{4}.\quad  Rings, $\sigma$-Rings, $\delta$-Rings}
\setlength\parindent{2em}
\setcounter{theoremn}{0}
\renewcommand{\thepage}{\S4-\arabic{page}}



\qquad 
Let $\XX$ be a nonempty set; 
let $\jS$ be a subset of $\jP (\XX)$ containing the empty set $-$then $\jS$ is said to be a 
\un{ring} if 
\[
S, \hsx T \in \jS 
\implies 
S \hsy \Delta \hsy T \in \jS 
\quad \text{and} \quad 
S \cap T \in \jS.
\]
Since 
\[
\begin{cases}
\ S \cup T \ = \ (S \hsy \Delta \hsy T) \hsy \Delta \hsy (S \cap T)
\\[8pt]
\ S - T \ = \ S \hsy \Delta \hsy (S \cap T)
\end{cases}
,
\]
a ring is closed under the formation of finite unions and differences and, in fact, is characterized by these requirements.  
An \un{algebra} is a ring containing $\XX$.  
Trivially, $\{\emptyset\}$ is a ring while $\{\emptyset, \XX\}$ and $\jP (\XX)$ are algebras.
\\

{\small\bf Example} \ 
(Kolmogoroff) \ 
Any ring is a semiring.  
We have seen in \S3 that every lattice gives rise in a natural manner to a semiring; 
in turn, every semiring gives rise in a natural manner to a ring.  
Thus let $\jS$ be a semiring and consider the class $\jKol (\jS)$ of all sets of the form 
$\bigcup\limits_{i = 1}^m \ S_i$, the $S_i$ being elements of $\jS$ which, without loss of generality, 
can be taken pairwise disjoint (cf. Lemma 2 (\S3)) $-$then we claim that $\jKol (\jS)$ is a ring.  
Indeed, if 
$S = \bigcup\limits_{i = 1}^m \  S_i$, 
$T = \bigcup\limits_{j = 1}^n \ T_j$ 
are disjoint unions of elements $S_i, \hsx T_j \in \jS$, then so is
\[
S \cap T 
\ = \ 
\bigcup\limits_{i = 1}^m \ 
\bigcup\limits_{j = 1}^n 
(S_i \cap T_j).
\]
As for $S \hsy \Delta \hsy T$, use Lemma 1 (\S3) to write
\[
\begin{cases}
\ S_i \ = \ 
\bigcup\limits_{j = 1}^n 
(S_i \cap T_j) 
\hsx \cup \hsx 
\bigcup\limits_{k = 1}^{r_i} \ 
S_{i k}
\\[8pt]
\ T_j \ = \ 
\bigcup\limits_{i = 1}^m \ 
(S_i \cap T_j) 
\hsx \cup \hsx 
\bigcup\limits_{k = 1}^{r_j} \ 
T_{j k}
\end{cases}
\qquad \text{(disjoint union)}.
\]
Then we have
\[
S \hsy \Delta \hsy T
\ = \ 
\bigcup\limits_{i = 1}^m \ 
\left(
\bigcup\limits_{k = 1}^{r_i} \ 
S_{i k}
\right)
\hsx \cup \hsx 
\bigcup\limits_{j = 1}^n 
\left(
\bigcup\limits_{k = 1}^{r_j} \ 
T_{j k}
\right),
\]
which again is a disjoint union of elements in $\jS$.  
Accordingly, the class $\jKol (\jS)$ is a ring.
\\[-.25cm]

[Note: \ 
Suppose that $P(\XX) = \{X_i \hsy : \hsy i \in I\}$ is a partition of $\XX$ $-$then 
the class consisting of $\emptyset$ and the $X_i$ is a semiring.  
Therefore the class formed by the empty set and all nonempty finite unions of the components of $P(\XX)$ is a ring.]
\\

The justification of the term ``ring of sets'' lies in the following remarks.  
In $\jP (\XX)$ itself, introduce operations of addition and multiplication via the stipulations
\[
\begin{cases}
\ S + T \ \equiv S \hsx \Delta \hsx T
\\[8pt]
\ S \cdot T \ \equiv S \hsx \cap \hsx T
\end{cases}
.
\]
Then by an elementary if slightly tedious verification, one checks that $\jP (\XX)$ thus equipped is a commutative ring 
with zero element
$\emptyset$ and multiplicative identity $\XX$.  
It is a point of some importance that these operations, when viewed as maps 
\[
\jP (\XX) \times \jP (\XX)  \ra \jP (\XX),
\]
are jointly continuous, i.e., $\jP (\XX)$ is a topological ring; 
on the other hand, these operations, when viewed as maps
\[
\jP (\XX)_\tS \times \jP (\XX)_\tS \ra \jP (\XX)_\tS
\]
are separately continuous.
\\

Utilizing now the customary algebraic terminology, a subring of $\jP (\XX)$ is a subset containing the zero element, 
i.e., $\emptyset$, and closed under addition and multiplication or still, under symmetric differences and intersections; 
in other words, subring of $\jP (\XX) = $ ring of subsets of $\XX$.  
In addition, a subalgebra of $\jP (\XX)$is a subring containing the multiplicative identity, i.e., $\XX$; 
in other words: \ subalgebra of $\jP (\XX) = $ algebra of subsets of $\XX$.  
\\[-.25cm]

[Note: \ 
A ring (algebra) of sets is evidently a Boolean ring (algebra).  
It must be stressed, however, that a ring $\jS$ may well admit a multiplicative identity, thus is a Boolean algebra, but is
not an algebra, the point being that generally $\XX \notin \jS$.  
Consider, e.g., $\jS = \jP (S)$, $S$ a nonempty proper subset of $\XX$.  
Accordingly, we shall use the term \un{ring with unit} to refer to a ring $\jS$ possessing a multiplicative identity; 
in particular, therefore, every algebra is a ring with unit.  
It is easy to check that a ring $\jS$ is a ring with unit iff $\bigcup \ \jS \in \jS$.  
If 
$\bigcup \ \jS \notin \jS$, 
then the class $\widehat{\jS}$ consisting of all $S$, $\bigcup \ \jS - S$ $(S \in \jS)$ is a ring with unit containing $\jS$.  
Finally, it should be recalled that every Boolean ring is of characteristic 2, hence may be regarded as an algebra 
over the field $\Z_2$.]
\\

The usual algebraic notions then admit easy descriptive interpretations.  
Consider e.g., the notion of an ideal $\sI$ in the ring $\jS$ $-$then, descriptively, $\sI$ can be characterized as a 
nonempty subclass of $\jS$ which is closed under the formation of finite unions and is hereditary in the sense that 
$I \in \sI$, $S \in \jS$, $S \subset I \implies S \in \sI$.  
The corresponding quotient $\jS / \sI$ is a Boolean ring, elements $S, \hsx T \in \jS$ being equivalent 
$\mod \sI$ iff $S \hsx \Delta \hsx T \in \sI$ or still, iff $S = (T - I) \cup J$ $(I, \hsx J \in \sI)$.
\\


{\small\bf Lemma 1} \ 
Let $\jS$ be a ring; 
let $\sI \neq \jS$ be an ideal $-$then
\\[-.5cm]

\qquad (1) \quad 
$\sI$ is contained in a maximal ideal; 
\\[-.25cm]

\qquad (2) \quad 
$\sI$ is maximal iff $\sI$ is prime; 
\\[-.25cm]

\qquad (3) \quad 
$\sI$ is the intersection $\bigcap\limits_{\fp \supset \sI} \ \fp$, 
\quad $\fp$ prime.
\\[-.25cm]

[There is nothing to be gained by giving a proof in extenso.  
The point is this.  
$\jS$ need not have a multiplicative identity and, as is well known, if a ring does not have a multiplicative identity, 
then, e.g., generic ideals need not be contained in maximal ideals, maximal ideals need not be prime, 
prime ideals need not be maximal, etc.  
But $\jS$ is a Boolean ring, hence carries compensating structure.  
To illustrate, consider (1).  
Since $\sI \neq \jS$, $\exists \ S_0 \in \jS$, $S_0 \notin \sI$.  
Let $\fm$ be any ideal in $\jS$ maximal with respect to the property that $\fm \supset \sI$, $S_0 \notin \fm$ 
(Zorn's lemma ensures the exisitence of $\fm$) $-$then $\fm$ is in fact a maximal ideal, 
as can be checked without difficulty $(S_0^2 = S_0!)$.
Statement (2) is also easy, as is (3).]
\\

A \un{$\sigma$-ring} is a ring $\jS$ which is closed under the formation of countable unions, i.e., 

\[
\{S_i \hsy : \hsy i = 1, 2, \ldots \} \subset \jS 
\implies 
\bigcup \ S_i \in \jS,
\]
or still, $\jS = \jS_\sigma$.  
A \un{$\sigma$-algebra} is a $\sigma$-ring containing $\XX$.  
A  \un{$\delta$-ring} is a ring $\jS$ which is closed under the formation of countable intersections, i.e., 
\[
\{S_i \hsy : \hsy i = 1, 2, \ldots \} \subset \jS 
\implies 
\bigcap \ S_i \in \jS,
\]
or still, $\jS = \jS_\delta$.  
A \un{$\delta$-algebra} is a $\delta$-ring containing $\XX$.  
\\[-.5cm]

A \un{$\sigma$-ideal} ($\delta$-ideal) 
is an ideal in a ring which is closed under the formation of countable unions (intersections).
\\

{\small\bf Example} \ 
Let $\XX$ be a topological space $-$then the class $\jS$ of all subsets of $\XX$ having the Baire property is a 
$\sigma$-algebra containing the $\sigma$-ideal of all first category subsets of $\XX$.
\\[-.4cm]

[Note: \ 
Recall that a set $S \subset \XX$ is said to have the 
\un{Baire property} 
if there exists an open set $G$ such that $S - G$ and $G - S$ are of the first category.]
\\

A $\sigma$-ring is a $\delta$-ring.  
To see this, put 
$S = \bigcup \ S_i$ $(S_i \in \jS)$ $-$then

\[
\bigcap \ S_i 
\ = \ 
S - \bigcup
\left(
S - S_i
\right).
\]
Consequently, if $\{S_i\}$ is a sequence of sets in a $\sigma$-ring $\jS$, then

\[
\varlimsup \ S_i \in \jS, 
\quad 
\varliminf \ S_i \in \jS.
\]
In particular: \ 
A $\sigma$-ring is necessarily closed in $\jP (\XX)_\tS $.  
Furthermore, due to the separate continuity of the operations

\[
\jP (\XX)_\tS \times \jP (\XX)_\tS \ra \jP (\XX)_\tS, 
\]
the closure in $\jP (\XX)_\tS$ of a ring is again a ring, thus is actually a $\sigma$-ring.
\\

{\small\bf Example} \ 
There are $\delta$-rings which are not $\sigma$-rings.  
For instance, take $\XX = \R^n$ and consider the class $\jS$ of all relatively compact subsets.
\\

{\small\bf Lemma 2} \ 
Let $\jS$ be a ring $-$then $\jS$ is a $\delta$-ring iff for every $S_0 \in \jS$, the set 
$\{S \in \jS \hsy : \hsy S \subset S_o\}$ is a $\sigma$-algebra in $S_0$.
\\[-.25cm]

[We omit the elementary verification.]
\\

It follows from Lemma 2 that every $\delta$-ring which admits a multiplicative identity is necessarily a $\sigma$-ring.
\\[-.5cm]

A ring $\jS$ is said to be \un{complete} if $\jS$ is closed under the formation of arbitrary nonempty unions.  
A complete ring is evidently also closed under the formation of arbitrary nonempty intersections.  
If $\jS$ is complete, then $\jS$ is a ring with unit $\bigcup \ \jS$; 
of course 
$\bigcup \ \jS \neq \XX$ in general, hence $\jS$  need not be an algebra.
\\


{\small\bf Example} \ 
Let $\XX$ be a set of cardinality $\aleph_0$; 
let $\jS$ be a $\sigma$-ring in $\XX$ $-$then $\jS$ is complete.
\\

Consider $\jP (\XX)$, equipped with the topology of pointwise convergence $-$then a net $\{S_i\}$ in $\jP (\XX)$ 
is convergent with limit $S$, say, iff it is order convergent, i.e., 
\[
\bigcap\limits_i \ \bigcup\limits_{j \geq i} \ S_j 
\ = \ 
\bigcup\limits_i \ \bigcap\limits_{j \geq i} \ S_j ,
\]
the order limit being exactly $S$.
\\[-.25cm]

This being so, suppose that $\jS$ is a complete ring in $\XX$ $-$then $\jS$ is closed in $\jP (\XX)$.  
If $\jS$ is a ring but is not complete, then the closure $\bar{\jS}$ of $\jS$ in $\jP (\XX)$ is a complete ring in $\XX$, 
the \un{completion} of $\jS$.
Every complete subring of $\jP (\XX)$ containing $\jS$ must contain $\bar{\jS}$, therefore the completion of $\jS$ 
is the minimal complete ring in $\XX$ containing $\jS$ or still, the complete ring generated by $\jS$ (cf. \S6).
\\

{\small\bf Example} \ 
Let $\jS$ be a ring in $\XX$.  
Suppose that $\forall \ x \in \XX$, $\{x\} \in \jS$ $-$then the completion of $\jS$ is $\jP (\XX)$.
\\

Let $\jS$ be a ring $-$then a nonempty subset $\tA \in \jS$ is said to be an \un{atom} if, 
apart from the empty set, A properly contains no other elements of $\jS$.  
We write $\At (\jS)$ for the class of all atoms in $\jS$.
\\[-.5cm]

If every nonempty $S \in \jS$ contains an atom, then $\jS$ is said to be \un{atomic}; 
on the other hand, if no nonempty $S \in \jS$ contains an atom, then $\jS$ is said to be \un{antiatomic}.
\\

{\small\bf Example} \ 
Let $\XX$ be a Hausdorff topological space, $\XX_\isol$ its set of isolated points $-$then $\XX$ can be written as a 
disjoint union 
$\XX = \XX_\perf \cup \XX_\scat$, 
where 
$\XX_\perf$ 
is the perfect kernel of $\XX$, i.e., the union of all subsets of $\XX$ which are dense in themselves, 
and 
$\XX_\scat \supset \XX_\isol$
is the corresponding complement.  
$\XX_\perf$ is closed while 
$\XX_\scat$ is open; 
one of them may, of course, be empty.  
Assume now that $\XX$ is in addition, locally compact and totally disconnected.  
Consider the ring $\jS$ of all open and compact subsets of $\XX$ $-$then 
$\At (\jS) = \{\{x\} \hsy : \hsy x \in \XX_\isol\}$, so
\[
\begin{cases}
\ \jS \ \text{is atomic iff} \ \XX = \bar{\XX}_\isol
\\[8pt]
\ \jS \ \text{is antiatomic iff} \ \XX = \bar{\XX}_\perf
\end{cases}
.
\]
In this connection, note that $\XX = \XX_\perf$ iff $\XX_\scat = \emptyset$ 
but 
$\XX = \bar{\XX}_\isol$ 
does not imply that 
$\XX_\perf = \emptyset$, 
as can be seen by example.  
It is also easy to envision intermediate situations, a particularly transparent case being when 
$\XX$ is extremally disconnected.
\\


Any complete ring $\jS$  is atomic, there being an easy characterization of the atoms.  
Thus define an equivalence relation in $\bigcup \ \jS$ by requiring that $x$ be equivalent to $y$ iff 
every set in $\jS$  which contains $x$ also contains $y$.  
The equivalence class $[x]$ $(x \in \bigcup \ \jS)$ belongs to $\jS$, as can be seen by noting that
\[
[x] 
\ = \ 
\bigcap\limits_{x \in S} \ S
\qquad (S \in \jS).
\]
The atoms of $\jS$ are just the $[x]$ $(x \in \bigcup \ \jS)$.  
Every nonempty $S \in \jS$ is partitioned by the atoms which it contains.
\\[-.25cm]

Let now $\jS$ be an arbitrary ring in $\XX$ $-$then there is a canonical map
\[
\phi \hsy : \hsy \jS \ra \jP (\At (\jS)), 
\]
namely the rule which assigns to each $S \in \jS$ the class $\phi(S)$ of all atoms $A \subset S$.  
It is clear that $\phi$ is a homomorphism of rings.  
Furthermore: \ 
\\[-.25cm]

\qquad (1) \quad 
If $\jS$ is atomic, then $\phi$ is injective.  
Indeed, if $S, \hsx T \in \jS$, $S \neq T$, then $S - T \neq \emptyset$, say, thus 
$\exists \ A \in \At (\jS)$, $A \subset S - T$, and so $A \in \phi(S)$, $A \notin \phi (T)$.
\\[-.25cm]


\qquad (2) \quad 
If $\jS$ is complete, then $\phi$ is surjective.  
Indeed, if $\{A_i\}$ is any class of atoms, then $\bigcup \ A_i \in \jS$ and $\phi (\bigcup \ A_i) = \{A_i\}$.  
\\[-.25cm]

We have seen above that every complete ring is atomic.  
Therefore, in this case, $\phi$ is an isomorphism of rings.  
We remark that $\phi$ is then even a complete isomorphism in that it preserves arbitrary unions and intersections. 
\\[-.5cm]

In passing, it should be noted that the trace of a 
ring ($\sigma$-ring, $\delta$-ring) is again a 
ring ($\sigma$-ring, $\delta$-ring), the same also being true of complete rings.
\\

\newpage 
\[
\textbf{\un{Notes and Remarks}}
\]

The theory presented in this $\S$ can be approached more generally, viz. from the point of view of abstract Boolean rings and 
Boolean algebras; cf. 
R. Sikorski\footnote[1]{\vspace{.11 cm}\un{Boolean Algebras}, Springer-Verlag, Berlin (1969).}, 
as well as 
D. Ponasse and 
J-C. Carrega\footnote[2]{\vspace{.11 cm}\un{Alg\`ebre et Topologie Bool\'eennes}, Masson, Paris (1979).}.  
The terminology, particularly in the older literature, is tangled.  
Specifically, what we have termed a lattice is frequently called a ring while what we have termed a ring is frequently called a field; 
cf. 
F. Hausdorff\footnote[3]{\vspace{.11 cm}\un{Grundz\"uge der Mengenlehre}, Veit \& Comp., Leipzig, (1914).}, 
(see pp. 14-16), the German being \unz{Ring} and \unz{K\"orper}, respectively.  
To compound the confusion, 
M. Fr\'echet\footnote[4]{\vspace{.11 cm}\textit{Bull. Soc. Math. France}, \textbf{43} (1915), 248-265.}, 
refers to a $\sigma$-ring as a \unz{famille additive d'ensembles}, whereas 
O. Nikodym\footnote[5]{\vspace{.11 cm}\textit{Fund. Math.}, \textbf{15} (1930), 131-179.} 
understands by \unz{corps d'ensembles} a $\sigma$-algebra.  
There are other permutations and combinations too; 
e.g., 
R. de Possel\footnote[6]{\vspace{.11 cm}\textit{J. Math. Pures Appl. (9)}, \textbf{15} (1936), 391-409.}, 
has suggested \unz{tribe} (\unz{tribu} in French) for $\sigma$-ring, a \unz{clan} then being a ring.  
In the sense employed in the text, the term ring appears in 
J. v. Neumann\footnote[7]{\vspace{.11 cm}\un{Functional Operators}, 
Annals of Mathematics Studies vol 21. Princeton,  (1950).}, (see p. 84). 
That semirings lead naturally to rings was pointed out by 
A. Kolmogoroff\footnote[8]{\vspace{.11 cm}\textit{Ann. Math.}, \textbf{103} (1930), 654-696.}.  
Ideals in rings have been investigated systematically by 
A. Tarski\footnote[9]{\vspace{.11 cm}\textit{Fund. Math.}, \textbf{32} (1939), 45-63.},
A. Tarski\footnote[9]{\vspace{.11 cm}\textit{Fund. Math.}, \textbf{33} (1945), 51-65.}, 
A. Tarski\footnote[1]{\vspace{.11 cm}\textit{Soc. Sci. Lett. Varsovie C. R. C1. III Sci. Math. Phys.}, \textbf{30} (1937), 151-181.}.  
The notion of atom is generally attributed to 
M. Fr\'echet\footnote[2]{\vspace{.11 cm}\textit{Fund. Math.}, \textbf{5} (1924), 206-251.} 
although it can be traced back to 
E. Sch\"oder\footnote[3]{\vspace{.11 cm}\un{Vorlesungen \"uber die Algebra der Logik}, 
 II (Bd. I), B.G. Teubner, Leipzig, 1891}, (see \S47).  
The fact that every complete ring is isomorphic to the power set of its atoms is due to 
Lindenbaum and Tarski; 
cf. 
A. Tarski\footnote[4]{\vspace{.11 cm}\textit{Fund. Math.}, \textbf{24} (1935), 177-198.}.


\chapter{
$\boldsymbol{\S}$\textbf{4}.\quad  Exercises}
\setlength\parindent{2em}
\setcounter{theoremn}{0}
\renewcommand{\thepage}{\S4-\text{E-}\arabic{page}}



(1) \quad 
Take $\XX = \R$.  
For $n = 0, 1, \ldots,$ let $\jD_n$ be the class consisting of the empty set and all nonempty finite disjoint unions of 
dyadic left closed and right open intervals of order $n$, i.e., the 
$
\bigg[
\frac{k}{2^n}, \frac{k+1}{2^n}
\bigg[
$.
Verify that $\jD_n$ is a ring.  
Noting that $\jD_0 \subset \jD_1 \subset \ldots$, put $\jD = \bigcup \ \jD_n$, 
the class of all finite unions of dyadic left closed and right open intervals of any order.  
Verify that $\jD$ is a ring.  
Formulate and prove a multidimensional generalization.
\\[-.25cm]

[Observe that
\[
\bigg[
\frac{k}{2^n}, \frac{k+1}{2^n}
\bigg[
\ = \ 
\bigg[
\frac{2k}{2^{n+1}}, \frac{2k+1}{2^{n+1}}
\bigg[
\ \cup \ 
\bigg[
\frac{2k + 1}{2^{n+1}}, \frac{2k+2}{2^{n+1}}
\bigg[
\  .]
\]

(2) \quad 
Let $\XX$ be a topological space $-$then the class $\jS$ comprised of all sets $S \subset \XX$ whose boundary is nowhere dense is an algebra of subsets of $\XX$.
\\

(3) \quad 
Let $\XX$ be a nonempty set $-$then the class $\jS$ comprised of all sets $S \subset \XX$ such that 
either $\card(S) < \aleph_\alpha$ or $\card(\XX - S) < \aleph_\alpha$ is an algebra of subsets of $\XX$.
\\

(4) \quad
Given a ring $\jS$, consider the following conditions: 
\\[-.25cm]

\qquad ($\tC_1$) \quad 
Every subset of $\jS$ consisting of nonempty, pairwise disjoint elements if finite or countable.
\\[-.25cm]

\qquad ($\tC_2$) \quad 
Every subset of $\jS$ consisting of nonempty, pairwise comparable elements is finite or countable.
\\[-.25cm]

Show by example that there exist infinite rings which satisfy ($\tC_1$) but not ($\tC_2$) and vice-versa.
\\

(5)  \quad
Let $\jS$ be a subset of $\jP(\XX)$ containing the empty set.  
Suppose that $\jS$ is multiplicative 
$-$then the following are equivalent: 
\\[-.25cm]

\qquad (i) \quad 
The class $\jS$ is a semiring;
\\[-.25cm]

\qquad (ii) \quad 
The class consisting of all sets expressible as a finite union of pairwise disjoint sets from $\jS$ is a ring; 
\\[-.25cm]

\qquad (iii) \quad 
Given elements $S_1, \ldots, S_m$ of $\jS$, there exist pairwise disjoint elements $T_1, \ldots, T_n$ of $\jS$ 
such that each $S_i$ is the union of certain of the $T_j$.
\\

(6) \quad 
Let $\jS$ be a semiring.  
Consider the class of all sets of the form 
$\bigcup\limits_{i = 1}^\infty \ S_i$, the $S_i$ being elements of $\jS$, which, without loss of generality, 
can be taken pairwise disjoint (cf. Lemma 2 (\S3)).  
Show by example that this class need not be a ring.
\\

(7) \quad 
True or False? \ $\jP(\XX)_\ts$ is a topological ring, i.e., the operations of addition and multiplication
\[
\jP(\XX)_\ts \times \jP(\XX)_\ts \ra \jP(\XX)_\ts
\]
are jointly continuous.
\\[-.25cm]

[Is Exer. 10 (\S1) relevant here?]
\\

(8) \quad 
Let $\jS$ be a ring $-$then the following are equivalent: 
\\[-.25cm]

\qquad (i) \quad 
$\jS$ admits a nonprincipal prime ideal;
\\[-.25cm]

\qquad (ii) \quad 
$\jS$ admits a nonprincipal ideal; 
\\[-.25cm]

\qquad (iii) \quad 
$\jS$ is infinite.
\\[-.25cm]

[If (iii) is in force, then $\jS$ must possess countably many nonempty pairwise disjoint elements (cf. Exer. 5 (\S2)).]
\\

(9) \quad
In a ring with unit, there is a natural one-to-one correspondence between ideals and filters, the two concepts being dual 
to one another; 
under this correspondence, prime ideals are matched with ultrafilters.
\\[-.25cm]

[Let $\jS$ be a ring with unit $\bigcup \ \jS$ $-$then the correspondence in question is simply the complementation relative to 
$\bigcup \ \jS$.]
\\

(10) \quad 
Let $\XX$ be a locally compact, totally disconnected, Hausdorff space; 
let $\jS$ be the ring of open and compact subsets of $\XX$.  
Is $\jS$ a $\sigma$-ring?
\\

(11) \quad 
Let $\XX$ and $\YY$ be nonempty sets; 
let $f : \XX \ra \YY$ be a map $-$then
there is an induced map $f^{-1} : \jP(\YY) \ra  \jP(\XX)$ .  
Show that if $\jT$ is a ring ($\sigma$-ring) in $\YY$, 
then $\{f^{-1} (T) \hsy : \hsy T \in \jT\}$ is a ring ($\sigma$-ring) in $\XX$, 
and if $\jS$ is a ring ($\sigma$-ring) in $\XX$, then $\{T \subset Y  \hsy : \hsy f^{-1} (T) \in \jS \}$ is a ring 
($\sigma$-ring) in $\YY$.  
Are these assertions true if ring ($\sigma$-ring) is replaced by algebra ($\sigma$-algebra)?
\\

(12) \quad 
Let $\jS$ be a $\sigma$-ring in $\XX$ not containing $\XX$ $-$then the classes
\[
\begin{cases}
\ \{A \subset \XX \hsy : \hsy A \in \jS \ \text{or} \ \XX - A  \in \jS \}
\\[8pt]
\ \{A \subset \XX \hsy : \hsy S \in \jS \implies A \cap S \in \jS\}
\end{cases}
\]
are $\sigma$-algebras in $\XX$ containing $\jS$, the latter containing the former.
\\

(13) \quad
Prove that there does not exist an infinite $\sigma$-algebra $\jS$ with countably many members.  
Can $\sigma$-algebra be replaced by $\sigma$-ring in this assertion?
\\[-.25cm]

[Bear in mind Exer. 5 (\S2).]
\\

(14) \quad
Let $\jS_1 \subset \jS_2 \subset\ldots$ be a strictly increasing chain of subsets of $\jP(\XX)$.  
Show that if the $\jS_i$ are algebras in $\XX$, then the union $\bigcup \ \jS_i$ is again an algebra in $\XX$ 
but if the $\jS_i$ are $\sigma$-algebras in $\XX$,  then the union $\bigcup \ \jS_i$ is never a $\sigma$-algebra in $\XX$.  
What happens if , instead, the $\jS_i$ are rings?
\\[-.25cm]

[To discuss the second assertion, first show that there exists a sequence $\{S_i\}$ of nonempty, pairwise disjoint sets 
$S_i \hsy : \hsy S_i \in \jS_{i + 1} - \jS_i$ $\forall \ i$ (change the indexing if necessary).  
This done, proceed by contradiction and assume that $\bigcup \ \jS_i$ is a  $\sigma$-algebra $-$then eventually the 
\[
\fN_i 
\ = \ 
\big\{ S \subset \N \ : \ \bigcup\limits_{j \in S} \ S_j \ \in \ \jS_i \big\}
\]
are  $\sigma$-algebras in $\N$.]
\\

(15)  \quad
True or False? \ Let $\aleph$ be an infinite cardinal.  
Let $\XX$ be a set of cardinality $\aleph$; 
let $\jS$ be a ring in $\XX$ which is closed under the formation of unions of cardinality $\leq \aleph$ 
$-$then $\jS$ is complete.
\\[-.25cm]

(16) \quad 
Let $P(\XX) = \{X_i : i \in I\}$ be a partition of $\XX$ $-$then the class consisting of the empty set and all possible nonempty 
unions of the $X_i$ is a complete algebra.  
Conversely, let $\jS$ be a complete algebra $-$then there exists a partition 
$P(\XX) = \{X_i : i \in I\}$
of $\XX$ such that the class consisting of the empty set and all possible nonempty unions of the $X_i$ is $\jS$.
\\[-.25cm]

[Note: \ 
The correspondence between partitions and complete algebras is evidently one-to-one.]
\\

(17) \quad 
A ring $\jS$ such that it and all its subrings are atomic is called \unz{superatomic}.
\\[-.25cm]

True or False? \ There exist infinite superatomic rings.
\\

(18) \quad 
Let $\jS$ be a ring $-$then the following are equivalent: 
\\[-.25cm]

\qquad (i) \quad 
There exists a prime ideal containing $\jAt(\jS)$; 
\\[-.25cm]

\qquad (ii) \quad 
There exists a proper ideal containing $\jAt(\jS)$; 
\\[-.25cm]

\qquad (iii) \quad 
There exists an infinite class $\{S_i\} \subset \jS$ of nonempty, pairwise disjoint sets $S_i$ and a set $S \in \jS$ 
such that $\bigcup \ S_i \subset S$.
\\[-.25cm]

[What, if any, is the connection between the three conditions figuring here and the three which appear in Exer. 8?]
\\

(19) \ 
Let $\jS$ be a ring; let $\aleph$ be an infinite cardinal $-$then the following are equivalent: 

\qquad (i) \quad 
$\jS$ is complete and the cardinality of $\jAt(\jS)$ is $\aleph$; 
\\[-.25cm]

\qquad (ii) \quad 
$\jS$ is closed under the formation of unions of cardinality $\leq \aleph$ and 
$\aleph$ is the largest cardinal for which there exists a class $\jX \subset \jS$ of cardinality $\aleph$ comprised of 
nonempty, pairwise disjoint sets.
\\

(20) \quad 
Let $\jS$ be a ring $-$then the following are equivalent: 
\\[-.25cm]

\qquad (i) \quad 
$\jS$ is complete and $\jAt(\jS)$ is countable.
\\[-.25cm]

\qquad (ii) \quad 
$\jS$ is an infinite $\sigma$-ring with the property that every class $\jX \subset \jS$ of nonempty, 
pairwise disjoint sets is finite or countable.
\\[-.25cm]

[What additional fact must be cited in order to make this exercise a corollary to the preceding exercise?]
\\[-.25cm]

Taking into account Exer. 16, explicate the significance of this result for the collection of countable partitions of $\XX$.
\\

(21) \quad
Construct an example of an atomic ring $\jS$ possessing elements which cannot be written as a union of atoms.
\\

(22) \quad 
True or False? 
\\[-.25cm]

\qquad (a) \quad 
If $\jS$ is an antiatomic ring, then every nonempty $S \in \un{\jS}$ contains $\aleph_0$ nonempty,
pairwise disjoint sets $S_i \in \un{\jS}$.
\\[-.25cm]

\qquad (b) \quad 
If $\jS$ is an antiatomic $\sigma$-ring, then every nonempty $S \in \un{\jS}$ contains $\aleph_1$ nonempty,
pairwise disjoint sets $S_i \in \un{\jS}$.
\\[-.25cm]

(23) \quad 
Let $\jL_\XX$ stand for the collection of $\sigma$-algebras on $\XX$.  
Given $\jS^\prime$, $\jS^{\prime\prime} \in \jL_\XX$, write 
$\jS^\prime \leq \jS^{\prime\prime}$ if $\jS^\prime \subset \jS^{\prime\prime}$ $-$then, 
with this definition of order, $\jL_\XX$ is a complete lattice with largest and smallest elements.  
However, in general, $\jL_\XX$ is neither distributive nor modular.  
If $\card(\XX) \leq \aleph_0$, then $\jL_\XX$ is isomorphic to the partition lattice on $\XX$ (cf. Exer. 3 (\S2)), 
thus is complemented but, as can be shown, this fails if $\card(\XX) > \aleph_0$.
\\[-.25cm]

[Note: \ 
It is necessary to admit here the notion of generated $\sigma$-algebra (see \S6).  
For details (and additional information), see 
K. Bhaskara Rao and 
B. Rao\footnote[4]{\vspace{.11 cm}\textit{Dissertationes Math.}, \textbf{190} (1981), 1-68.}
.]
\\

(24) \quad 
Every abstract lattice is isomorphic to a sublattice of $\jL_\XX$ for some $\XX$.
\\[-.25cm]

[Combine the theorem of Whitman (Exer. 3 (\S2)) with Exer. 16.]


\chapter{
$\boldsymbol{\S}$\textbf{4}.\quad  Problem}
\setlength\parindent{2em}
\setcounter{theoremn}{0}
\renewcommand{\thepage}{\S4-\text{P-}\arabic{page}}



\textbf{I} \ TOPOLOGICAL REPRESENTATION OF BOOLEAN RINGS
\\[-.25cm]

Rings of sets and their quotients are the Boolean rings of primary importance in analysis.  
To deal with both simultaneously, it is most economical to consider an arbitrary Boolean ring.  
Such rings were studied intensively by Stone in the 1930's.  
The foundational results of this theory, a sketch of which will be given below, can be regarded as but simple exercises in the modern theory of schemes.  
Accordingly, the reader who is familiar with the language of contemporary algebraic geometry should have no difficulty filling in the omitted details.  
\\[-.25cm]

By a \unz{Boolean space}, we shall understand a topological space $\XX$ whose topology $\jT$ is locally compact, totally disconnected, and Hausdorff.  
Open subsets of a Boolean space are Boolean spaces, as are the closed subsets.  
Associated with every Boolean space $\XX$ is a ring $\jA(\XX)$, viz. the ring of open and compact subsets of $\XX$.  
The prime ideals in $\jA(\XX)$ are parameterized by the points $x \in \XX : \fp_x$ = elements of $\jA(\XX)$ not containing $x$.  
\\[-.25cm]

[Note: \ 
Owing to the Urysohn metrization theorem, a compact Boolean space is metrizable iff the cardinality of $\jA(\XX)$ is $\leq \aleph_0$.]
\\[-.25cm]

Let $\jA$ be a Boolean ring $-$then attached to $\jA$ is the set $\Spec(\jA)$ of all prime (= maximal) ideals of $\jA$.  
Given $f \in \jA$, put
\[
\jA_f \
\ = \ 
\{\fp \in \Spec(\jA) \hsy : \hsy f \notin \fp\}.
\]
Then the map
\[
\jA \ra \fp(\Spec(\jA))
\]
which assigns to each $f$ in $\jA$ the set $\jA_f$  in $\Spec(\jA)$ is an injective homomorphism of rings.  
The range of this map is a multiplicative class, hence is a base for a topology on $\Spec(\jA) = \bigcup \jA_f$, 
the so-called \unz{spectral topology}.  
In the spectral topology, $\Spec(\jA)$ is a locally compact, totally disconnected, Hausdorff space, i.e., is a Boolean space, 
the $\jA_f$ then being the ring of open and compact subsets of $\Spec(\jA)$.  
Because $\Spec(\jA)$ is compact iff $\jA$ admits a multiplicative identity, 
in the noncompact case, 
compactifying $\Spec(\jA)$ by the Alexandroff procedure is equivalent to formally passing from $\jA$ to the Boolean ring 
$\widehat{\jA}$ obtained by adjunction of a unit. 
If $\jA$ is infinite, then the weight of  $\Spec(\jA)$ is the cardinality of $\jA$; 
If $\jA$ is finite, then $\jA$ has $2^n$ elements and therefore $\Spec(\jA)$ is a discrete space with $n$ elements.
\\[-.25cm]

The set $\Spec(\jA)$, equipped with the spectral topology, is called the \unz{Stone space} of $\jA$.  
We shall denote it by the symbol $\ST(\jA)$.  
Evidently the Stone spaces of isomorphic Boolean rings are homeomorphic and conversely.
\\[-.25cm]

[Note: \ 
In reality, $\ST(\jA)$ comes supplied with a sheaf of rings.  
However, this additional structure, while fundamental from the scheme-theoretic point of view, 
plays no explicit role in the present considerations, 
the ring $\jS \jT (\jA)$ of open and compact subsets of $\ST(\jA)$ being its replacement.]
\\[-.25cm]

If $\XX$ is a Boolean space, then the Stone space of $\jA(\XX)$ can be identified with $\XX$.
\\[-.25cm]

\textbf{Examples}
\\[-.25cm]

(1) \quad
Let $\XX$ be an infinite set, equipped with the discrete topology.  
Let $\jS$ be the algebra consisting of the finite and cofinite subsets of $\XX$.  
Fix a point $\infty$ which is not in $\XX$ $-$then the map 
$\Phi : \jS \ra \jP(\XX \cup \{\infty\})$ defined by the rule
\[
\begin{cases}
\ \Phi(S) = S \hspace{1.75cm} \text{if $S$ is finite}
\\[11pt]
\ \Phi(S) = S \cup \{\infty\} \quad \text{if $S$ is infinite}
\end{cases}
\]
sets up an isomorphism between $\jS$ and an algebra $\jS_\infty$ of subsets of  $\XX \cup \{\infty\}$.  
Topologize $\XX \cup \{\infty\}$ by taking the class $\jS$ as a basis 
$-$then $\XX \cup \{\infty\}$ can be viewed as the Stone space of $\jS$ or still, 
the Stone space of $\jS$ is the Alexandroff compactification of $\XX$.
\\[-.25cm]

(2) \quad
Let $\XX$ be a set equipped with the discrete topology $-$then the Stone space of $\jP(\XX)$ has cardinality
\[
\begin{cases}
\ 2^{2^{\card(\XX)}} \hspace{0.6cm} \text{if $\XX$ is infinite}
\\[11pt]
\ \card(\XX) \quad \text{if $\XX$ is finite}
\end{cases}
\]
and can be identified with the Stone-C\v ech compactification of $\XX$.
\\[-.25cm]

(3) \quad
Let $\jA$ be a commutative ring with unit; 
let $\sI(\jA)$ be the set of idempotents of $\jA$ $-$then $\sI(\jA)$ is a Boolean algebra, the operations being
\[
\begin{cases}
\ f + g \equiv f + g - 2 \hsy f \hsy g
\\[8pt]
\ f \cdot g \equiv f \hsy g
\end{cases}
.
\]
\\[-1cm]

Suppose now that $\jA$ is regular in the sense of von Neumann, i.e., that every principal ideal is idempotent.  
Consider $\XX - \Spec(\jA)$ $-$then, topologized in the usual way, $\XX$ is a compact Boolean space and $\jA(\XX)$ 
is isomorphic to $\sI(\jA)$, implying, therefore, that $\XX$ can be regarded as the Stone space of $\sI(\jA)$.
\\[-.25cm]

The fact that $\jA$ is isomorphic to $\jS \jT (\jA)$ means that purely algebraic notions per $\jA$ can be 
reinterpreted vis-\`a-vis topological properties of the 
corresponding open and compact sets in $\ST(\jA)$.  
On the other hand, it is to be emphasized that this correspondence may break down when it becomes a question of infinite operations. 
For example, $\jS \jT (\jA)$ need not be a $\sigma$-ring even if $\jA$ is.
\\[-.25cm]

We shall write \bBR for the category whose objects are Boolean rings $\jA$, $\jB$, \ldots, and whose morphisms are the 
ring homomorphisms $\Phi : \jA \ra \jB$ such that $\Phi(\jA) \not\subset \fq$ $\forall \ \fq \in \Spec(\jB)$.  
Any morphism $\Phi : \jA \ra \jB$ of Boolean rings induces a continuous map $\Psi: \ST(\jB) \ra \ST(\jA)$ 
of corresponding Stone spaces.  This map is, moreover, proper.
\\[-.25cm]

[Note: \ 
We remark that if $\jA$ and $\jB$ are both Boolean algebras, then the condition that $\Phi : \jA \ra \jB$ be a morphism of 
Boolean rings is equivalent to the requirement that $\Phi : \jA \ra \jB$ be a homomorphism of rings taking the 
multiplicative  identity of $\jA$ to the multiplicative  identity of $\jB$.]
\\[-.25cm]

We shall write \bBS for the category whose objects are Boolean spaces $\XX$, $\YY$, \ldots, and whose morphisms are the 
proper continuous maps $\Psi : \XX \ra \YY$. 
Any morphism $\Psi : \XX \ra \YY$ of Boolean spaces induces a morphism $\Phi : \jA(\YY) \ra \jA(\XX)$ of Boolean rings.
\\[-.25cm]

\textbf{Example} \ 
Let  $\jA$ be a Boolean ring without a multiplicative  identity, 
$\widehat{\jA}$ the  Boolean ring obtained by adjunction of a unit $-$then the canonical injection 
$\jA \hookrightarrow \widehat{\jA}$ is not a morphism in \bBR.  
Put $\XX = \Spec(\jA)$, $\widehat{\XX} = \Spec(\widehat{\jA})$ $-$then the canonical injection 
$\XX \hookrightarrow \widehat{\XX}$ is not a morphism in \bBS.
\\[-.25cm]


These considerations can best be interpreted functorially.
\\[-.25cm]

(1) \quad
$\exists$ a contravariant functor 
\[
\jF_\bRS \hsy : \hsy  \bBR \ra \bBS. 
\]
Here
\[
\begin{cases}
\ \jA \ra \ST(\jA)
\\[11pt]
\ \Phi \in \Hom (\jA, \jB) \ra \Psi \in \Hom(\ST(\jB), \ST(\jA))
\end{cases}
.
\]
\\[-.25cm]

(2) \quad
$\exists$ a contravariant functor 
\[
\jF_\bSR \hsy : \hsy  \bBS \ra \bBR. 
\]
Here
\[
\begin{cases}
\ \XX \ra \jA(\XX)
\\[11pt]
\ \Phi \in \Hom (\XX, \YY) \ra \Psi \in \Hom(\jA(\YY), \jA(\XX))
\end{cases}
.
\]
\\[-.25cm]

Call $1_\bBR$, $1_\bBS$ the identity functors in \bBR, \bBS $-$then it is easy to check that 
$\jF_\bSR \circ \jF_\bRS$ is isomorphic to $1_\bBR$
and 
$\jF_\bRS \circ \jF_\bSR$ is isomorphic to $1_\bBS$.  
The categories \bBR and \bBS are therefore dual.  
\\[-.25cm]

\unz{Ref} \quad 
The results discussed above are surveyed in 
M. Stone\footnote[2]{\vspace{.11 cm}\textit{Bull. Amer. Math. Soc.}, \textbf{44} (1938), 807-816.}
the complete account being given in 
M. Stone\footnote[3]{\vspace{.11 cm}\textit{Trans. Amer. Math. Soc.}, \textbf{40} (1936), 37-111.}
and 
M. Stone\footnote[4]{\vspace{.11 cm}\textit{Trans. Amer. Math. Soc.}, \textbf{41} (1937), 375-481.}


\chapter{
$\boldsymbol{\S}$\textbf{5}.\quad  Products and Sums}
\setlength\parindent{2em}
\setcounter{theoremn}{0}
\renewcommand{\thepage}{\S5-\arabic{page}}



\qquad
Let $\XX$ and $\YY$ be nonempty sets $-$then by
\[
\begin{cases}
\ \pi_\XX : \XX \times \YY \ra \XX
\\[8pt]
\ \pi_\YY : \XX \times \YY \ra \YY
\end{cases}
,
\]
we shall understand the \un{projections} of $\XX \times \YY$ onto $\XX$ and $\YY$, respectively.  
Given a subset $E$ of $\XX \times \YY$ and points $x \in \XX$, $y \in \YY$, put
\[
\begin{cases}
\ E_x \ = \ \pi_\YY^{} \hsy [\pi_\XX^{-1} (x) \cap E]
\\[8pt]
\ E^y \ = \ \pi_\XX^{} \hsy [\pi_\YY^{-1} (y) \cap E]
\end{cases}
,
\]
the
\[
\begin{cases}
\ \text{vertical}
\\[4pt]
\ \text{horizontal}
\end{cases}
\]
\un{sections} of $\E$ over 
$
\begin{cases}
\ x
\\
\ y
\end{cases}
.
$
It is easy to check that

\[
\begin{cases}
\ \left(\bigcup \ E_i \right)_x \ = \ \bigcup\ \left(E_i \right)_x
\\
\ \left(\bigcup \ E_i \right)^y \ = \ \bigcup\ \left(E_i \right)^y
\end{cases}
\qquad 
\begin{cases}
\ \left(\bigcap \ E_i \right)_x \ = \ \bigcap\ \left(E_i \right)_x
\\
\ \left(\bigcap \ E_i \right)^y \ = \ \bigcap\ \left(E_i \right)^y
\end{cases}
\]
\\[-1cm]

\[
\begin{cases}
\ (\XX \times \YY - E)_x \ = \ \YY - E_x
\\[8pt]
\ (\XX \times \YY - E)^y \ = \ \XX - E^y
\end{cases}
.
\]
\\[-1cm]

Let $S \subset \XX$, $T \subset \YY$ $-$then the \un{rectangle} $R$ determined by $S$ and $T$ is the Cartesian product 
$S \times T \hookrightarrow \XX \times \YY$, $S$ and $T$ being its \un{sides}.
One has 
$\chisubR = \chisubS \cdot \chisubT$.  
It is clear that a rectangle is empty iff one of its sides is empty.  
Furthermore, if 
$R_1 = S_1 \times T_1$ and 
$R_2 = S_2 \times T_2$ are nonempty rectangles, then $R_1 \subset R_2$ iff $S_1 \subset S_2$ and $T_1 \subset T_2$.  
Consequently, two nonempty rectangles are equal iff both of their sides are equal. 
\\

There are some simple identities governing the manipulation of rectangles  
which we had best record explicitly as they will be used 
tacitly in what follows.
\\[-.25cm]

\qquad $\bigcup$: 
\[
\left(
\bigcup\limits_{i \in I} \ S_i
\right) 
\times
\left(
\bigcup\limits_{j \in J} \ T_j
\right) 
\ = \ 
\bigcup\limits_{(i, j) \in I \times J} \ S_i \times T_j
\]

[In particular: 
\[
(S_1 \cup S_2) \times (T_1 \cup T_2) 
\ = \ 
(S_1 \times T_1) \cup (S_1 \times T_2) \cup (S_2 \times T_1) \cup (S_2 \times T_2).]
\]
\\[-1cm]

\qquad $\bigcap$: 
\[
\left(
\bigcap\limits_{i \in I} \ S_i
\right) 
\times
\left(
\bigcap\limits_{j \in J} \ T_j
\right) 
\ = \ 
\bigcap\limits_{(i, j) \in I \times J} \ S_i \times T_j
\]

[In particular: 
\[
(S_1 \cap S_2) \times (T_1 \cap T_2) 
\ = \ 
(S_1 \times T_1) \cap (S_2 \times T_2).]
\]


\qquad $-$: 
\[
\begin{cases}
\ (S_1 - S_2) \times T 
\ = \ 
(S_1 \times T) - (S_2 \times T)
\\[8pt]
\  S \times (T_1 - T_2)
\ = \ 
(S \times T_1) - (S \times T_2)
\end{cases}
\]

\[
(S_1 \times T_1) - (S_2 \times T_2)
\begin{cases}
\ [(S_1 - S_2) \times (T_1 \cap T_2)] 
\hsx \cup \hsx
\ [S_1 \times (T_1 - T_2)]
\\[8pt]
\ [(S_1 - S_2) \times T_1]
\hsx \cup \hsx
[(S_1 \cap S_2) \times (T_1 - T_2)]
\end{cases}
\]

[In particular: \ 
The difference of two rectangles can be written as the disjoint union of two other rectangles.]
\\

Consider now the natural map 
\[
\jP (\XX) \times \jP (\YY)
\ra 
\jP (\XX \times \YY),
\]
namely the rule assigning to each pair $(S, T)$ the rectangle $R = S \times T$.  
As this map is evidently bilinear, it must factor canonically

\[
\begin{tikzcd}[sep=large]
\jP (\XX) \times \jP (\YY)
\ar{r}
\ar{rd}
&\jP (\XX)  \otimes  \jP (\YY)
\ar{d}
\\
&\jP (\XX \times \YY)
\end{tikzcd}
.
\]
Here, the tensor product is taken over $\Z$ or still, since it amounts to the same, over $\Z_2$.  
After a moments reflection, the reader will agree that the vertical arrow is actually an injection, 
its range being the class consisting of those sets in $\XX \times \YY$ which can be written as a finite union of rectangles.  
Because the image of the ring $\jP (\XX)  \otimes  \jP (\YY)$ contains all singletons, the associated completion is 
$\jP (\XX \times \YY)$.
\\[-.25cm]

To illustrate these remarks, suppose that $\jS$ is a subring of $\jP (\XX)$ and that $\jT$ is a subring of $\jP(\YY)$ 
$-$then since everything in sight is flat, 
\[
\jS  \otimes  \jT
\hookrightarrow
\jP (\XX)  \otimes  \jP (\YY).
\]
Accordingly, $\jS  \otimes  \jT$ may be regarded as the class of all subsets of $\XX \times \YY$ of the form
\[
\bigcup\limits_{i = 1}^m \ (S_i \times T_i) 
\qquad (S_i \in \jS, \hsx T_i \in \jT),
\]
it not being restrictive to suppose that any such union is even disjoint. 
\\[-.25cm]

Generally, if $\jS$ is a nonempty class of subsets of $\XX$ and if $\jT$ is a nonempty class of subsets of $\YY$, 
then we shall write $\jS \ \boxtimesdmc \ \jT$ for the class of rectangles $R = S \times T$ $(S \in \jS, \hsx T \in \jT)$.  
In other words $\jS \ \boxtimesdmc \ \jT$ is simply the image of $\jS \times \jT$ under the natural map
\[
\jP (\XX) \times \jP (\YY) \ra 
\jP (\XX \times \YY).
\]
\\[-1cm]


Observe that: 
\\[-.5cm]

(1) \quad 
If $\jS$ and $\jT$ are multiplicative classes, then $\jS \ \boxtimesdmc \ \jT$ is a multiplicative class. 
\\[-.25cm]

(2) \quad 
If $\jS$ and $\jT$ are additive classes, then $\jS \ \boxtimesdmc \ \jT$ need not be an additive class. 
\\

{\small\bf Lemma 1} \ 
Let $\jS$ and $\jT$ be semirings $-$then $\jS \ \boxtimesdmc \ \jT$ is a semiring.
\\[-.5cm]

[We omit the verification.]
\\

Suppose that $\jS$ and $\jT$ are rings $-$then $\jS \ \boxtimesdmc \ \jT$ is a semiring but rarely a ring.  
However, if we apply the Kolmogoroff procedure to $\jS \ \boxtimesdmc \ \jT$ (cf. \S4), the result will be a ring, viz. 
$\jS  \otimes  \jT$.  
\\[-.5cm]

Suppose that $\jS$ and $\jT$ are $\sigma$-rings 
$-$then $\jS$ is necessarily closed in $\jP (\XX)_S$ and 
$\jT$ is necessarily closed in $\jP (\YY)_S$.  
Nevertheless, 
$\jS  \otimes  \jT$
is not necessarily closed in $\jP (\XX \times \YY)_S$, hence ordinarily fails to be a $\sigma$-ring.
\\

{\small\bf Example} \ 
Take $\XX = \YY$ of cardinality $\aleph_0$ and let $\jS = \jT$ be the class of all subsets of cardinality $\leq \aleph_0$ 
$-$then the diagonal $D$ belongs to the closure of 
$\jS  \otimes  \jS$ in 
$\jP (\XX \times \XX)_S$ 
but is certainly not in 
$\jS  \otimes  \jS$ itself.
\\

If $\jS$ and $\jT$ are rings, then in what follows we shall write 
$\jS \ \ov{\otimes} \ \jT$ 
for the closure of $\jS  \otimes  \jT$ in 
$\jP (\XX \times \YY)_S$.  
Needless to say, 
$\jS \ \ov{\otimes} \  \jT$ 
is a $\sigma$-ring; of course, 
$\jS  \otimes  \jT \neq \jS \ \ov{\otimes} \  \jT$ 
in general, even if both $\jS$ and $\jT$ are $\sigma$-rings (cf. supra).
\\

{\small\bf Lemma 2} \ 
Let $\jS$ and $\jT$ be $\sigma$-rings; 
let $E \in \jS \ \ov{\otimes} \  \jT$
$-$then 
\[
\begin{cases}
\ E_x \in \jT \hspace{0.5cm} \forall \ x \in \XX
\\[8pt]
\ E^y \in \jS \hspace{0.5cm} \forall \ y \in \YY
\end{cases}
.
\]
\\[-1cm]

[One need only note that the class of all subsets of $\XX \times \YY$ with the stated property contains 
$\jS  \otimes  \jT$ 
and is closed in 
$\jP (\XX \times \YY)_S$
.]
\\

Here is a corollary.  
Let 
$R = S \times T$ 
be a nonempty rectangle in 
$\XX \times \YY$ 
$-$then 
$R \in \jS \ \ov{\otimes} \  \jT$ 
iff $S \in \jS$ and $T \in \jT$.
\\[-.5cm]

[Note:  \ 
The converse to Lemma 2 is false as can be seen by a slight alteration of the preceeding example, namely this time take 
$\XX = \YY$ 
of cardinality $> \aleph_0$ 
and, with 
$\jS = \jT$ as there, consider the diagonal $D$.]
\\

{\small\bf Example} \ 
Take $\XX = \YY$.
Consider the following question: \ 
Is 
$\jP (\XX)  \otimes  \jP (\XX)$ 
dense in 
$\jP (\XX \times \XX)_S$
?  
The answer depends on the cardinality of $\XX$.
\\[-.25cm]

\qquad (1) \quad 
Suppose that $\card(\XX) \geq \fc$ $-$then 
$\jP (\XX)  \otimes  \jP (\XX)$ 
is not dense in 
$\jP (\XX \times \YY)_s$.  
\\[-.25cm]


\un{Re (1)} \ 
Proceed by contradiction $-$then of necessity, the diagonal $D$ would belong to 
$\jP (\XX) \ \ov{\otimes} \  \jP (\XX)$.  
Therefore, in view of a simple property of the sequential modification (cf. \S1), one could find a ring $\jS$ in $\XX$ of 
cardinality $\leq \aleph_0$ such that $D$ actually belongs to 
$\jS \ \ov{\otimes} \  \jS$.  
Denote by 
$\sigma\dash\jRin (\jS)$ the closure of $\jS$ in $\jP (\XX)_s$ $-$then, thanks to Lemma 2,
\[
\forall \ x \in \XX \hsy : \hsy \{x\} \in \sigma\dash\jRin (\jS).
\]
Let $S_1, \hsy S_2, \ldots$ be an enumeration of the elements of $\jS$ $-$then we claim that the characteristic function 
$f : \XX \ra \Cx$ of the $S_i$ (cf. Prob. IV (\S1)), 
\[
f(x) 
\ = \ 
2 \hsx \cdot \hsx 
\sum\limits_{i = 1}^\infty \ 
\chisubSi (x) / 3^i 
\qquad (x \in \XX),
\]
is one-to-one, hence that $\card(\XX) \leq \fc$.  
Indeed, if $f(x) = f(y)$, then $\forall \ i$, $x \in S_i$ iff $y \in S_i$.  
But the class of all subsets $S \subset \XX$ such that either $\{x, y\} \subset S$ or $\{x, y\} \cap S = \emptyset$ 
is a $\sigma$-ring containing $\jS$, thus contains the singletons and so $x = y$, as claimed.
\\[-.25cm]

[Note: \ 
For a somewhat different approach to this result, see Exer. 21 (\S6).]
\\[-.25cm]

\un{Re (2)} \ 
There is no loss of generality in taking $\XX$ to be a subset of $\R$.  
If $\card(\XX) \leq \aleph_0$, then the assertion is clear.  
We shall therefore suppose that 
$\card(\XX) = \aleph_1$.  
For the purposes at hand, let us agree that a \un{curve} in $\XX \times \XX$ is simply any set of the form
\[
\{(x , f(x)) \hsy : \hsy x \in \dom (f)\}, 
\quad 
\{(x , g(x)) \hsy : \hsy x \in \dom (g)\}
\]
where
\[
\dom (f) \subset \XX, 
\quad
\dom (g) \subset \XX, 
\]
and 
$f \hsy : \hsy \dom (f)  \ra \XX$, 
$g \hsy : \hsy \dom (g)  \ra \XX$ 
are functions.  
Every curve is in the closure of 
$\jP (\XX)  \otimes \hsx   \jP (\XX)$
in 
$\jP (\XX \times \XX)_s$.  
To see this, note that
\[
\{(x , f(x)) \hsy : \hsy x \in \dom (f)\}
\ = \ 
\bigcap\limits_{m = 1}^\infty \ E_m
\]
where
\[
E_m
\ = \ 
\bigcup\limits_{i = -\infty}^\infty \ E_{i \hsy m}
\]
with
\[
E_{i \hsy m}
\ = \ 
\left\{
x \in \dom (f) \hsy : \hsy \frac{i}{m} \leq f(x) < \frac{i+1}{m}
\right\} 
\times \ 
\XX 
\ \cap \ 
\left[
\frac{i}{m}, \frac{i+1}{m}
\right[
\]
and similarly for $g$.  
To prove (2), therefore, it need only be shown that 
$\XX \times \XX$
can be written as a countable union of curves.  
To this end, well-order 
$\XX \hsy : \hsy \{x_\alpha \hsy : \hsy \alpha < \Omega\}$.  
Divide 
$\XX \times \XX$ 
into complementary sets $E$ and $F$ by the definitions

\[
\begin{cases}
\ E =  \{(x_\alpha, x_\beta) \hsy : \hsy \beta < \alpha\}
\\[8pt]
\ F =  \{(x_\alpha, x_\beta) \hsy : \hsy \alpha < \beta\}
\end{cases}
.
\]
It is clear that the vertical sections of $E$ are finite or countable, as are the horizontal sections of $F$.  
For each $x \in \XX$, arrange $E_x$ into a sequence $\{x_n\}$, it being understood that the sequence is to be completed 
in an arbitrary way if it is finite to begin with.  
Define now functions 
$f_n : \XX \ra \XX$ by the prescription $f_n(x) = x_n$.  
Analogous considerations apply to the horizontal sections $F^x$ of $F$ leading to functions 
$g_n : \XX \ra \XX$.  
Taken together, the curves
\[
\{(x , f_n(x)) \hsy : \hsy x \in \XX\}, 
\quad 
\{(x , g_n(x)) \hsy : \hsy x \in \XX\} 
\]
cover $\XX \times \XX$.
\\[-.5cm]

[Note: \ 
The last part of the preceeding argument is virtually the same as that needed in the first part of Prob. IV (\S2).]
\\

\un{Re (3)} \ 
On the basis of (2), this is immediate.
\\

[Note: \ 
Actually, one can get away with less here in that Martin's axiom alone suffices to force the conclusion if 
$\aleph_1 < \card(\XX) \leq \fc$; 
cf. Kunen, 
\un{Inaccessibility} \un{Properties of Cardinals}, Ph. D. Thesis, Stanford University, 1968.]
\\


Partitions in $\XX$ and $\YY$ are closed related to partitions in $\XX \times \YY$ and vice versa. 
\\[-.25cm]

{\small\bf Lemma 3} \ 
Let $R = S \times T$ be a nonempty rectangle; 
let 
$\{R_k = S_k \times T_k\}$
be a class of nonempty rectangles $-$then the $R_k$ partition $R$ iff
\\[-.25cm]

\qquad (i) \quad 
$R = \bigcup \ R_k$; 
\\[-.5cm]

\qquad (ii) \quad 
$S = \bigcup \ S_k$, 
\quad 
$T = \bigcup \ T_k$; 
\\[-.5cm]

\qquad (iii) \quad 
$\forall \ k \neq \ell$
\\[-.5cm]

\[
\begin{cases}
\ S_k \cap S_\ell \neq \emptyset 
\implies T_k \cap T_\ell = \emptyset
\\[4pt]
\quad \text{or}
\\[4pt]
\ T_k \cap T_\ell \neq \emptyset 
\implies S_k \cap S_\ell = \emptyset
\end{cases}
.
\hspace{2cm}
\]
\\[-1cm]

[We omit the verification.]
\\

Let $R$ be a nonempty rectangle $-$then a partition 
$P(R) = \{R_k \hsy : \hsy k \in K\}$ of $R$ by rectangles is said to be a  \un{network} on $R$ if

\[
\begin{cases}
\ \text{the $\pi_\XX (R_k)$ partition $\pi_\XX (R)$}
\\[4pt]
\quad \text{and}
\\[4pt]
\ \text{the $\pi_\YY (R_k)$ partition $\pi_\YY (R)$}
\end{cases}
.
\]
\\[-1cm]

[Note: 
Here we are admitting a small solecism in that 
repetitions may, of course, be present in the classes $\pi_\XX (R_k)$, $\pi_\YY (R_k)$.]
\\

{\small\bf Lemma 4} \ 
Let 
$\jS \subset \jP (\XX)$, 
$\jT \subset \jP (\YY)$
be multiplicative classes; 
let 
$R = S \times T \in \jS \ \boxtimesdmc \  \jT$ 
be a nonempty rectangle.  
Suppose that $P(R)$ is a finite 
$\jS \ \boxtimesdmc \  \jT$-partition of $R$ $-$then there exists a partition in 
$\Par_{\jS \ \boxtimesdmc \  \jT} (R)$ which refines $P(R)$ and is a network on $R$.
\\[-.25cm]

\un{Proof} \ 
It can be assumed that $P(R)$ is not a network on $R$.  
Denoting the components of $P(R)$ by $R_k$, let $S_k = \pi_\XX (R_k)$, $T_k = \pi_\YY (R_k)$ $-$then 
$S = \bigcup \ S_k$, $T = \bigcup \ T_k$.  
Consider the $S_k$.  
Define an equivalence relation on $\XX$ by stipulating that $x_1$ be equivalent to $x_2$ iff
\[
\forall \ y \in \Y, \ (x_1, y) \sim (x_2, y), 
\]
the latter equivalence being that corresponding to $P(R)$.  
Given $x \in \XX$, the equivalence class $[x]$ determined by $x$ is simply the intersection of the $S_k$ containing $x$.  
All told, therefore, this procedure produces a finite $\jS$-partition 
$P(S) = \{S_i \hsy : \hsy i \in I\}$ of $S$.  
Work with the $T_k$ in an analogous fashion to produce a finite $\jT$-partition 
$P(T) = \{T_j \hsy : \hsy j \in J\}$ of $T$.    
The $S_i \times T_j$ then
constitute a finite 
$\jS \ \boxtimesdmc \  \jT$ 
partition of $R$, refining $P(R)$ and forming a network on $R$.
\\

Retaining the notation from Lemma 4, suppose that $P(R)$ is a countable 
$\jS \ \boxtimesdmc \  \jT$-partition 
of $R$.  
We then ask: 
Does there exist a partition in $\sigma$-$\Par_{\jS \ \boxtimesdmc \  \jT} (R)$ which refines $P(R)$ and is a network on $R$? 
Unfortunately, even after imposing about as much additional structure on $\jS$ and $\jT$ as can be reasonably expected, 
the answer will in general be negative.
\\

{\small\bf Examples} \ 
\\[-.5cm]

\qquad (1) \quad 
Take $\XX = [-1, 1[$, $\YY = [0, +\infty[$.  
Let $\jS$ be the class consisting of all left closed and right open subintervals of $\XX$; 
let $\jT = \jP (\YY)$ $-$then $\jS$ is a semiring and $\jT$ is a complete ring. 
Consider the countable 
$\jS \ \boxtimesdmc \  \jT$ 
of 
$\XX \times \YY$ 
by the rectangles 
\\

$
[-1, 1[ \times [0, 1[ \quad
\begin{cases}
\left[-1, -\frac{1}{n}\right] \times [n - 1, n[
\\[8pt]
\left[-\frac{1}{n}, \frac{1}{n}\right] \times [n - 1, n[
\\[8pt]
\left[\frac{1}{n}, 1\right] \times [n - 1, n[
\end{cases}
\qquad (n > 1).
$
\\[0.25cm]

\noindent
Because $0 \in \left[-\frac{1}{n}, \frac{1}{n}\right]$ $\forall \ n$, it is impossible to find a countable 
$\jS \ \boxtimesdmc \ \jT$ 
$-$network on $\XX \times \YY$ which refines this partition.
\\

\qquad (2) \quad 
Take $\XX = ]0, 1[$, $\YY = ]0, 1[ \hsx \cap \Q$.    
Let 
$\jS = \jP (\XX)$, 
$\jT = \jP (\YY)$ 
$-$then both $\jS$ and $\jT$ are complete rings.  
Consider the countable 
$\jS \ \boxtimesdmc \ \jT$-partition of 
$\XX \times \YY$ 
by the rectangles 

\[
\begin{cases}
\ ]0, q[ \hsx \times \{q\}
\\[8pt]
\ [q, 1[ \hsx \times \{q\}
\end{cases}
\qquad (0 < q < 1, q \in \Q).
\]
Suppose that the $S_i \times T_j$ $(i \in I, \hsx j \in J)$ refine this partition and form a network on $\XX \times \YY$ 
$-$then, of necessity,

\[
\begin{cases}
\ \card(I) > \aleph_0
\\[8pt]
\ \card(J) = \aleph_0
\end{cases}
,
\]
so $I \times J$ must be uncountable.
\\

Up until this point, the discussion has dealt exclusively with products involving two factors.  
The extension of the theory to $n > 2$ $(n \in \N)$ factors is purely formal, hence need not be considered in detail.  
We remark only that tacitly one makes throughout the usual conventions as regards the associativity of the relevant operations. 
\\

The situation for products involving an arbitrary number of factors is only slightly more complicated, 
it being a matter of setting up the definitions in a succinct fashion.  
Let, then, 
$\{X_i \hsy : \hsy i \in I\}$ 
be a class of nonempty subsets $X_i$ indexed by an infinite set $I$ $-$then we shall agree that a \un{rectangle} in 
$\prod \ X_i$ 
is a set of the form 
$\prod \ S_i$, 
where $S_i \subset X_i$ $\forall \ i$ and $S_i = X_i$ for all 
but a finite set of $i$.  
If 
$S = \prod \ S_i$ 
and if 
\[
\begin{cases}
\ S^\prime = \prod \ S_i^\prime 
\\[8pt]
\ S^{\prime\prime} = \prod \ S_i^{\prime\prime} 
\end{cases}
\]
are nonempty rectangles, then $S = S^\prime \cup S^{\prime\prime}$ with 
$S^\prime \cap S^{\prime\prime} = \emptyset$ 
iff there exists a unique index $i_0$ such that 
\[
\begin{cases}
\ i \neq i_0 \implies S_i = S_i^\prime = S_i^{\prime\prime}
\\[8pt]
\ i = i_0 \implies S_i = S_i^\prime \cup S_i^{\prime\prime},
\quad 
S_i^\prime \cap S_i^{\prime\prime} = \emptyset
\end{cases}.
\]
\\[-1cm]

Consider now the tensor product 
$ \otimes  \ \jP(X_i)$ 
$-$then, $\forall \ i$, $\exists$ a canonical homomorphism 
\[
\iota_i : \jP (X_i) \ra  \otimes  \ \jP(X_i),
\]
namely the rule which assigns to each $S_i \subset X_i$ the tensor whose $i^\nth$ entry is $S_i$ 
and whose $j^\nth$ entry is $X_j$ $(j \neq i)$.  
The subalgebra of 
$ \otimes  \ \jP(X_i)$ 
generated by the 
$\iota_i ( \otimes  \ \jP(X_i))$ 
is composed of all finite sums of elements of the form 
$ \otimes  \ S_i$, 
where 
$S_i = X_i$ 
except for a finite number of indices.  
Algebraists customarily refer to this subalgebra of 
$ \otimes  \ \jP(X_i)$  
as the \un{tensor product of the algebras} 
$\jP(X_i)$.  
We shall denote it by 
$\otimes^* \jP (X_i)$.  
Since the index set $I$ is infinite, it differs in general from 
$ \otimes  \ \jP(X_i)$.  
\\

[Note: \ 
Consideration of $\oplus^* \hsx \jP (X_i)$ is, of course, necessary 
from the categorical point of view.]
\\

Denote by 
$\prod^* \ \jP (X_i)$ 
that subset of 
$\prod \ \jP (X_i)$
consisting of the $(S_i)$ such that $S_i = X_i$ for all but a finite set of $i$.
There is a commutative triangle

\[
\begin{tikzcd}[sep=large]
\prod^* \ \jP (X_i)
\ar{r}
\ar{rd}
&\oplus^* \hsx \jP (X_i)
\ar{d}
\\
&\jP \left(\prod \ X_i \right)
\end{tikzcd}
.
\]

\noindent
The vertical arrow is an injection, its range being the class of those sets in $\prod \ X_i$ which can be written as a 
finite union of rectangles.  

Finally, we come to the one big difference between infinite as opposed to finite products, namely this: \ 
It is necessary to consider algebras $\jS_i \subset \jP (X_i)$ rather than just rings.  
The reason is easy enough to see.  
Indeed, if we proceed as above to form 
$\oplus^* \hsx \jS_i$, 
then each of the $\jS_i$'s must at least be rings with unit and to ensure compatibility, 
it is best to assume that they are actually algebras.  
Under these circumstances, 
\[
\oplus^* \hsx \jS_i 
\hookrightarrow 
\oplus^* \hsx \jP (X_i)
\]
meaning, therefore, that $\oplus^* \hsx \jS_i$ can be thought of as sitting
inside 
$\jP \left(\prod \ X_i \right)$, 
the characterization reading as in the finite case, i.e., the class of all finite disjoint unions of rectangles 
$\prod \ S_i$, where $S_i \in \jS_i$ $\forall \ i$.  
This being so, we shall then write 
$\ovs{\otimes}^* \hsx \jS_i$ for the closure of 
$\oplus^* \hsx \jS_i$ 
in 
$\jP \left(\prod \ X_i \right)_\tS$.  
Evidently, 
$\ovs{\otimes}^* \hsx \jS_i$ is a $\sigma$-algebra.  
\\

Keeping to the preceding notation, put 
$\XX = \prod \ X_i$, 
$\jS = \ov{\otimes} \hsx \jS_i$.  
Let $I = I_1 \cup I_2$ be a partition of $I$.
Let 
$X_1 = \prod\limits_{I_1} \ X_i$, 
$X_2 = \prod\limits_{I_2} \ X_i$; 
let 
$\jS_1 = \underset{I_1}{\ovs{\otimes}^*} \hsx \jS_i$, 
$\jS_2 = \underset{I_2}{\ovs{\otimes}^*} \hsx \jS_i$ 
$-$then $\XX$ may be indentified with 
$X_1 \times X_2$ 
and, when this is done, we have 
$\jS = \jS_1 \ov{\otimes} \jS_2$.  
Therefore, in a certain sense, we are right back at the beginning.
\\

{\small\bf Example} \ 
Let 
$\{X_i \hsy : \hsy i \in I\}$ 
be a class of compact Hausdorff spaces $X_i$ indexed by an infinite set $I$.  
Take for $\jS_i$ the algebra of open and compact subsets of $X_i$ $-$then 
$\otimes^* \ \jS_i$ 
is the algebra of open and compact subsets of $\prod \ X_i$.  
\\[-.25cm]

[Let us consider an important special case.  
Equip $\{0, 1\}$ with the discrete topology.  
Given any $i = 1, 2, \ldots,$ put $X_i = \{0, 1\}$ $-$then, in the product topology, 
$\ds 2^\N \ =  \ \prod\ X_i$ is a compact, totally disconnected, Hausdorff space of weight $\aleph_0$, 
the so-called \un{Cantor space}.  
Of course, the terminology arises from the fact that $2^\N$ is homeomorphic to $\tC$, viz. (cf. Prob. IV (\S1)): \ 
\[
(f \in 2^\N \mapsto 
\left(
2 \cdot 
\sum\limits_{i = 1}^\infty \ 
\frac{f(i)}{3^i} \in \tC
\right).
\]
Let $\jS_i$ be the algebra of all subsets of $X_i$ $-$then 
$\otimes^* \ \jS_i$ 
is the algebra of
open and compact subsets of 
$2^\N$ 
and 
$\ovs{\otimes}^* \ \jS_i$ 
is the $\sigma$-algebra of Borel subsets of 
$2^\N$  
(cf. \S6).]
\\

Let 
$\{X_i \hsy : \hsy i \in I\}$ 
be a class of nonempty sets $X_i$ indexed by a nonempty set $I$ (finite or infinite), the $X_i$ being, in addition, pairwise disjoint.  
Write
$\oplus \ \jP(X_i)$ 
for the direct sum of the $\jP(X_i)$.  
Suppose that $\forall \ i$, $\jS_i$ is a ring in $X_i$ $-$then the direct sum 
$\oplus \ \jS_i$ 
of the $\jS_i$ 
is a subring of 
$\oplus \ \jP(X_i)$ 
.  
The elements of 
$\oplus \ \jS_i$ 
may be viewed as those subsets $S$ of $\bigcup \ X_i$ with the property that 
$S \cap X_i \in \jS_i$ for all $i$, or still, as the class of all unions $\bigcap \ S_i$, where $S_i \in \jS_i$ $(\forall \ i)$.  
If each of the $\jS_i$ is a $\sigma$-ring, then so is 
$\oplus \ \jS_i$.
\\

[Note: \ 
If the $X_i$ are not initially pairwise disjoint, then this may always be arranged by looking at the $X_i \times \{i\}$.]
\\

{\small\bf Example} \ 
Let $\jS$ be a $\sigma$-ring in $\XX$.  
Fix a countable partition 
$P(\XX) = \{X_i \hsy : \hsy i \in I\}$ 
of $\XX$, where $X_i \in \jS$ $\forall \ i$.  
Put 
$\jS_i = \tr_{X_i} (\jS)$ $-$then 
$\jS = \oplus \ \jS_i$.
\\

\[
\textbf{\un{Notes and Remarks}}
\]

Just who was the first to consider products in abstraco is not completely clear.  
The following papers are relevant: 
H. Hahn\footnote[1]{\vspace{.11 cm}\textit{Ann. Scuola Norm. Sup. Pisa}, \textbf{2} (1933), 429-452.}. 
F. Maeda\footnote[2]{\vspace{.11 cm}\textit{T\^ohoku Math. J.}, \textbf{37} (1933), 446-453.}. 
Z Lomnicki and 
S. Ulam\footnote[3]{\vspace{.11 cm}\textit{Fund. Math.}, \textbf{23} (1934), 27-36.}. 
J. Ridder\footnote[4]{\vspace{.11 cm}\textit{Fund. Math.}, \textbf{24} (1934), 72-117.}. 
W. Feller\footnote[5]{\vspace{.11 cm}\textit{Bull. Int. Acad Youg.}, \textbf{28} (1934), 30-45.}. 
B Jessen\footnote[6]{\vspace{.11 cm}\textit{Acta. Math.}, \textbf{63} (1934), 249-323.}. 
The question of the density of $\jP(\XX)  \otimes  \ \jP(\XX)$ in $\jP(\XX \times \XX)_\ts$ is an old problem of Ulam and has been considered 
by a number of authors; 
cf. 
B. Rao\footnote[7]{\vspace{.11 cm}\textit{Acta. Math. Acad. Sci. Hungar.}, \textbf{22} (1971), 197-198.}. 
Lemma 4 is a variation on a well known theme; 
it is explicitly stated and proved in 
D. Goguadze\footnote[8]{\vspace{.11 cm}
{Kolmogoroff Integrals and Some of their Applications},
{\fontencoding{OT2}\selectfont
[Ob
Integralah
Kolmogorova
I 
Ih  
Nekotoryh
Prilozheniyah],
Meciereba
Tbilisi
}
, (1979), 152-153.}. 
This author goes on to claim (statement 13.8, p. 154) that if $\jS$ and $\jT$ are semirings, then Lemma 4 is true when ``finite'' is replaced by ``countable''.  
As we have seen in the text, this is false.  
It may have occurred to the reader that the language of category theory might be helpful at certain points in this \S; 
some comments in this direction may be found in 
L. Auslander and 
C. Moore\footnote[9]{\vspace{.11 cm}\textit{Mem. Amer. Math. Soc}, \textbf{62} (1966), 1-199.}. 


\chapter{
$\boldsymbol{\S}$\textbf{5}.\quad  Exercises}
\setlength\parindent{2em}
\setcounter{theoremn}{0}
\renewcommand{\thepage}{\S5-\text{E-}\arabic{page}}



\qquad
(1) \ 
True or False?
\\[-.5cm]

\qquad (a) \quad 
There exists a nonempty set $E$ such that $E \times E \subset E$.
\\[-.5cm]

\qquad (b) \quad 
There exists a nonempty set $E$ such that $E \subset E \times E$.
\\[-.25cm]

(2) \ 
Discuss the continuity of the natural maps

\[
\begin{cases}
\ \jP(\XX) \times \jP(\YY)  \ra \jP(\XX \times \YY) 
\\[8pt]
\ \jP(\XX)_\ts \times \jP(\YY)_\ts  \ra \jP(\XX \times \YY)_\ts 
\end{cases}
.
\]

(3)
Let $\jKol (?)$ be the ring obtained from the semiring ? via the Kolmogoroff procedure (cf. \S4).
\\[-.5cm]

\qquad
True or False?
If $\jS$ is a semiring in $\XX$ and if $\jT$ is a semiring in $\YY$ (so that $\jS \hsx \boxtimes \hsx \jT$ is a semiring in $\XX \times \YY$), then
\[
\jKol (\jS) \hsy \otimes \hsy \jKol (\jT)
\ = \ 
\jKol (\jS \hsx \boxtimes \hsx \jT).
\]
\\[-1cm]

(4) \ 
Let $\jS$ and $\jT$ be $\sigma$-rings; 
let $E \in \jS \ \ov{\otimes} \ \jT$ $-$then there exists $S \in \jS$, $T \in \jT$ such that $E \subset S \times T$.  
\\[-.25cm]

(5)  \ 
Let $\XX$ and $\YY$ be nonempty sets $-$then 
\[
\jP(\XX \times \YY)
\ = \ 
\jP(\XX) \ \ov{\otimes} \  \jP(\YY)
\]
if 
\[
\begin{cases}
\ \card(\XX) \leq \aleph_1
\\[8pt]
\ \card(\YY) \leq \aleph_1
\end{cases}
\qquad 
\text{(or even $\leq \fc$ under Martin's axiom),}
\]
but
\[
\jP(\XX \times \YY)
\ \neq \ 
\jP(\XX) \ \ov{\otimes} \  \jP(\YY)
\]
if both $\XX$ and $\YY$ are uncountable and at least one of them has cardinality $> \fc$.
\\[-.25cm]

(6) \ 
Suppose that $P(\XX) = \{X_i \hsy : \hsy i \in \sI\}$ is a partition of $\XX$; 
suppose that $P(\YY) = \{Y_j \hsy : \hsy j \in \sJ\}$ is a partition of $\YY$ $-$then 
the \unz{product} of $P(\XX)$ and $P(\YY)$ is that partition  $P(\XX)  \times P(\YY)$ of $\XX \times \YY$ 
whose components are the $X_i \times Y_j$.  
Check that a product is a network and that, conversely, a network is a product. 
\\[-.25cm]

(7) \ 
Suppose that there is attached to each $i$ in an uncountable set $\sI$ 
a nonempty set $X_i$ and a nontrivial $\sigma$-algebra $\jS_i \subset \sP(X_i)$ $-$then 
$\ovs{\otimes}^* \hsx  \jS_i$ is antiatomic.
\\[-.5cm]

[Note: \ 
This need not be true of course, if $\sI$ is countable.]
\\[-.25cm]

(8) \ 
Given a class of nonempty, pairwise disjoint sets $X_i$, 
let $\jS_i$ be an atomic ring in $X_i$ $-$then $\oplus \ \jS_i$ is an atomic ring in $\bigcup \ X_i$.


\chapter{
$\boldsymbol{\S}$\textbf{5}.\quad  Problem}
\setlength\parindent{2em}
\setcounter{theoremn}{0}
\renewcommand{\thepage}{\S5-\text{P-}\arabic{page}}



\vspace{-1cm}
\[
\text{PROJECTIONS}
\]
\\[-1cm]

Let $\XX$ and $\YY$ be nonempty sets $-$then by \unz{projection onto} $\XX$ we understand the map from 
$\jP(\XX \times \YY)$ onto $\jP(\XX)$ defined by the rule
\[
\Pro_\XX (E) 
\ = \ 
\{x \in \XX : E_x \neq \emptyset\}.
\]
\\[-1cm]

Verify that
\[
\begin{cases}
\ \Pro_\XX (\bigcup \ E_i) \ = \ \bigcup \Pro_\XX (E_i)
\\[11pt]
\ \Pro_\XX (\bigcap \ E_i) \ \subset \ \bigcap \Pro_\XX (E_i)
\end{cases}
,
\]
the second containment being strict, even for a decreasing sequence, although for rectangles it is true that
\[
\Pro_\XX ((S_1 \times T_1) \cap (S_2 \times T_2) \cap \ldots) 
\ = \ 
S_1 \cap S_2 \cap \ldots
\]
if $T_1 \cap T_2 \cap \ldots \neq \emptyset$.

Let $\jS$ be a nonempty class of subsets of $\XX$; 
let $\jT$ be a nonempty class of subsets of $\YY$ $-$then, for any nonempty $E$, 
\\[-.5cm]

\qquad 
$
E \in (\jS \hsx \boxtimes \hsx \jT)_\ts 
\implies 
\Pro_\XX (E) \in \jS_\ts
$
\\[-.25cm]

\qquad 
$
E \in (\jS \hsx \boxtimes \hsx \jT)_\td 
\implies 
\Pro_\XX (E) \in \jS_\td
$
\\[-.25cm]

\qquad 
$
E \in (\jS \hsx \boxtimes \hsx \jT)_\sigma 
\implies 
\Pro_\XX (E) \in \jS_\sigma
$
\\[-.25cm]

\qquad 
$
E \in (\jS \hsx \boxtimes \hsx \jT)_\delta 
\implies 
\Pro_\XX (E) \in \jS_\delta
$
\\[-.25cm]

\qquad 
$
E \in (\jS \hsx \boxtimes \hsx \jT)_{\delta \ts} 
\implies 
\Pro_\XX (E) \in \jS_{\delta \ts} 
$
\\[-.25cm]

\qquad 
$
E \in (\jS \hsx \boxtimes \hsx \jT)_{\td \sigma} 
\implies 
\Pro_\XX (E) \in \jS_{\td \sigma}
$
\\[-.25cm]

What can be said about the operations, e.g., $\ts \delta$, $\sigma \delta$, etc.?
\\[-.25cm]

\textbf{Example} \ 
Take $\XX = \YY = [0,1]$.  
Let $\jS$ be the class comprised of all closed subintervals of $\XX$; let $\jT = \sP(\YY)$ $-$then

\[
\Pro_\XX (\jS \hsx \boxtimes \hsx \jT)_{\ts \delta} 
\ = \ 
\sP(\XX)
\]
\\[-1cm]

So, the moral is that some assumptions will have to be imposed if a positive result is to be obtained.
\\[-.5cm]

This said, prove that if $\jT$ is countably compact (Prob. VIII (\S1)), then for any nonempty $E$ ,

\[
E \in (\jS \hsx \boxtimes \hsx \jT)_{\ts \delta} 
\implies 
\Pro_\XX (E) \in \jS_{\ts \delta}.
\]
\\[-1cm]

[Recall that the countable compactness of $\jT$ implies the countable compactness of 
$\jT_{\ts \td} = \jT_{\td \ts}$ (cf. op. cit.).  
With this in mind, establish the following lemma.  
If $E_1 \supset E_2 \supset \ldots$ $(E_i \subset \XX \times \YY \forall \ i)$, and if $\forall \ x \in \XX$, the 
class $\{)E_i)_x \hsy : \hsy i = 1, 2, \ldots\}$ is countably compact, then
\[
\ \Pro_\XX (\bigcap \ E_i) \ = \ \bigcap \Pro_\XX (E_i).]
\]

Maintaining the above hypothesis on $\jT$, it can also be shown that
\[
E \in (\jS \hsx \boxtimes \hsx \jT)_{\sigma \delta} 
\implies 
\Pro_\XX (E) \in \jS_\A.
\]
Here, the sub-$A$ refers to operation $A$ (cf. \S8).
\\[-.5cm]

[Note: \ 
$\jS_\A \supset \jS_{\sigma \delta}$ but the result cannot be improved to read 
$\Pro_\XX (E) \in \jS_{\sigma \delta}$ as may be seen by example.]
\\[-.25cm]

\unz{Ref} \quad 
El Marczewski and 
C. Ryll-Nardzewski\footnote[3]{\vspace{.11 cm}\textit{Fund. Math.}, \textbf{40} (1953), 160-164.}.

\chapter{
$\boldsymbol{\S}$\textbf{6}.\quad  Extension and Generation}
\setlength\parindent{2em}
\setcounter{theoremn}{0}
\renewcommand{\thepage}{\S6-\arabic{page}}



\qquad
Let $\XX$ be a nonempty set.  
Let $\jstar$ be a property of certain nonempty classes of subsets of $\XX$ 
$-$then $\jstar$ is said to be \un{extensionally attainable} if for every subset 
$\jS$ of $\jP (\XX)$, there exists a subset $\star (\jS)$ of $\jP (\XX)$ which 
\[
\begin{cases}
\ \text{(a) \ contains} \ \jS\\
\ \text{(b) \ possesses} \ \jstar
\end{cases}
\]
and, in addition, is minimal with respect to (a) and (b).  
$\star (\jS)$, if it exists, is said to be the $\star$-class \un{generated} by $\jS$.
\\

{\small\bf Lemma 1} \ 
Property $\jstar$ is extensionally attainable iff $\jP (\XX)$ has property $\jstar$ and the intersection of any 
nonempty collection of classes having property $\jstar$ also has property $\jstar$.
\\[-.5cm]

[We omit the elementary verification.]
\\

Suppose that $\jstar$ is extensionally attainable $-$then, for any $\jS$, 
\[
\jstar (\jS)
\ = \ 
\bigcap\ \jS_i,
\]
the $\jS_i$ running over all those classes which contain $\jS$ and which possess $\jstar$.  
\\[-.25cm]

Here are some typical examples of extensionally attainable properties:
\[
\begin{cases}
\ \text{$\jS$ \ has property \ $\jstar$ \ iff \ $\jS$  \ is a lattice}\\[4pt]
\ \text{$\jS$ \ has property \ $\jstar$ \ iff \ $\jS$  \ is a ring (algebra)}\\[4pt]
\ \text{$\jS$ \ has property \ $\jstar$ \ iff \ $\jS$  \ is a $\sigma$-ring ($\sigma$-algebra)}\\[4pt]
\ \text{$\jS$ \ has property \ $\jstar$ \ iff \ $\jS$  \ is a $\delta$-ring ($\delta$-algebra)}\\[4pt]
\ \text{$\jS$ \ has property \ $\jstar$ \ iff \ $\jS$  \ is a complete ring}
\end{cases}
.
\]

On the other had, the stipulations that
\[
\begin{cases}
\ \text{$\jS$ \ has property \ $\jstar$ \ iff \ $\jS$  \ is a ring with unit}\\[4pt]
\ \text{$\jS$ \ has property \ $\jstar$ \ iff \ $\jS$  \ is a semiring}
\end{cases}
\]
are not extensionally attainable.
\\

{\small\bf Examples}  \ 
\\[-.5cm]

(1) \ 
The intersection of two rings with unit need not be a ring with unit.
\\[-.25cm]

[Take $\XX = [0,3]$.  
If \ 
$
\begin{cases}
\ \jS\\
\ \jT
\end{cases}
$
is the class of all subsets of 
$
\begin{cases}
\ [0,2]
\\
\ [1,3]
\end{cases}
$
which are either finite or have a finite complement per 
$
\begin{cases}
\ [0,2]
\\
\ [1,3]
\end{cases}
,
$
then both $\jS$ and $\jT$ are rings with unit, but their intersection $\jS \cap \jT$ consists of all finite subsets of 
$[1,2]$, hence is not a ring with unit.]
\\

(2) \ 
The intersection of two semirings need not be a semiring.
\\[-.5cm]

[Take $\XX = \{1, 2, 3\}$ $-$then
\[
\begin{cases}
\ \jS \ = \ \{\emptyset, \{1\}, \{2, 3\}, \{1, 2, 3\}\}\\
\ \jT \ = \ \{\emptyset, \{1\}, \{2\}, \{3\}, \{1, 2, 3\}\}
\end{cases}
\]
are both semirings, but their intersection 
\[
\jS \cap \jT
\ = \
\{\emptyset, \{1\},  \{1, 2, 3\}\}
\]
is not.]
\\

Suppose that $\star$ is extensionally attainable $-$then $\star$ determines a map
\[
M_\star \hsy : \hsy \jP (\jP(\XX)) \ra \jP (\jP(\XX)),
\]
namely the rule which assigns to each $\jS$ its $\star$-class $\star (\jS)$  
The fixed points for this map are exactly those classes $\jS$ having property $\star$.
The central question to be considered now is this: 
Given $\jS$, describe $\star (\jS)$.  
Naturally, the description itself will depend on $\star$.  
In terms of $M_\star$, there is a variant in that typically a generic nonempty fiber $M_\star^{-1} (\jS_0)$ is fixed in advance, 
the point being that each $\jS$ in this fiber generates the $\star$-class $\jS_0$, 
i.e., $\star (\jS) = \jS_0$, implying, therefore, that $\jS_0$ can be studied in a variety of ways.  
\\[-.5cm]

[Note: \ 
In what follows, we shall leave it up to the reader to struggle with the empty class.]
\\[-.5cm]

Let us begin with a simple illustration.  
Take $\star$ to be the 
property: ? is a lattice.  
Given a nonempty set $\jS$, we then call $\star (\jS)$ the lattice generated by $\jS$, and denote it by $\jLat (\jS)$.  
In terms of $\jS$, $\jLat (\jS)$ is the class 
$\jS_\tsd (= \jS_\tds)$ with, if necessary, the empty set adjoined.  

A slightly more complicated situation arises when we take $\star$ to be the 
property: ? is a ring.  
Given any nonempty $\jS$, we then call $\star (\jS)$ the ring generated by $\jS$ and denote it by $\jRin (\jS)$.  
Viewed abstractly,  $\jRin (\jS)$ is simply the intersection of all rings in $\XX$ containing $\jS$.  
Thus, on algebraic grounds, $\jRin (\jS)$ can be described as the class of all finite symmetric differences 
$S_1 \hsx \Delta \hsx \ldots \hsx \Delta \hsx S_m$, each $S_i$ being in turn a finite intersection of sets belonging to $\jS$.  
Consequently, if $\jS$ is finite (countable), then so is $\jRin (\jS)$.
\\[-.5cm]

[Note: \ 
Other characterizations of $\jRin (\jS)$ may be found in Exer. 3.  
Trivially, every element of $\jRin (\jS)$ is contained in some element of $\jS_\ts$ (cf. Exer. 8).]
\\[-.25cm]

\textbf{Example} 
Let $\jS$ be a semiring $-$then
\[
\jRin (\jS)
\ = \ 
\jKol (\jS).
\]
\\[-1.25cm]

Take now for $\star$  the 
property: ? is a $\sigma$-ring ($\delta$-ring).  
Given any nonempty $\jS$, we then call $\star (\jS)$ the $\sigma$-ring ($\delta$-ring)
generated by $\jS$ and denote it by $\sigma$-$\jRin (\jS)$ ($\delta$-$\jRin (\jS)$).
Observe that the notation is unambiguous in that the $\sigma$-ring ($\delta$-ring) generated by $\jS$ 
is in fact the same as the $\sigma$-ring ($\delta$-ring) generated by $\jRin (\jS)$.  
Obviously, 
\[
\text{$\delta$-$\jRin (\jS)$}
\subset
\text{$\sigma$-$\jRin (\jS)$},
\]
$\sigma$-$\jRin (\jS)$ being in fact the class of all countable unions of elements from $\delta$-$\jRin (\jS)$, 
i.e., 
\[
\text{$\sigma$-$\jRin (\jS)$}
\ = \ 
[\text{$\delta$-$\jRin (\jS)$}]_\sigma.
\]

\textbf{Examples} 

\qquad (1) \quad 
Let $\XX$ be a topological space $-$then the $\sigma$-ring generated by the open 
(or, equivalently, closed) subsets of $\XX$ is called the $\sigma$-ring of 
\un{Borel sets} in $\XX$ and is denoted by $\jBo (\XX)$.


\qquad (2) \quad 
Let $\XX$ be a Hausdorff topological space $-$then the $\delta$-ring generated by the compact subsets of $\XX$ 
is called the $\delta$-ring of \un{bounded Borel sets} in $\XX$ and is denoted by $\jBo_\tb (\XX)$.
\\[-.5cm]

[Note: \ $\XX$ is taken to be Hausdorff here in order to ensure that every compact subset of $\XX$ is a Borel set 
(all compacta then being closed, of course).  
By comparison, observe that if $\XX$ is equipped with the indiscrete topology, then the Borel sets are $\emptyset$ and $\XX$, 
but every subset of $\XX$ is compact.]
\\[-.25cm]

One cannot, in general, describe the $\sigma$-ring generated by a class of sets in purely algebraic terms.  
There are, however, useful alternative procedures, essentially transfinite in nature.
\\[-.5cm]

We have already encountered one such.  
Indeed, given $\jS$, $\sigma$-$\jRin (\jS)$ is simply the closure of 
$\jRin (\jS)$ in $\jP (\XX)_\ts$ (cf. \S4) or still (cf. \S1), 
\[
\text{$\sigma$-$\jRin (\jS)$}
\ = \ 
\bigcup\limits_{\alpha < \Omega} \ u_\alpha (\jRin (\jS)).
\]
In this connection, let us recall that $u_\alpha (\jRin (\jS))$ is the class comprised of those sets $S \subset \XX$ for which 
there exists a sequence 
$\{S_i\} \subset \bigcup\limits_{\beta < \alpha} \ u_\beta (\jRin (\jS))$ such that 
$\lim \ S_i = S$.  
The $u_\alpha (\jRin (\jS))$ are rings which increase with $\alpha$.  
Consequently, inside $\sigma$-$\jRin (\jS)$ is a transfinite sequence of rings
\[
\jRin (\jS) 
\subset \cdots \subset 
u_\alpha (\jRin (\jS))
\subset \cdots
\qquad (\alpha < \Omega),
\]
whose union is precisely $\sigma$-$\jRin (\jS)$ itself.
\\[-.25cm]

\textbf{Example} \ 
Let $\jS$ be a ring in $\XX$; 
let $\jT$ be a ring in $\YY$ 
$-$then 
\[
\text{$\sigma$-$\jRin (\jS \hsx \otimes \hsx \jT)$} 
\ = \ 
\jS \ \ovs{\otimes} \ \jT.
\]
More generally, let $\jS_i$ be an algebra in $X_i$ $(i \in \sI$, $\sI$ infinite) $-$then
\[
\text{$\sigma$-$\jRin (\otimes^* \ \jS_i)$} 
\ = \ 
\ovs{\otimes}^{\hsy \ast} \ \jS_i.
\]

Starting from $\jS$, we shall now define by transfinite recursion a class $\jS_\alpha$ for each ordinal $\alpha < \Omega$.  
Thus putting $\jS_0 = \jS$, write
\[
\jS_\alpha
\ = \ 
\left(
\bigcup_{\beta < \alpha}^\infty \ \jS_{\beta} 
\right)_{r \hsy \sigma}
\qquad (\alpha < \Omega).
\]
Observe that the $\jS_\alpha$ increase with $\alpha$.
\\

\textbf{Lemma 2} \ 
\un{We have}

\[
\text{$\sigma$-$\jRin (\jS)$} 
\ = \ 
\bigcup\limits_{\alpha < \Omega} \ \jS_{\alpha}.
\]
\\[-1.25cm]

To see what the rationale behind the construction is, replace $\sigma$ by $s$ $-$then 
$\jS_0 = \jS$, 
$\jS_1 = \jS_\rs$, 
$\jS_2 = \jS_{\rs \hsy\rs}$, 
the ring generated by $\jS$ (cf. Exer. 3). 
\\[-.5cm]

[Note: \ 
Trivially, every element of $\sigma$-$\jRin (\jS)$ is contained in some element of $\jS_\sigma$ (cf. Exer. 8).]
\\[-.25cm]

\un{Proof of Lemma 2} \ 
There are two steps in the argument.
\\[-.25cm]

\qquad (1) \quad 
$\bigcup\limits_{\alpha < \Omega} \ \jS_\alpha$ is contained in $\sigma$-$\jRin (\jS)$.
\\[-.25cm]

\qquad (2) \quad 
$\bigcup\limits_{\alpha < \Omega} \ \jS_\alpha$ is a $\sigma$-ring.
\\[-.25cm]

\noindent
\un{Re (1)} \quad 
By definition $\jS_0 = \jS \subset \sigma$-$\jRin (\jS)$; 
in addition, $\emptyset \in \jS_1$.  
Proceeding by transfinite induction, assume that $\jS_\beta \subset \sigma$-$\jRin (\jS)$ 
for every $\beta < \alpha$ and consider a typical element $S \in \jS_\alpha$ $-$then 
$S$ is a countable union, say $\bigcup \ S_i$, where each $S_i$ has the form 
$A_i$ or $A_i - B_i$, with 
\[
A_i \hsx B_i \in 
\bigcup\limits_{\beta < \alpha} \ \jS_\beta
\subset 
\text{$\sigma$-$\jRin (\jS)$}.
\]
Thus $S_i \in \sigma$-$\jRin (\jS)$ and so 
$S = \bigcup \ S_i \in \sigma$-$\jRin (\jS)$, which implies that 
$S_\alpha \subset \sigma$-$\jRin (\jS)$.  
This completes the proof of (1).
\\[-.5cm]

\noindent
\un{Re (2)} \quad 
Let $\{S_i\}$ be a sequence in 
$\bigcup\limits_{\alpha < \Omega} \ \jS_\alpha$ 
$-$then we claim that 
$\bigcup \ S_i \in \bigcup\limits_{\alpha < \Omega} \ \jS_\alpha$.  
To prove it, note that for each $i$ there is an $\alpha_i$ such that $S_i \in \jS_{\alpha_i}$.  
Select, as is possible, an $\alpha < \Omega$ such that $\alpha_i < \alpha$ $\forall \ i$ $-$then
\[
\bigcup \ S_i 
\in 
\left(
\bigcup_{i = 1}^\infty \ \jS_{\alpha_i} 
\right)_{r \hsy \sigma}
\subset 
\jS_\alpha
\subset
\bigcup\limits_{\alpha < \Omega} \ \jS_\alpha,
\]
as claimed.  
In an entirely analogous manner, one can show that if 
$S$, $T \in \bigcup\limits_{\alpha < \Omega} \ \jS_\alpha$, then 
$S - T \in \bigcup\limits_{\alpha < \Omega} \ \jS_\alpha$.  
This completes the proof of (2). 
\\[-.25cm]

Hence the lemma. //
\\[-.25cm]

The transfinite description of $\sigma$-$\jRin (\jS)$ provided by Lemma 2 carries with it an added bonus in that an estimate for the cardinality of 
$\sigma$-$\jRin (\jS)$ can easily be obtained.  
To this end, we can suppose that $\card(\jS)\geq 2$ since 
\[
\text{$\sigma$-$\jRin (\{\jS\})$} \ = \ 
\begin{cases}
\ \{S, \emptyset\} \hspace{0.5cm} S \neq \emptyset\\[4pt]
\ \ \{\emptyset\} \hspace{0.75cm} S = \emptyset
\end{cases}
.
\]
Our estimate then reads:
\[
\card(\text{$\sigma$-$\jRin (\jS)$}) 
\ \leq \ 
\card(\jS)^{\aleph_0}.
\]
Indeed, the assumption that $\card(\jS) \geq 2$,  in conjunction with consideration of the ways in which the sets 
$\bigcup \ S_i \in \jS_1$ can be formed (at most $\card(\jS)^2$ choices for each $S_i$), 
leads at once to the conclusion that 
$\card(\jS_1) \leq (\card(\jS)^2)^{\aleph_0} = \card(\jS)^{\aleph_0}$.  
Utilizing now transfinite induction, suppose that 
$\card(\jS_\beta) \leq \card(\jS)^{\aleph_0}$ 
for all $\beta$ such that $1 \leq \beta < \alpha$, where $1 < \alpha < \Omega$ $-$then
\[
\card
\left(
\bigcup\limits_{\beta < \alpha} \ \jS_\beta
\right)
\ \leq \ 
\aleph_0 \cdot \card\left(\jS\right)^{\aleph_0}
\ = \ 
\card\left(\jS\right)^{\aleph_0}
\]
and so, arguing as above, it follows that 
$\card(\jS_\alpha) \leq \card\left(\jS\right)^{\aleph_0}$.  
Consequently, for every $\alpha$ with $0 \leq \alpha < \Omega$, 
$\card(\jS_\alpha) \leq \card(\jS)^{\aleph_0}$.  
All told therefore, 

\begin{align*}
\text{$\card(\sigma$-$\jRin (\jS))$} \ 
&=\ 
\card
\left(
\bigcup\limits_{\alpha < \Omega} \ \jS_\alpha
\right)
\\[15pt]
&\leq \ 
\aleph_1 \cdot \card\left(\jS\right)^{\aleph_0}
\\[15pt]
&\leq \ 
2^{\aleph_0} \cdot \card\left(\jS\right)^{\aleph_0}
\\[15pt]
&= \ 
\card\left(\jS\right)^{\aleph_0}.
\end{align*}

[Note: \ 
If $\jS$ is finite, then, of course, $\sigma$-$\jRin (\jS)$ is finite, there being the estimate
\[
\#(\text{$\sigma$-$\jRin (\jS)$}) 
\ \leq \ 
2^{2^{\# (\jS)}}
\]
which is even attainable under the obvious conditions.]
\\[-.25cm]

\textbf{Example} \ 
Let $\XX$ be a topological space with weight $\aleph_0$ $-$then the cardinality of the class of Borel sets in $\XX$ cannot 
exceed the cardinality of the continuum.  
In fact, the cardinality in question is the same as that of the $\sigma$-ring generated by the open sets and this cannot exceed 
$\fc^{\aleph_0} = \fc$.  
Specialize and suppose in addition that $\XX$ is a metric space which is complete and perfect, so that $\card (\XX) = \fc$.  
Because there are then $\fc$ open
sets, the cardinality of the class of Borel sets in $\XX$ is exactly $\fc$, thus is $< 2^\fc$, the cardinality of $\jP (\XX)$.
\\[-.25cm]

Let $\star$ be the property: 
$? = ?_\sigma$ 
and 
$? = ?_\delta$.  
It is clear that 
$\star$ is extensionally attainable.  
Given any nonempty $\jS$, we then write
$\jS_\tB $
for
$\star (\jS)$
and refer to 
$\tM_\star$ 
as 
\un{operation B}.  
Obviously, 
$\jS_\tB = \jS_{\tB\tB}$ 
and
\[
\begin{cases}
\ 
\jS_\sigma
\subset 
\jS_\tB
\\[4pt]
\ 
\jS_\delta
\subset 
\jS_\tB
\end{cases}
\qquad
\begin{cases}
\ 
\jS_\tB 
\ = \ 
\jS_{\tB\hsy\sigma}
\ = \ 
\jS_{\sigma\hsy\tB}
\\[4pt]
\ 
\jS_\tB 
\ = \ 
\jS_{\tB\hsy\delta}
\ = \ 
\jS_{\delta\hsy\tB}
\end{cases}
.
\]

The topological interpretation of 
$\jS_\tB$ 
is very simple.  
Indeed, 
$\jS_\tB$ 
is nothing more nor less than the closure in 
$\jP(\XX)_\ts$ 
of 
$\jS_{\ts\hsy\td} = \jS_{\td\hsy\ts}$, 
thus in particular, is the closure of 
$\jLat (\jS)$ 
in 
$\jP(\XX)_\ts$ 
if 
$\emptyset \in \jS$.
\\[-.5cm]

[Note: \ 
The reader will agree that the closure of 
$\jS$ 
itself in 
$\jP(\XX)_\ts$ 
will, in general, be a proper subset of 
$\jS_\tB$.]
\\[-.25cm]

There is an equally straightforward transfinite description of 
$\jS_\tB$ .  
Namely, put 
$\jB^{(0)} (\jS) = \jS$, 
$\jB_{(0)} (\jS) = \jS$ 
and define via transfinite recursion the classes
$\jB^{(\alpha)}  (\jS)$, 
$\jB_{(\alpha)}  (\jS)$
by writing
\[
\begin{cases}
\ \ds
\jB^{(\alpha)}  (\jS)
\ = \ 
\Big[
\bigcup\limits_{\beta < \alpha} \ 
\jB_{(\beta)}  (\jS)
\Big]_\sigma
\\[18pt]
\ \ds
\jB_{(\alpha)}  (\jS)
\ = \ 
\Big[
\bigcup\limits_{\beta < \alpha} \ 
\jB^{(\beta)}  (\jS)
\Big]_\delta
\end{cases}
\qquad (\alpha < \Omega)
.
\]
The 
$\jB^{(\alpha)}  (\jS)$, 
$\jB_{(\alpha)}  (\jS)$
evidently increase with $\alpha$ and for 
$
\begin{cases}
\ \alpha \geq 1
\\
\ 
\alpha > 1
\end{cases}
$
\[
\begin{cases}
\ 
\Big[
\jB^{(\alpha)}  (\jS)
\Big]_\sigma
\ = \ 
\jB^{(\alpha)}  (\jS)
\\[18pt]
\Big[
\jB^{(\alpha)}  (\jS)
\Big]_\td
\ = \ 
\jB^{(\alpha)}  (\jS)
\end{cases}
,
\]

\[
\begin{cases}
\ 
\Big[
\jB_{(\alpha)}  (\jS)
\Big]_\delta
\ = \ 
\jB_{(\alpha)}  (\jS)
\\[18pt]
\Big[
\jB_{(\alpha)}  (\jS)
\Big]_\ts
\ = \ 
\jB_{(\alpha)}  (\jS)
\end{cases}
.
\]
In addition, if 
\[
\jB  (\jS: \alpha) 
\ = \ 
\jB^{(\alpha)}  (\jS) 
\hsx \cap \hsx 
\jB_{(\alpha)}  (\jS), 
\]
then
\[
\jB^{(\alpha)}  (\jS) 
\hsx \cup \hsx 
\jB_{(\alpha)}  (\jS)
\hsx \subset \hsx 
\jB  (\jS: \alpha + 1).
\]
\\[-.75cm]

Our hierarchy may be visualized as follows: 
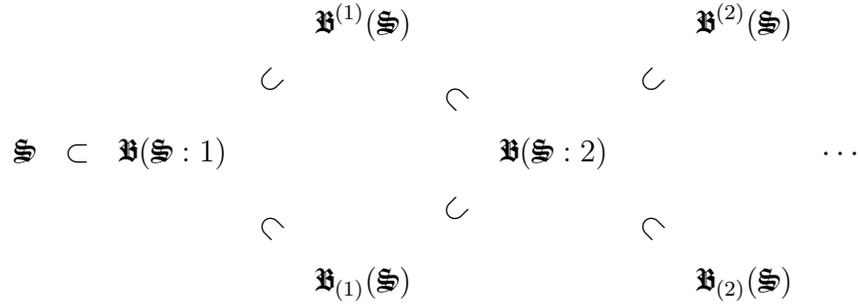
\[
\begin{tikzcd}[sep=0.05em]
&&&&
\jB^{(1)}  (\jS) 
&&&&
\jB^{(2)}  (\jS) 
\\
&&&
\rotatebox{45}{$\subset$}
&&
\rotatebox{-45}{$\subset$}
&&
\rotatebox{45}{$\subset$}
\\
\jS 
&
\subset 
&
\jB(\jS: 1)
&&&&
\jB(\jS: 2)
&&&
\cdots
\\
&&&
\rotatebox{-45}{$\subset$}
&&
\rotatebox{45}{$\subset$}
&&
\rotatebox{-45}{$\subset$}
\\
&&&&
\jB_{(1)}  (\jS) 
&&&&
\jB_{(2)}  (\jS) 
\end{tikzcd}
\]
\\[-.25cm]

\textbf{Lemma 3} \ 
\un{We have}
\[
\jS_\tB 
\ = \ 
\begin{cases}
\ \ds
\bigcup\limits_{\alpha < \Omega} \ 
\jB^{(\alpha)}  (\jS)
\\[18pt]
\ \ds
\bigcup\limits_{\alpha < \Omega} \ 
\jB_{(\alpha)}  (\jS)
\end{cases}
.
\]

[One need only imitate the argument used in the proof of Lemma 2.]
\\

There is a variant on the preceding definitions which is frequently encountered in the literature.  
To describe it, let us recall that any ordinal $\alpha$ can be written uniquely in the form 
$\alpha = \lambda + n$, where $\lambda$ is a limit ordinal or zero and $n$ is a nonnegative integer  
($\alpha$ then being termed \un{odd} or \un{even} according to the parity of $n$).  
This being so, put 
$\jB^{[0]} (\jS) = \jS$, 
$\jB_{[0]} (\jS) = \jS$, 
and define via transfinite recursion the classes 
$\jB^{[\alpha]} (\jS)$, 
$\jB_{[\alpha]} (\jS)$
by writing

\[
\jB^{[\alpha]} (\jS) 
\ = \ 
\begin{cases}
\ 
\jB^{(\alpha)}  (\jS) 
\quad \text{if $\alpha$ is odd}
\\[11pt]
\
\jB_{(\alpha)}  (\jS)
\quad \text{if $\alpha$ is even} 
\end{cases}
\qquad (\alpha < \Omega), 
\]

\[
\jB_{[\alpha]} (\jS)
\ = \ 
\begin{cases}
\ 
\jB^{(\alpha)}  (\jS) 
\quad \text{if $\alpha$ is even}
\\[11pt]
\
\jB_{(\alpha)}  (\jS)
\quad \text{if $\alpha$ is odd} 
\end{cases}
\qquad (\alpha < \Omega).
\]
Then it is again the case that
\[
\jS_\tB
\ = \ 
\begin{cases}
\ \ds
\ \bigcup\limits_{\alpha < \Omega} \ \jB^{[\alpha]} (\jS)\\[18pt]
\ \ds
\ \bigcup\limits_{\alpha < \Omega} \ \jB_{[\alpha]} (\jS)\\[18pt]
\end{cases}
.
\]
Note too that if for some $\alpha \geq 1$, 
$\jB^{[\alpha]} (\jS) = \jB^{[\alpha+1]} (\jS)$ (or 
$\jB_{[\alpha]} (\jS) = \jB_{[\alpha+1]} (\jS)$), 
then 
$\jB^{[\alpha]} (\jS) = \jS_\tB$ 
(or 
$\jB_{[\alpha]} (\jS) = \jS_\tB$).
For of the two classes 
$\jB^{[\alpha]} (\jS)$ and 
$\jB^{[\alpha+1]} (\jS)$  
(or 
$\jB_{[\alpha]} (\jS)$ and
$\jB_{[\alpha+1]} (\jS)$), 
one is closed under countable unions while the other is closed under countable intersections, hence, when they coincide, 
$\jB^{[\alpha]} (\jS)$ 
(or 
$\jB_{[\alpha]} (\jS)$) must give $\jS_\tB$.
\\[-.25cm]

\textbf{Example} \ 
By the \un{Kolmogoroff number} \ $\tK(\jS)$ of $\jS$, we understand the smallest ordinal $\alpha$ such that 
$\jB^{[\alpha]} (\jS) = \jS_\tB$.  
The apparent asymmetry in the definition is, of course, essentially illusory.  
There are initial and terminal possibilities, namely, if $\jS = \jS_\tB$ to begin  with, then $\tK(\jS)= 0$, 
whereas, if 
$\jB^{[\alpha]} (\jS) \neq \jS_\tB$ $\forall \ \alpha < \Omega$, then we agree to take 
$\tK(\jS) = \Omega$.  
Two problems can then be posed.
\\[-.25cm]

\qquad (1) \quad 
Given $\jS$, determine $\tK(\jS)$.
\\[-.25cm]

\qquad (2) \quad 
Given $\alpha$, find an $\jS$ such that $\tK(\jS) = \alpha$.
\\[-.25cm]

\noindent
Here, we shall deal with the second, setting aside the systematic consideration of the first for now.  
Let us mention in passing, however, that examples for which
$\tK(\jS) = \Omega$ do in fact abound, the simplest instance being the case when $\jS$ is the class of all open 
(or closed) subintervals of the line.  
In Exer. 14 (\S1), it was pointed out that there exist easy examples of classes $\jS$ such that
$\tK(\jS) = 0$, 1 and 2, but to get an example when $\tK(\jS) = 3$ turned out to be surprisingly difficult, 
at least if one works on the line, the point being that the classical solution utilizes the continuum hypothesis 
(but see the paper of Maly\v sev referenced below).  
Actually, operating within ZFC alone, it is possible to give a complete answer to (2) in that $\forall \ \alpha < \Omega$, 
there exists a nonempty set $\XX$ and a nonempty class $\jS$ contained in $\jP(\XX)$ 
such that $\tK(\jS) = \alpha$.  
While interesting, we shall forgo the details, settling instead for an indication.  
To begin with, it is best to generalize the problem, 
replacing $\jP(\XX)$ by a complete Boolean algebra $\jA$ and then introducing a 
notion of Kolmogoroff number $\tK(\jA)$ for $\jA$.  
This done, the crucial step in the argument consists of proving that $\forall \ \alpha < \Omega$, there exists a complete Boolean algebra $\jA$ 
satisfying the countable chain condition with $\tK(\jA) = \alpha$.  
Thanks to the Loomis-Sikorski theorem, any $\sigma$-complete Boolean algebra is isomorphic to a $\sigma$-algebra of subsets of some set $\XX$ 
modulo a $\sigma$-ideal.  
Accordingly, $\jA$ can be represented as a certain quotient per a certain $\XX$ and finally, using the fact that $\tK(\jA) = \alpha$, 
one produces without difficulty a subset $\jS$ of $\jP(\XX)$ with the property that $\tK(\jS) = \alpha$.
\\[-.25cm]

[This result is due to Kunen; cf. 
A. Miller\footnote[2]{\vspace{.11 cm}\textit{Ann. Amer. Logic}, \textbf{16} (1979), 233-267.}.]
\\[-.25cm]

For a fairly simple example of a class $\jS$ such that $\tK(\jS) = 3$ 
(and not involving the continuum hypothesis), see 
V. Maly\v sev\footnote[3]{{\fontencoding{OT2}\selectfont
[V. Malyxev]}
\vspace{.11 cm}
\textit{Vestnik Moskov. Univ. Ser. I Mat. Meh.}, \textbf{6} (1965), 8-10.}.
\\[-.25cm]


On the basis of the definitions, 
\[
\jS_\tB 
\subset
\text{$\sigma$-$\jRin (\jS)$},
\]
the containment being strict in general.  
Indeed, 
$\sigma$-$\jRin (\jS) = \jS_{r \hsy \tB}$ 
but it need not be true that 
$\jS_{r \hsy \tB} = \jS_{\tB \hsy r}$, say.  
For example, take 
$\XX = \{1, 2, 3\}$ and let $\jS = \{\{1\}, \{1, 2\}, \{1, 2, 3\}\}$ $-$then
$\{1, 3\} \in \jS_{r \hsy \tB} -  \jS_{\tB \hsy r}$.
\\[-.5cm]

[Note: \ It can even happen that all the classes 
$\jS_\tB$, $\jS_{\tB \hsy r}$, $\jS_{\tB \hsy r \hsy r}, \ldots$ are distinct.]
\\[-.25cm]

\textbf{Lemma 4} \ 
$\jS_\tB = \sigma$-$\jRin (\jS)$ \ \un{iff} \ $\jS_r \subset \jS_\tB$.
\\[-.5cm]

[The necessity is clear.  
As for the sufficiency, observe that
\begin{align*}
\jS_r  \subset \jS_\tB 
&\implies
\jS_{r \hsy \tB} \subset \jS_{\tB \hsy \tB} = \jS_\tB
\\[11pt]
&\implies
\text{$\sigma$-$\jRin (\jS)$} \subset \jS_\tB.]
\end{align*}
\\[-1cm]

\textbf{Example} \ 
Let $\XX$ be  a topological space $-$then, traditionally, one writes
\[
\begin{cases}
\ \fF = \text{closed subsets of $\XX$}\\[4pt]
\ \fG = \text{open subsets of $\XX$}
\end{cases}
,
\]
the classical resolutions 

\[
\begin{cases}
\ \ds
\fF_\tB 
\ = \  \bigcup\limits_{\alpha < \Omega} \ \jB^{[\alpha]} (\fF) 
\ \equiv \ 
\fF_\alpha
\\[18pt]
\ \ds
\fG_\tB 
\ = \  
\bigcup\limits_{\alpha < \Omega} \ \jB_{[\alpha]} (\fG) 
\ \equiv \ 
\fG_\alpha
\end{cases}
\]
then being
\[
\begin{cases}
\ 
\fF \subset \fF_\sigma \subset \fF_{\sigma \hsy \delta} \cdots
\\[8pt]
\ 
\fG \subset \fG_\delta \subset \fG_{\delta \hsy \sigma} \cdots
\end{cases}
.
\]
The associated $\sigma$-rings (actually $\sigma$-algebras)
\[
\begin{cases}
\ \text{$\sigma$-$\jRin (\fF)$}\\[4pt]
\ \text{$\sigma$-$\jRin (\fG)$}
\end{cases}
\]
are equal, yielding, by definition, the Borel subsets of $\XX$.  
We then ask: Is
\[
\begin{cases}
\ \fF_\tB 
\ = \ 
\text{$\sigma$-$\jRin (\fF)$}?
\\[8pt]
\ \fG_\tB 
\ = \ 
\text{$\sigma$-$\jRin (\fG)$}?
\end{cases}
\]
Thanks to Lemma 4, these questions are equivalent, i.e., 
\[
\fF_\tB 
\ = \ 
\text{$\sigma$-$\jRin (\fF)$}
\iff
\fG_\tB 
\ = \ 
\text{$\sigma$-$\jRin (\fG)$}.
\]
To be specific, we shall work with $\fG$ $-$then, in decreasing order of strengh, the relation
\[
\fG_\tB 
\ = \ 
\text{$\sigma$-$\jRin (\fG)$}
\]
is forced by the following conditions.
\\[-.25cm]

\qquad $(\tC_1)$ \quad 
Every closed subset of $\XX$ is in $\fG_\delta$.
\\[-.5cm]

\qquad $(\tC_2)$ \quad 
Every closed subset of $\XX$ is in $\fG_\alpha$ for some fixed $\alpha$. 
\\[-.5cm]

\qquad $(\tC_3)$ \quad 
Every closed subset of $\XX$ is in $\fG_\alpha$ for some $\alpha$, but no
fixed $\alpha$ suffices.
\\[-.25cm]

Are there topological spaces $\XX$ satisfying these conditions?  
It is easy to meet 
$(\tC_1)$: 
Simply require that $\XX$ be perfectly normal (in particular, a metric space); Exer. 28 is also relevant.  
Turning to 
$(\tC_2)$, 
we claim that $\forall \ \alpha$ $(1  <  \alpha < \Omega)$ there exists a topological space $\XX_\alpha$ such that 
$\alpha$ is the smallest ordinal for which every closed subset $\XX_\alpha$ is in $\fS_\alpha$.  
Here is the construction.  
Fix $\alpha$, $1 \ < \alpha < \Omega$ $-$then, as a consequence of certain generalities established in 
\S3  
there exists a subset 
$S_\alpha$ of $\XX = \R$ (usual topology) which is in $\fG_\alpha$ 
but is not in $\fG_\beta$
for any $\beta < \alpha$.  
This being so, take for our space $\XX_\alpha$ the real line topologized by specifying that the open sets are to be all sets of the form 
$U \cup V$, where $U$ is open in the usual topology and $V$ is any subset of $\XX_\alpha - S_\alpha$.  
It is not difficult to see that $\XX_\alpha$ is normal and Hausdorff, and has the required properties.  
As for 
$(\tC_3)$, 
it is in fact possible to construct an example having the requisite property, at least if the continuum hypothesis is admitted 
(cf. Problem V. (\S6)).  
There is one final point to be considered:  
Do there exist examples of topological spaces $\XX$ such that 
\[
\fG_\tB 
\ \neq \ 
\text{$\sigma$-$\jRin (\fG)$}?
\]
The answer is an emphatic ``yes''!  
Consider 
\[
\XX
\ = \ 
[0, \Omega]
\quad \text{in the order topology}.
\]
\\[-1.5cm]

\noindent
or 
\[
\XX
\ = \ 
[0, 1]^{[0, 1]}
\quad \text{in the product topology}.
\]
In the first case, $\{\Omega\}$, while closed, is not in $\fS_\tB$; 
in the second case, $\{c\}$ ($c$ a constant), while closed, is not in $\fS_\tB$.  
Note that in both cases, $\XX$ is a compact Hausdorff space.
\\[-.5cm]


Let $\jS$ be a nonempty subset of $\jP(\XX)$ $-$then by 
$\jS_{\ts_\td}$ $(\jS_{\sigma_\td})$
we understand the class of subsets of $\XX$ comprised of all nonempty finite (finite or countable) disjoint unions of sets in $\jS$.
\\[-.25cm]

\textbf{Example} \ 
Let $\jS$ be a semiring $-$then
\[
\jRin (\jS) 
\ = \ 
\jS_{\ts_\td},
\]
but $\jS_{\ts_\td}$ need not be a ring (cf. Exer. 6 (\S4)).
\\[-.25cm]

Given a nonempty $\jS$, the notations of generated algebra and $\sigma$-algebra are clear, 
as are the notions $\jAlg (\jS)$ and $\sigma$-$\jAlg (\jS)$.
We have
\[
\begin{cases}
\ \jRin (\jS) \subset \jAlg (\jS) 
\\[8pt]
\ \sigma\dash\jRin (\jS) \subset \sigma\dash\jAlg (\jS)
\end{cases}
,
\]
with 
\[
\begin{cases}
\ \jAlg (\jS)
 \ = \  
 \{S : S \in \jRin (\jS)   \quad \text{or} \quad \XX - S \in    \jRin (\jS)\}     
\\[8pt]
\ \sigma\dash\jAlg (\jS) 
\ = \ 
\{S : S \in \sigma\dash\jRin (\jS)   \quad \text{or} \quad \XX - S \in    \sigma\dash\jRin (\jS)\} 
\end{cases}
,
\]
\\[-.5cm]

\noindent
that is, 
\[
\begin{cases}
\ \jAlg (\jS) 
\ = \  
\jRin (\jS) , \{\XX\})
\\[8pt]
\ \sigma\dash\jAlg (\jS) 
\ = \  
\sigma\dash\jRin (\jS , \{\XX\})
\end{cases}
.
\]
\\[-.5cm]

On algebraic grounds alone, it is plain that 
\[
\jAlg (\jS)
\ = \ 
[\jS \hsx \cup \hsx \jS_\tc]_{\td\ts}.
\]
Actually, slightly more is true, viz.
\[
\jAlg (\jS)
\ = \
[\jS \hsx \cup \hsx \jS_\tc]_{\td\ts_\td}.
\]

Topologically, $\sigma$-$\jAlg (\jS)$ can be viewed as the closure of $\jAlg (\jS)$ in $\jP(\XX)_\ts$.  
On the other hand, thanks to Lemma 4, 
\[
[\jS \hsx \cup \hsx \jS_\tc]_\tB
\ = \ 
\sigma\dash\jAlg (\jS),
\]
leading, thereby, to the attendent transfinite descriptions.
\\[-.5cm]

Let $\star$ be the property: \ 
? = $?_{\sigma_\td}$ and 
? = $?_\delta$.
It is clear that $\star$ is extensionally attainable.  
Given any nonempty $\jS$, we then write $\jS_{\tB_d}$ for $\star(\jS)$ and refer ro $\tM_\star$ as 
\un{operation $\tB_\td$}.  
Obviously,
$\jS_{\tB_d} = \jS_{\tB_d \hsy \tB_d}$ and 
\[
\begin{cases}
\ \jS_{\sigma_\td} \subsetx \jS_{\tB_d}    \\[4pt]
\ \jS_\delta \subsetx \jS_{\tB_d}  
\end{cases}
\hspace{2.6cm} 
\begin{cases}
\ \jS_{\tB_d}  = \jS_{\tB_d \hsy \sigma_\td}    = \jS_{\sigma_\td \hsy \tB_d} 
\\[8pt]
\ \jS_{\tB_d} = \jS_{\tB_d \hsy \delta}   =  \jS_{\delta \hsy\tB_d}
\end{cases}
,
\]
with
\[
\jS_{\tB_d} \subsetx \jS_\tB,
\]
the containment being strict in general, as can be seen by taking 
$\XX = \{1, 2, 3\}$ and letting 
$\jS = \{\{1\}, \{1, 2\}, \{1, 3\}\}$.  
We shall leave
it up to the reader to discuss the transfinite aspects of operation $\tB_\td$.
\\[-.25cm]

\textbf{Lemma 5} \quad 
$\jS_{\tB_d} = \sigma\dash\jAlg (\jS)$  \ \un{iff} \ $\jS_\tc \subsetx \jS_{\tB_d}$.
\\[-.5cm]

[The necessity is clear.  
As for the sufficiency, observe that 
\\[-.25cm]

\qquad $\jS_\tc \subsetx  \hsx\jS_{\tB_\td}$  
\begin{align*}
&\implies \quad
[\jS_{\tB_\td} \hsx \cap \hsx \jS_{\tB_\td \tc}]_\tB
\ = \ 
\jS_{\tB_\td} \hsx \cap \hsx \jS_{\tB_\td \tc}
\\[11pt]
%
&\implies \quad
[\jS \hsx \cup \hsx \jS_\tc]_\tB 
\subsetx 
\jS_{\tB_\td}
\\[11pt]
%
&\implies \quad
\sigma \dash \jAlg (\jS) \subsetx \jS_{\tB_\td}.]
\end{align*}
\\[-1cm]

\textbf{Example} \ 
Let $\XX$ be a topological space $-$then by a 
\unz{zero set} 
in $\XX$ we mean any set of the form $f^{-1} (0)$, where $f : \XX \ra \R$ is continuous.  
In this connection, observe that it is not restrictive to suppose that $f(\XX) \subset [0,1]$.  
The complements in $\XX$ of the zero sets are called the 
\unz{cozero sets}.  
Agreeing to write
\[
\gZ 
\ = \ 
\text{zero sets in $\XX$},
\]
we have $\jZ \subset \jF$, the containment being strict in general (cf. Exer. 30), but there being coincidence if, 
e.g., $\XX$ is perfectly normal.  
Note that

\[
\begin{matrix*}[l]
(1)  \quad \emptyset, \ \XX \in \gZ, 
&&&&(2) \quad \gZ = \gZ_\tS,
\\[11pt]
(3)  \quad \gZ = \gZ_\delta,
&&&&(4) \quad \gZ_\tc \subset \gZ_\sigma.
\end{matrix*}
\]
In addition, given disjoint $Z_1$, $Z_2 \in \gZ$, there exist disjoint $U_1$, $U_2 \in \gZ_\tc$ such that 
\[
Z_1  \subset U_1, 
\quad 
Z_2 \subset U_2.
\]
This said, the 
\unz{Baire sets} in $\XX$ are by definition the elements of the $\sigma$-algebra $\jBa (\XX)$ generated by $\gZ$.  
Every Baire set is a Borel set but, in general, not vice-versa (cf. Exer. 32).  
Owing to Lemma 4 and property (4) supra,
\[
\gZ_\tB 
\ = \ 
\jBa (\XX).
\]
Because
\[
\gZ_\tc \subset \gZ_\sigma 
\quad \text{iff} \quad 
\gZ \subset \gZ_{\tc \hsy \delta},
\]
it follows from Lemma 5 that 
\[
(\gZ_\tc)_{\tB_\td} \ = \ \jBa (\XX).
\]
It is also true that
\[
\gZ_{\tB_\td} \ = \ \jBa (\XX),
\]
although this is not immediate.  
On the basis of Lemma 5 again, our assertion is equivalent to the statement that 
$\gZ_\tc \subset \gZ_{\tB_\td}$.
\\[-.25cm]

\un{Claim} \ 
Take $\XX = \R$ $-$then
\[
\gZ_\tc 
\subsetx 
\gZ_{\sigma_\td \hsy \delta \hsy \sigma_\td}.
\]
\\[-1.25cm]

[To appreciate the subtlety of this point, the reader may find it instructive to prove directly that 
$]a,b[$ $(a < b)$ does not belong to $\gZ_{\sigma_\td \hsy \delta }$!]
\\[-.25cm]

Admit the claim $-$then, for any topological space $\XX$, 
\[
\gZ_\tc 
\subsetx 
\gZ_{\sigma_\td \hsy \delta \hsy \sigma_\td}
\]
and, consequently, 
$\gZ_\tc \subset \gZ_{\tB_\td}$, as desired.
Indeed, if $U \in \gZ_\tc$, then there exists a continuous function $f : \XX \ra [0,1[$ such that 
\[
U 
\ = \
\{x \in \XX \hsy : \hsy f(x) \in \hsx ]0,1[\}.
\]
Now, in vew of the claim, $]0,1[ \ \in  \gZ_{\sigma_\td \hsy \delta \hsy \sigma_\td}$ (per $\R$), and so 
\[
U 
\ = \ 
f^{-1} (\hsy ]0,1[ \hsy) \in  \gZ_{\sigma_\td \hsy \delta \hsy \sigma_\td}.
\]
\\[-1cm]

\noindent
\un{Proof of Claim} \ 
Let $U \in \gZ_\tc$ $-$then $U$ is open, hence is a finite or countable union of open, pairwise disjoint intervals.  
Accordingly, there is no loss of generality in supposing that 
$U = ]a,b[$ $(a < b)$.  
Let $\{I_m\}$ be a sequence of closed, pairwise disjoint intervals in $\R$ whose union is dense in $]a,b[$.  
\\[-.25cm]

\noindent
Put
\[
S 
\ = \ 
]a,b[ \ - \  \cup I_m.
\]
Then the closure $\bar{S}$ of $S$ in $\R$ is a closed, nowhere dense subset of $[a.b]$, and 
$\bar{S} - S$ is a countable set $\{x_n\}$ consisting of $a$, $b$ and the endpoints of $I_m$.  
Since $]a,b[$ is
\[
S \cup I_1 \cup I_2 \ldots \hsx, 
\]
the union being countable and disjoint, and
\[
S 
\ = \ 
\bigcap \ (\bar{S} - \{x_1, \ldots, x_n\}),
\]
it will be enough to prove that $\bar{S} - \{x_1, \ldots, x_n\} \in \gZ_{\sigma_\td}$.  
However, because $\bar{S}$ is nowhere dense, thus 0-dimensional, one can certainly write 
$\bar{S} - \{x_1, \ldots, x_n\}$ as a countable disjoint union of sets which are closed in 
$\bar{S}$, and so too in $\R$.
\\[-.5cm]

[Note: \ 
Suppose that $\XX$ is a perfectly normal topological space $-$then, of course,
\[
\jF 
\ = \ 
\jG_\tc \subset \jG_\delta 
\ \implies \ 
\jG_{\tB_\td}
\ = \ 
\sigma \dash \jAlg (\jG).
\]
Furthermore, in this case,
\[
\jG 
\ = \ 
\jF_\tc 
\ = \ 
\gZ_\tc
\ \subset \ 
\gZ_{\sigma_\td \hsy \delta \hsy \sigma_\td}
\ = \ 
\jF_{\sigma_\td \hsy \delta \hsy \sigma_\td}
\]

\[
\implies \quad
\jF_{\tB_\td}
\ = \ 
\sigma \dash \jAlg (\jF).
\]
Here, therefore, 
\[
\jBa(\XX) \ = \ 
\begin{cases}
\ \jF_{\tB_\td} \\[8pt]
\ \jG_{\tB_\td} 
\end{cases}
.
\]
We remark, in passing, that perfect normality, while sufficient, is not necessary in order to draw these conclusions (cf. Exer. 33).]
\\[-.25cm]

Suppose that $\jS$ is a $\sigma$-lattice containing $\XX$.  
Put
\[
\Sigma_0 \hsy (\jS) \ = \ \jS, 
\quad 
\Pi_0\hsy  (\jS) \ = \ \jS_\tc
\]
and define via transfinite recursion the classes $\Sigma_\alpha \hsy (\jS)$, $\Pi_\alpha \hsy (\jS)$ by
writing
\[
\begin{cases}
\ \Sigma_\alpha \hsy (\jS) 
\ = \ 
\bigg[
\bigcup\limits_{\beta < \alpha} \ \Pi_\beta\hsy  (\jS)
\bigg]_\sigma\\[18pt]
\ \Pi_\beta \hsy (\jS)
\ = \ 
\bigg[
\bigcup\limits_{\beta < \alpha} \ \Sigma_\beta \hsy (\jS)
\bigg]_\delta
\end{cases}
\qquad (\alpha < \Omega).
\]
If $\beta < \alpha$, then 
\[
\begin{cases}
\ \Sigma_\beta \hsy (\jS) 
\subsetx
\ \Pi_\alpha \hsy (\jS)
\\[8pt]
\ \Pi_\beta  \hsy(\jS)
\subsetx
\Sigma_\alpha \hsy (\jS)
\end{cases}
,
\]
and if $\alpha > 1$ and $\beta < \alpha$, then 
\[
\begin{cases}
\ \Sigma_\beta\hsy  (\jS) 
\subsetx
\ \Sigma_\alpha \hsy (\jS)
\\[8pt]
\ \Pi_\beta\hsy  (\jS)
\subsetx
\Pi_\alpha\hsy  (\jS)
\end{cases}
.
\]
Therefore
\[
\bigcup\limits_{\alpha < \Omega} \ \Sigma_\alpha \hsy (\jS)
\ = \ 
\bigcup\limits_{\alpha < \Omega} \ \Pi_\alpha \hsy (\jS)
\]

\noindent
the $\sigma$-algebra generated by $\jS$.  
Note too that
\[
\begin{cases}
\ \Sigma_\alpha (\jS) \quad \text{is a $\sigma$-lattice}
\\[8pt]
\ \Pi_\alpha (\jS)      \quad \text{is a $\delta$-lattice}
\end{cases}
\]
with 
\[
\begin{cases}
\ \Sigma_\alpha (\jS) 
\ = \
[\Pi_\alpha (\jS)]_\tc
\\[8pt]
\ \Pi_\alpha (\jS)
\ = \
[\Sigma_\alpha (\jS)]_\tc
\end{cases}
.
\]
It is customary to refer to the sets in 
\[
\begin{cases}
\ \Sigma_\alpha (\jS) 
\\[8pt]
\ \Pi_\alpha (\jS)
\end{cases}
\quad \text{as} \quad 
\begin{cases}
\ \text{\un{additive of class $\alpha$}}
\\[8pt]
\ \text{\un{multiplicative of class $\alpha$}}
\end{cases}
,
\]
the sets in the intersection
\[
\Delta_\alpha(\jS)
\ = \ 
\Sigma_\alpha (\jS)  \cap \Pi_\alpha  (\jS) 
\]
then being \un{ambiguous of class $\alpha$}.  
Evidently, $\Delta_\alpha(\jS)$ is an algebra.
\\[-.25cm]

Our hierarchy may be visualized as follows: 
\\[-.25cm]

\[
\begin{tikzcd}[sep=0.25em]
&\Sigma_1 (\jS)
&&\Sigma_2 (\jS)
\\
&{\hspace{-2cm}\rotatebox{45}{$\subset$}}
&{\hspace{-2cm}\rotatebox{-45}{\raisebox{.35cm}{$\subset$}}}
&{\hspace{-2cm}\rotatebox{45}{$\subset$}}
\\
\Delta_1(\jS)
&&\Delta_2(\jS)
&&\cdots
\\
&{\hspace{-2cm}\rotatebox{-45}{$\subset$}}
&{\hspace{-2cm}\rotatebox{45}{\raisebox{-.35cm}{$\subset$}}}
&{\hspace{-2cm}\rotatebox{-45}{$\subset$}}
\\
&\Pi_1 (\jS)
&&\Pi_2 (\jS)
\end{tikzcd}
\]
\\[-0.5cm]

[Note: \ 
It need not be true, of course, that
\[
\begin{cases}
\ \Sigma_0 (\jS) 
\subsetx \Sigma_1 (\jS) 
\\[8pt]
\ \Pi_0 (\jS)
\subsetx \Pi_1 (\jS) 
\end{cases}
.
\]
However, the assumption $\jS \subsetx \jS_{\tc \sigma}$ would guarantee this.]
\\[-.25cm]

\textbf{Examples} \ 
\\[-.5cm]

\qquad (1) \quad
Let $\XX$ be a topological space $-$then the preceding considerations are applicable with 
$\jS = \fG$, the associated $\sigma$-algebra being $\jBo(\XX)$.
\\[-.25cm]

\qquad (2) \quad
Let $\XX$ be a topological space $-$then the preceding considerations are applicable with 
$\jS = \gZ_\tc$, the associated $\sigma$-algebra being $\jBa(\XX)$.
\\[-.25cm]

For each $\alpha < \Omega$, put
\[
\Lambda_\alpha  (\jS)
\ = \ 
u_\alpha (\jS).
\]

\textbf{Lemma 6} \ 
\un{Suppose that} $\jS \subsetx \jS_{\tc \hsy \sigma}$ $-$then 

\[
\begin{cases}
\ \Lambda_{2 \hsy n}  (\jS) \hspace{0.4cm} = \ \Sigma_{2 \hsy n} (\jS)
\hspace{.9cm} (n = 0, 1, 2, \ldots)
\\[8pt]
\ \Lambda_{2  n + 1}  (\jS) \ = \ \Pi_{2  n + 1} (\jS)
\hspace{0.5cm} (n = 0, 1, 2, \ldots)
\end{cases}
\]
\un{and}
\[
\Lambda_\alpha (\jS) 
\ = \ 
\Delta_{\alpha + 1} (\jS)
\quad (\omega \leq \alpha < \Omega).
\]
\\[-1cm]

[Note: \ 
There is also a dual result whereby, working with $\jS_\tc$ (instead of $\jS$), 
one picks off 
$\Sigma_\text{odd} (\jS)$
and 
$\Pi_\text{even} (\jS)$, 
the contention as regards the $\Delta_{\alpha + 1} (\jS)$ being unchanged.]
\\[-.25cm]

The proof, while not difficult, is a bit lengthy.
\\[-.25cm]

We shall first deal with the case of finite $n$.  
If $n = 0$, then, by definition, 
$\Lambda_0 (\jS) = \jS = \Sigma_0 (\jS) $.
On the other hand, if $S \in \Lambda_1 (\jS)$, then $S = \lim \ S_i$, where $\{S_i\}$ is a sequence in $\jS$.
In particular: 
\[
S 
\ = \ 
\varlimsup \ S_i 
\ = \ 
\bigcap\limits_{i = 1}^\infty \ 
\left(
\bigcup\limits_{m = i}^\infty \ S_m
\right).
\]
Because $\jS$ is a $\sigma$-lattice, 
$\ds\bigcup\limits_{m = i}^\infty \ S_m \in \jS$ $\forall \ i$, 
hence $S \in \Pi_1(\jS)$.  
To go the other way, take an $S \in \Pi_1(\jS)$ $-$then
\[
S 
\ = \ \bigcap \ S_i 
\ = \ 
\lim \ (S_1 \cap \cdots \cap S_i) 
\quad (S_i \in \jS)
\]
belongs to $\Lambda_1 (\jS)$.  
Proceeding by induction, suppose now that $n \geq 0$ and that our assertion is true for $n$  
$-$then it must be shown that
\[
\begin{cases}
\ \Lambda_{2 \hsy n + 2}  (\jS) \ = \ \Sigma_{2 \hsy n + 2} (\jS)
\\[8pt]
\ \Lambda_{2  n + 3}  (\jS) \ = \ \Sigma_{2  n + 3} (\jS)
\end{cases}
.
\]
Let us consider the first of these relations, the argument for the second being similar.  
If $S \in \Lambda_{2 \hsy n + 2}  (\jS)$, then $S = \lim \ S_i$, where $\{S_i\}$ is a sequence in 
\[
\Lambda_0  (\jS) \cup \cdots \cup \Lambda_{2 \hsy n + 1}  (\jS)
\ = \ 
\Lambda_{2 \hsy n + 1}  (\jS)
\ = \ 
\Pi_{2 \hsy n + 1}  (\jS)
\quad \text{(by induction).}
\] 
In particular: 

\[
S 
\ = \ 
\varliminf \ S_i 
\ = \ 
\bigcup\limits_{i = 1}^\infty \ 
\left(
\bigcap\limits_{m = i}^\infty \ 
S_m
\right).
\] 
Because $\Pi_{2 n + 1} (\jS)$ is a $\delta$-lattice, 
$\ds\bigcap\limits_{m = i}^\infty \ S_m \in \Pi_{2 \hsy n + 2}  (\jS)$ $\forall \ i$, hence 
$S \in \Sigma_{2 \hsy n + 1}  (\jS)$.  
To go the other way, take an $S \in \Sigma_{2 n + 2} (\jS)$ $-$then
\[
S 
\ = \ 
\bigcup \ S_i
\ = \ 
\lim \  (S_1 \cup \ldots \cup S_i),
\] 
where

\[
S_i \in 
\bigcup\limits_{j < 2 n + 2} \ \Pi_j(\jS) 
\ = \ 
\Pi_{2 n + 1} (\jS)
\ = \ 
\Lambda_{2 \hsy n + 1}  (\jS)
\quad \text{(by induction),}
\] 
that is, $S$ belongs to $\Lambda_{2 \hsy n + 2}  (\jS)$.
\\[-.25cm]

Passing to the transfinite assertion, suppose initially that $\alpha = \omega$.
If $S \in \Lambda_\omega  (\jS)$, then $S = \  \lim S_i$, where $S_i \in \Lambda_{m_i}  (\jS)$, say.
The claim is that

\[
\begin{cases}
\ S \in \Sigma_{\omega + 1} (\jS) = [\Pi_\omega (\jS)]_\sigma 
\\[8pt]
\ S \in \Pi_{\omega + 1} (\jS) = [\Sigma_\omega (\jS)]_\delta 
\end{cases}
.
\]
This, however, is immediate provided we take into account the relations

\[
\begin{cases}
\ 
\Lambda_{2 \hsy n}  (\jS) 
\ = \ 
\Sigma_{2 \hsy n}  (\jS) 
\ \subset \ 
\Pi_{2 \hsy n + 1}  (\jS) 
\\[8pt]
\ 
\Lambda_{2 \hsy n + 1}  (\jS) 
\ = \ 
\Pi_{2 \hsy n + 1}  (\jS) 
\ \subset \ 
\Sigma_{2 \hsy n + 2}  (\jS) 
\end{cases}
\]
and the fact that here
\[
\bigcap\limits_{i = 1}^\infty \ 
\left(\bigcup\limits_{m = i}^\infty \ S_m\right)
\ = \ 
\bigcup\limits_{i = 1}^\infty \ 
\left(\bigcap\limits_{m = i}^\infty \ S_m\right).
\] 
The other direction is slightly more complicated.  
Take an
$S \in \Delta_{\omega + 1} (\jS)$ $-$then there exist sequences 
$\{S_{i, j}^\prime\}$, 
$\{S_{i, j}^{\prime\prime}\}$
with 
\[
\begin{cases}
\ S_{i, j}^\prime \in \Sigma_0 (\jS) \cup \Sigma_1 (\jS) \cup \cdots
\\[8pt]
\ S_{i, j}^{\prime\prime} \in \Pi_0 (\jS) \cup \Pi_1 (\jS) \cup \cdots
\end{cases}
\]
such that 
\[
\begin{cases}
\ S \ = \ \bigcup\limits_i \ \bigcap\limits_j \ S_{i, j}^\prime
\\[18pt]
\ S \ = \ \bigcap\limits_i \ \bigcup\limits_j \ S_{i, j}^{\prime\prime}
\end{cases}
.
\]
Evidently, without loss of generality, it can be assumed that 
\[
S_{i, j}^\prime \ \supset \ S_{i, j+ 1}^\prime,  
\quad 
S_{i, j}^{\prime\prime} \ \subset \  S_{i, j+1}^{\prime\prime}.
\]
Consequently (cf. Prob. I (\S1)), 
\[
S \ = \ 
\lim(
(S_{1, j}^\prime \cap S_{1, j}^{\prime\prime}) 
\cup
(S_{2, j}^\prime \cap S_{1, j}^{\prime\prime} \cap S_{2, j}^{\prime\prime} ) 
\cup \cdots \cup 
(S_{j, j}^\prime \cap S_{1, j}^{\prime\prime} \cap \cdots \cap S_{j, j}^{\prime\prime})
).
\]
Each term inside the limit sign belongs to 
$\Lambda_0 (\jS) \cup \Lambda_1 (\jS) \cup \cdots$, 
implying, therefore, that $S \in \Lambda_\omega (\jS)$.  
Proceeding by transfinite induction, suppose for now that $\alpha$ is $> \omega$ and $< \Omega$ 
and that our assertion is true for $\omega \leq \beta < \alpha$.  
If $S \in \Lambda_\alpha (\jS)$, then $S = \lim \ S_i$, where
$S_i \in \Lambda_{\alpha_i} (\jS)$, say $(\omega \leq \alpha_i < \alpha)$.  
Because
\[
\Lambda_{\alpha_i} (\jS) 
\ = \ 
\Lambda_{\alpha_i + 1} (\jS)
\quad \text{(by induction)},
\]
and $\alpha_i + 1 \leq \alpha$, each $S_i$ belongs to $\Delta_\alpha (\jS)$, so the usual 
$\varlimsup, \  \varliminf$ 
representation forces $S$ into $\Delta_{\alpha + 1} (\jS)$.  
To finish up, take an $S \in \Delta_{\alpha + 1} (\jS)$ $-$then as above, there exist sequences 
$\{S_{i, j}^\prime\}$, 
$\{S_{i, j}^{\prime\prime}\}$
with 
\[
\begin{cases}
\ S_{i, j}^\prime \in \Sigma_{\xi_{i, j}} (\jS)
\hspace{0.5cm} (\omega \leq \xi_{i, j} < \alpha)
\\[8pt]
\ S_{i, j}^{\prime\prime} \in \Pi_{\eta_{i, j}} (\jS)
\hspace{0.5cm} (\omega \leq \eta_{i, j} < \alpha)
\end{cases}
\]
such that 
\[
\begin{cases}
\ S \ = \ \bigcup\limits_i \ \bigcap\limits_j \ S_{i, j}^\prime
\\[18pt]
\ S \ = \ \bigcap\limits_i \ \bigcup\limits_j \ S_{i, j}^{\prime\prime}
\end{cases}
,
\]
it not being restrictive to assume that 
\[
S_{i, j}^\prime
\ \supset \ 
S_{i, j + 1}^\prime, 
\quad 
S_{i, j}^{\prime\prime}
\ \subset \ 
S_{i, j+1}^{\prime\prime}.
\]
Let us distinguish two cases.
\\[-.25cm]

\qquad (A) \quad
$\alpha$ is an ordinal of the first kind, i.e., $\alpha$ possesses an immediate predecessor, say 
$\alpha = \beta + 1$ $-$then 
\[
\begin{cases}
\ \omega \ \leq \ \xi_{i, j} \ \leq \beta
\\[8pt]
\ \omega \ \leq \ \eta_{i, j} \ \leq \beta
\end{cases}
\quad \implies \ 
S_{i, j}^\prime, \hsx S_{i, j}^{\prime\prime} \in \Delta_{\beta + 1} (\jS)
\]

\noindent
$\implies$
\[
(S_{1, j}^\prime \cap S_{1, j}^{\prime\prime}) 
\cup 
(S_{2, j}^\prime \cap S_{1, j}^{\prime\prime} \cap S_{2, j}^{\prime\prime}) 
\cup \cdots \cup
(S_{j, j}^\prime \cap S_{1, j}^{\prime\prime} \cap \cdots \cap S_{j, j}^{\prime\prime}) 
\in 
\Delta_{\beta + 1} (\jS).
\]
But
\[
\Delta_{\alpha} (\jS)
\ = \ 
\Delta_{\beta + 1} (\jS)
\ = \ 
\Delta_{\beta} (\jS)
\quad \text{(by induction)},
\]
and so, $S \in \Lambda_\alpha (\jS)$, as desired.
\\[-.25cm]

\qquad (B) \quad 
$\alpha$ is an ordinal of the second kind, i.e., $\alpha$ possesses no immediate predecessor, 
thus is a limit ordinal, say $\alpha = \lambda$.  
\\[-.5cm]

\noindent
Put
\[
\zeta_j \ = \ \sup \ 
\begin{cases}
\ \xi_{1,j}, \ldots, \xi_{j,j}
\\[8pt]
\ \eta_{1,j}, \ldots, \eta_{j,j}
\end{cases}
.
\]
Then $\zeta_j  \ < \  \lambda$ $\forall \ j$
\\[-.5cm]

\noindent
$\implies$
\[
(S_{1, j}^\prime \cap S_{1, j}^{\prime\prime}) 
\cup 
(S_{2, j}^\prime \cap S_{1, j}^{\prime\prime} \cap S_{2, j}^{\prime\prime}) 
\cup \cdots \cup
(S_{j, j}^\prime \cap S_{1, j}^{\prime\prime} \cap \cdots \cap S_{j, j}^{\prime\prime}) 
\in 
\Delta_{\zeta_j + 1}  (\jS).
\]
However, as $\lambda$ is a limit ordinal, $\zeta_j + 1 < \lambda$ $\forall \ j$, hence
\[
\Delta_{\zeta_j + 1} (\jS) 
\ = \ 
\Lambda_{\zeta_j} (\jS) 
\quad \text{(by induction)}
\]
from which it follows that $S \in \Lambda_\alpha (\jS)$, as desired.
\\[-.5cm]

The proof of Lemma 6 is therefore complete.
\\[-.5cm]

[Note: \ 
It must be stressed that the assumption $\jS \subset \jS_{\tc \sigma}$ is crucial for the validity of this result.]
\\[-.25cm]

Suppose still that $\jS \subset \jS_{\tc \sigma}$ $-$then, thanks to Lemma 5,
\[
\jS_{\tB_\td}
\ = \ 
\sigma\dash\jAlg (\jS).
\]
Furthermore, $\forall \ \alpha > 0$:

\[
\begin{cases}
\jB^{(\alpha)} (\Delta_1 (\jS) ) 
\ = \ 
\Sigma_\alpha (\jS) 
\\[8pt]
\jB_{(\alpha)} (\Delta_1 (\jS) ) 
\ = \ 
\Pi_\alpha (\jS) 
\end{cases}
.
\]
\\[-1cm]

We shall conclude this \S \  with a brief discussion of  relativization and localization.
\\[-.25cm]

Suppose that $\star$ is extensionally attainable.  
Let $X_0$ be a subset of $\XX$ $-$then, given any nonempty $\jS$, we ask: 
Is 
\[
\star(\tr_{X_0} (\star(\jS))) 
\ = \ 
\tr_{X_0} (\star(\jS)) ?
\]
Generally, this need not be the case.  
But it will be true under the following assumptions: 
\\[-.25cm]

\qquad (1) \quad 
$
\star(\tr_{X_0} (\star(\jS))) 
\ = \ 
\tr_{X_0} (\star(\jS));
$
\\[-.25cm]

\qquad (2) \quad 
$
\{S \subset \XX \hsy : \hsy S \cap X_0 \in \star(\tr_{X_0} (\jS))\}
$
is a $\star$-class.
\\[-.25cm]

\noindent
Indeed, from (1) we get that
\[
\tr_{X_0} (\star(\jS)) 
\supset
\star(\tr_{X_0} (\jS))
\]
whereas from (2) we get that 
\[
\tr_{X_0} (\star(\jS)) 
\subsetx
\star (\tr_{X_0} (\jS)).
\]
Evidently, the properties
\[
\begin{cases}
\ \text{? \ is a lattice}\\[4pt]
\ \text{? \ is a ring ($\sigma$-ring, $\delta$-ring)}
\end{cases}
\]
are instances where conditions (1) and (2) are met.
\\[-.25cm]

\textbf{Example} \ 
Borel sets relativize.  
Thus, suppose that $\XX$ is a topological space with ambient topology $\jT$.  
Let $X_0$ be a subset of $\XX$ $-$then, by definition, the class $\tr_{X_0} (\jT)$ 
is the relative topology on $X_0$, and, by the above, we have
\[
\tr_{X_0} (\jBo(\XX)) 
\ = \ 
\sigma\dash \jRin (\tr_{X_0} (\jT))) 
\ = \ 
\jBo(X_0).
\]
\\[-1cm]

\textbf{Example} \ 
Baire sets need not relativize.  
To produce an example, we shall work within the Stone-\v Cech compactification 
$\beta \N$ of $\N$.  
Choose, as is possible, a class $\{S_i\}$ of $\fc$ infinite subsets of $\N$ such that 
\[
\card(S_i \cap S_j) 
\ < \ 
+\infty
\quad \forall \ i \neq j.
\]
This done, call $\bar{S}_i$ the closure of $S_i$ in $\beta \N$ $-$then the
$\bar{S}_i - \N$ are pairwise disjoint, open and closed subsets of $\beta \N - \N$.  
Put
\[
S 
\ = \ 
\bigcup\limits_i \ 
(\bar{S}_i - \N)
\]
and consider the subspace $\XX= \N \cup S$ of $\beta \N$.  
Since $\beta \N - \N$ is a zero set in $\beta \N$, $S$ is a zero set, hence a Baire set in $\XX$.  
Now
\[
\card (\jBa(\XX)) 
\ \leq \ 
2^{\aleph_0}, 
\]
$\XX$ being separable.  
On the other hand, it is clear that

\[
\card (\jBa(S)) 
\ \geq \ 
2^{2^{\aleph_0}}.
\]
Accordingly, not every Baire set of $S$ is a Baire set of $X$, and so here Baire sets do not relativize. 
\\[-.25cm]

Under certain conditions, however, Baire sets will relativize.  
Thus, suppose that $\XX$ is a topological space $-$then a subspace $X_0$ of $\XX$ is said to be 
\un{$\gZ$-embedded} in $\XX$ if $\forall$ zero set $Z_0$ in $X_0$ 
$\exists$ a zero set $Z$ in $\XX$ 
such that 
$Z_0 = Z \cap X_0$, i.e., if, in an obvious notation, 
\[
\tr_{X_0} (\gZ) 
\ = \ 
\gZ_0.
\]
But then
\\[-.25cm]

$\tr_{X_0} (\jBa(\XX)) 
\ = \ 
\sigma\dash \jRin (\tr_{X_0} (\gZ)) 
\ = \ 
\jBa(X_0).
$
\\[-.25cm]

For orientation, let us consider some specific instances of $\gZ$-embeddings.
\\[-.25cm]

(1) \quad 
Let $\XX$ be a completely regular, Hausdorff topological space $-$then $\XX$ is 
$\gZ$-embedded in its Stone-\v Cech compactification $\beta \XX$.
\\[-.5cm]

[This follows from the definitions.]
\\[-.25cm]

(2) \quad 
Let $\XX$ be a normal topological space $-$then every closed subset $X_0$ of $\XX$ is $\gZ$-embedded in $\XX$.
\\[-.5cm]

[Bear in mind the Tietze extension theorem.]
\\[-.25cm]

(3) \quad 
Let $\XX$ be a compact Hausdorff space $-$then every Baire set $X_0$ of $\XX$  is $\gZ$-embedded in $\XX$.
\\[-.5cm]

[In fact, $X_0$ is necessarily Lindel\"of.]
\\[-.25cm]

[Note: \ 
A systematic discussion of $\gZ$-embedding may be found in 
R. Blair and A. Hager\footnote[2]{\vspace{.11 cm}\un{ Math. Z.} \textbf{136} (1974), pp. 41-52.}
.
\\[-.25cm]

Let $\jS$ be a nonempty subset of $\jP(\XX)$ $-$then by the 
\un{localization} $\jS_\loc$ of $\jS$ we mean the class consisting of all $X_0 \subset \XX$ 
for which 
\[
\tr_{X_0} (\jS)
\subset 
\jS.
\]
Obviously, $\XX \in \jS_\loc$, so $\jS_\loc$ is nonemtpy.  
In addition, if $\jS$ is multiplicative, then $\jS \subset \jS_\loc$.
\\[-.25cm]

Suppose that $\jS$ is a ring ($\sigma$-ring, $\delta$-ring) $-$then $\jS_\loc$ is an 
algebra ($\sigma$-algebra, $\delta$-algebra).
\\[-.25cm]

\textbf{Example} \ 
Let $\XX$ be a Hausdorff topological space.  
Let $\jK$ be the class of all compact subsets of $\XX$ $-$then it is easy to see that 
\[
S \in [\jBo_\tb (\XX) ]_\loc 
\quad \text{iff} \quad 
S \cap \K \in \jBo (K) 
\quad \forall \ K \in \jK.
\]
Consequently, 
\[
\jBo (\XX) 
\subset 
[\jBo_\tb (\XX) ]_\loc, 
\]
the containment being strict in general (cf. Exer. 40), but there being coincidence if, e.g., $\XX$ is $\sigma$-compact.
\\[-.25cm]


Localization need not commute with generation. 
\\[-.25cm]

\textbf{Example} \ 
In general, 
\[
\sigma\dash \jRin (\jS_\loc)
\ \neq  \ 
[\sigma\dash \jRin (\jS)]_\loc \hsy.
\]
Thus, take $\XX = \N$ and set $\jS = \{\{n\} : n \in \N\}$ $-$then 
$\jS_\loc  = \{\XX\}$, hence, in this case,
\begin{align*}
\text{$\sigma$-$\jRin (\jS_\loc)$} \ 
&=\ \{\emptyset, \XX\}
\\[11pt]
&\neq \ 
\jP (\XX) 
\\[11pt]
&=\ 
\text{$[\sigma$-$\jRin (\jS)]_\loc$}.
\end{align*}
\\[-.75cm]

\[
\textbf{\un{Notes and Remarks}}
\]

The term ``extensionally attainable'' has been borrowed from 
T. Hildebrandt\footnote[2]{\vspace{.11 cm}\textit{\un{Introduction to the Theory of Integration}},
Academic Press, New York, (1963).}.
If $\star$ is an extensionally attainable property, then some authors would refer to $\star (\jS)$ as the 
\un{$\star$-stabilization} of $\jS$.  
The generation of lattices and rings was discussed already by 
F. Hausdorff\footnote[3]{\vspace{.11 cm}\textit{\un{Grundz\"uge der Mengenlehre}}, Veit \& Comp., Leipzig, (1914).}.
The transfinite approach to operation B has its origins in 
E. Borel\footnote[4]{\vspace{.11 cm}\textit{\un{Lecons sur la Th\'eorie des Fonctions}}, Gauthier-Villars, Paris, (1898).}.
although this author evidently did not believe in transfinite numbers.  
The general formulation is due to 
F. Hausdorff (op. cit. pp. 304-306), 
further details and refinements being presented by him in 
F. Hausdorff\footnote[2]{\vspace{.11 cm}\textit{\un{Math. Ann.}}, \textbf{77} (1916), 241-256.} 
and later on in his famous 
F. Hausdorff\footnote[3]{\vspace{.11 cm}\textit{\un{Mengenlehre}}, Walter de Gruyter, Berlin, (1927).} 
The axiomatic approach to Borel sets in terms of a generated $\sigma$-ring was stressed
by 
W. Sierpi\'nski\footnote[4]{\vspace{.11 cm}\textit{\un{Bull. Acad. Sci. Cracovie}}, \textbf{A} (1918), pp. 29-34.} 
Lemmas 4 and 5 are results of Sierpi\'nski; cf. respectively
W. Sierpi\'nski\footnote[4]{\vspace{.11 cm}\textit{\un{Annales Soc. Polon. Math.}},  \textbf{6} (1927), pp. 50-53.}
and
W. Sierpi\'nski\footnote[4]{\vspace{.11 cm}\textit{\un{Fund. Math.}},  \textbf{12} (1928), pp. 206-210.}.
For an excellent account of the theory as it stood around 1930 and which is still very readable even now, consult 
H. Hahn\footnote[5]{\vspace{.11 cm}\textit{\un{Reelle Funktionen}},  Akademische Verlagsgesellschaft M.B.H., Leipzig, (1932), pp. 258-276.}.
Given $\alpha$ $(1 < \alpha < \Omega)$, the existence of a topological space $\XX_\alpha$ such that $\jF \subset \jG_\alpha$ was first 
noted by 
S. Willard\footnote[6]{\vspace{.11 cm}\textit{\un{Fund. Math.}}, \textbf{71} (1971), pp. 187-191.}.
The definition in the text of a Baire set is apparently due to 
E. Hewitt\footnote[7]{\vspace{.11 cm}\textit{\un{Fund. Math.}}, \textbf{37} (1950), pp. 161-187.}.
The reader is warned that while we consider the definitions in the text of Borel set and Baire set to be the most natural, 
other writers might use these terms for very different entities.  
E.g.: 
In some treatments, the Borel sets in a Hausdorff topological space are taken to be the $\sigma$-ring generated by the compact sets, 
the Baire sets then being the $\sigma$-ring generated by the compact$G_\delta$'s.  
The fact that $\jBa (\XX)$ can be produced from $\jZ$ by operation $\tB_d$ was established by 
J. Jayne\footnote[8]{\vspace{.11 cm}\textit{\un{Mathematika}}, \textbf{24} (1977), 241-256.}.  
In this connection, it should be kept in mind that there is a theorem in general topology which says that no nonempty, open subset 
of a connected compact Hausdorff space $\XX$ can be written as a countable disjoint union of nonempty, closed subsets of $\XX$; 
cf. 
K. Kuratowski\footnote[1]{\vspace{.11 cm}\un{Topology} Vol II, Academic Press, New York, (1968) p. 173.}. 
The origin of the notation 
$\Sigma_\alpha (\jS)$, $\Pi_\alpha (\jS)$ lies in recursive function theory; 
it was introduced by 
J. Addison\footnote[2]{\vspace{.11 cm}\un{Fund. Math.}, \textbf{46} (1959), pp. 123-135.}.
The procedure itself, however, can be traced back to 
F. Hausdorff\footnote[3]{\vspace{.11 cm}\un{Math. Z.}, \textbf{5} (1919), pp. 292-309.}.
Emphasis on the $\Lambda_\alpha (\jS)$ was placed by 
Ch. de la Vall\'ee Poussin\footnote[4]{\vspace{.11 cm}\un{Int\'egrales de Lebesgue}, \un{Fonctions d'Ensembles}, \un{Classes de Baire}Gauthier-Villars, Paris, (1916), p. 37.}.
The connection between the two, 
i.e., Lemma 6, was found by 
W. Sierpi\'nski\footnote[7]{\vspace{.11 cm}\un{Fund. Math.} \textbf{19} (1932), pp. 257-264.}; 
see also 
J. Albuquerque\footnote[8]{\vspace{.11 cm}\un{Portugual. Math.} \textbf{4} (1943-1945), pp. 161-198, pp. 217-224.}. 
The notion of localization appears explicitly in 
I. Segal\footnote[9]{\vspace{.11 cm}\un{Amer. J. Math.} \textbf{73} (1951), pp. 275-313.}, 
although it is implicit in earlier writings.  
N. Dinculeanu\footnote[1]{\vspace{.11 cm}\un{Vector Measures}, Pergamon Press, London, (1967).}
defines the Borel sets in a locally compact Hausdorff space as the localization of the $\delta$-ring generated by the compact sets, 
Baire sets being defined similarly as the localization of the $\delta$-ring generated by the compact $G_\delta$'s.
\chapter{
$\boldsymbol{\S}$\textbf{6}.\quad  Exercises}
\setlength\parindent{2em}
\setcounter{theoremn}{0}
\renewcommand{\thepage}{\S6-\text{E-}\arabic{page}}


\vspace{.25cm}
\qquad 
(1) \quad 
Let $\star$ be the property:  ? is a topology.  
Verify that $\star$ is extensionally attainable.  
Given any nonempty $\jS$, $\star (\jS)$ is called the topology generated by $\jS$ and is denoted by $\jTop (\jS)$.  
Verify that $\jTop (\jS) = \jS_{\td \hsy \Sigma}$ with, if necessary, $\emptyset$ and $\XX$ adjoined.
\\

(2) \quad
Given a ring $\jS$, a ring with unit containing $\jS$ is the class
\[
\widehat{\jS} 
\ = \ 
\jRin (\jS, \big\{ \hsy\bigcup \ \jS \big\}).
\]
If $\jT$ is a ring with unit containing $\jS$, then $\bigcup \ \jT \ \supset \bigcup \ \jS$.  
Nevertheless, show by example that there exists a ring $\jS$ and a ring with unit $\jT$ such that 
\\[-.25cm]

\[
\begin{cases}
\ \jT \supset \jS
\\[8pt]
\ \jT \not\supset \ \widehat{\jS} 
\end{cases}
.
\]
\\[-.75cm]

[Take $\XX = [0,2]$.  Let $\jS$ be the class consisting of all first category subsets of $[0,1]$.  
Consider
\[
\jT 
\ = \ 
\jRin (\jS, \{[0,2]\}) \hsx .]
\]
\\[-.75cm]

(3) \quad
Let $\jS$ be nonempty $-$then we have: 
\\[-.25cm]

\qquad (i) \quad 
$\jRin(\jS) = \jS_{\trx \td \ts} = \jS_{\trx \ts \td}$ ;
\\[-.25cm]

\qquad (ii) \quad 
$\jRin(\jS) = \jS_{\td \ts \trx \ts} = \jS_{\ts \td \trx \ts}$ ;
\\[-.25cm]

\qquad (iii) \quad 
$\jRin(\jS) = \jS_{\trx \ts \trx \ts}$.
\\[-.25cm]

Show by example that $\jS_{\trx \ts \trx \ts} \neq  \jS_{\trx \ts \trx}$ in general. 
\\[-.25cm]

[Take $\XX = \{1, 2, 3, 4, 5\}$ and let $\jS = \{\{2, 4\}, \{1, 2, 3\}, \{1, 4, 5\}\}$.]
\\


(4) \quad 
Let $\star$ be the property: 
$? = ?_\trx$ and $? = ?_\Sigma$. 
Verify that $\star$ is extensionally attainable.  
Given any nonempty $\jS$, show that 
\[
\star (\jS) 
\ = \ 
\jS_{\trx \hsy \Sigma \hsy \trx \hsy \Sigma}.
\]

[It is enough to prove that
\[
\jS_{\trx \hsy \Sigma \hsy \trx \hsy \Sigma \hsy \trx}
\ = \ 
\jS_{\trx \hsy \Sigma \hsy \trx \hsy \Sigma}.
\]
Incidentally, observe that $\Sigma$ cannot, in general, be replaced by $\sigma$ here; 
on the other hand, in view of Exer. 3 (iii), the substitution of s for $\Sigma$ does lead to a true statement.]
\\

(5)  \quad 
True or False? \ Suppose that $\emptyset \in \jS$, $\jS = \jS_\td$, and $\jS_\ts = \jRin (\jS)$ $-$then 
$\jS$ is a semiring.
\\[-.5cm]

[Compare with Exer. 5 (\S4).]
\\

(6) \quad 
Let $\XX$ be a topological space $-$then the ring generated by the open subsets of $\XX$ is called the class of 
\unz{constructible sets} in $\XX$.  
Verify that $S \subset \XX$ is constructible iff $S$ can be written as a finite union of locally closed subsets of $\XX$.
\\

(7) \quad 
Let $\jS$ be nonempty $-$then $\jRin (\jS)$ $(\sigma\dash\jRin (\jS))$ is the union of the rings ($\sigma$-rings) 
generated by the subsets of $\jS$ of cardinality $< \aleph_0$ $(\leq \aleph_0)$.
\\

(8) \quad 
Let $\jS$ be nonempty $-$then every set in $\jRin (\jS)$ $(\sigma\dash\jRin (\jS))$ can be covered by a finite (countable) 
union of sets in $\jS$.
\\[-.5cm]

[The class of all sets which can be covered  by a finite (countable) union of sets in $\jS$ is a ring ($\sigma$-ring).]
\\

(9) \quad 
Let $\XX$ be a nonempty set.  Suppose that $\jS$ is a $\sigma$-algebra in $\XX$ admitting a generating subclass $\jS_0$ 
of cardinality $\leq \aleph_0$ with the property that for all $x \neq y$ there exists an $S_0 \in \jS_0$ such that either 
$x \in S_0$ and $y \notin S_0$ or $x \notin S_0$ and $y \in S_0$. 
Under these conditions, prove that $\XX$ can be equipped with the structure of a separable metric space in which the Borel 
sets are precisely the elements of $\jS$.
\\[-.5cm]

[Let $\jS_0 = \{S_1, S_2, \ldots\}$ be an enumeration of $\jS_0$.  
Consider the metric $d$ defined by the rule
\[
d(x, y) 
\ = \ 
\sum \ 
\bigg(
\hsx
\frac{|\chisubSi (x) - \chisubSi (y)|}{2^i }
\hsx
\bigg)
.]
\]
\\[-.75cm]

(10) \quad 
Let $\XX = [0, \Omega]$, equipped with the order topology $-$then the Borel sets in $\XX$ consist of those subsets 
$S$ of $\XX$ such that either $S$ or $\XX - S$ contains an unbounded, closed subset of $[0, \Omega[$.  
Is every subset of $\XX$ a Borel set?
\\[-.5cm]

[The class of unbounded, closed subsets of $[0, \Omega[$ is closed under countable intersections; 
accordingly, the class in question is a $\sigma$-ring containing the Borel sets.  
To obtain equality, let $S$ be an unbounded, closed subset of $[0, \Omega[$ $-$then it need only be shown that every subset 
$T$ of $\XX - S$ is Borel.  
There is no loss of generality in supposing that $0 \in S$, $\Omega \notin T$.  
Given $\alpha \in S$, let $\alpha^\prime$ be the first successor to $\alpha$ in $S$.  
Define a set-valued function $f$ on $S$ by the prescription
\[
f(\alpha) 
\ = \ \{\beta \in T \hsy : \hsy \alpha < \beta < \alpha^\prime\}.
\]
Then $f(S) = T$.  
For each $\alpha$ such that $f(\alpha) \neq \emptyset$, fix an enumeration $\{f(\alpha)_n\}$ of the elements of $f(\alpha)$.  
Write
\[
T_n 
\ = \ 
 \bigcup_{\alpha \in S} \ \{f(\alpha)_n\}.
\]
The $T_n$ are Borel and $T = \bigcup \ T_n$.]
\\

(11) \quad 
Let $\XX$ be a topological space $-$then every Borel set in $\XX$ has the property of Baire.
\\

(12) \quad 
Let $\XX$ be a metric space $-$then $\XX$ is separable iff 
$\forall \ \varepsilon > 0$, 
$\jBo (\XX)$ 
is generated by the open balls of radius $\leq \varepsilon$.  
Show by example that there exists a nonseparable metric space $\XX$ in which the open balls 
\[
\begin{cases}
\ \text{do generate $\jBo (\XX)$}
\\
\ \text{do not generate $\jBo (\XX)$}
\end{cases}
.
\]
\\[-.75cm]

(13) \quad 
Let $\XX$ be a topological space, all of whose points are closed; 
let $S$ be a discrete subspace of $\XX$ $-$then $S$ is a Borel subset of $\XX$.
\\[-.5cm]

[In fact, $S$ is constructible.]
\\

(14) \quad 
Let $\XX$ be a Hausdorff topological space $-$then
the $\sigma$-ring generated by the compact subsets of $\XX$ is, by definition, the class of 
\unz{$\sigma$-bounded Borel sets} in $\XX$.  
Justify this terminology by proving that a Borel set in $\XX$ is $\sigma$-bounded iff it is contained in a countable union of 
compact sets.  
Hence or otherwise, infer that if $\XX$ is
\[
\begin{cases}
\ \text{$\sigma$-compact}
\\
\ \text{compact}
\end{cases}
,
\]
then
\[
\begin{cases}
\ \jBo (\XX) = [\jBo_\tb (\XX)]_\sigma
\\[8pt]
\ \jBo (\XX) = \jBo_\tb (\XX)
\end{cases}
.
\]
\\[-.75cm]

(15)  \quad
Let $\XX = [0, \Omega[$,  equipped with the order topology.  
Characterize explicitly the elements of the $\delta$-ring of bounded Borel sets in $\XX$.
\\

(16) \quad 
Let $\XX$ be a Hausdorff topological space.  Give a transfinite description of $\jBo_\tb (\XX)$.
\\

(17) \quad 
Let $\XX$ be a Hausdorff topological space.
Let $X_0$ be a compact subset of $\XX$ $-$then the bounded Borel sets in $\XX$, 
when relativized to $X_0$, give the bounded Borel sets in $X_0$, i.e., 
\[
\tr_{X_0} (\jBo_\tb (\XX)) 
\ = \ 
\jBo_\tb (X_0).
\]
Is this true if $X_0$ is not compact?
\\

(18) \quad 
Let $\XX$ be a Hausdorff topological space.  
Let $\jK = \{K\}$ be a class of compact subsets of $\XX$ such that 
\[
\begin{cases}
\ \forall \ K_1, K_2 \in \jK, \ \exists \ K_3 \in \jK \ \text{st} \ 
\begin{cases}
\ K_1 \subset K_3
\\
\ K_2 \subset K_3
\end{cases}
\\[26pt]
\ \forall \ \text{compact} \ C \subset \XX, \exists \ K \in \jK  \ \text{st} \ C \subset K
\end{cases}
.
\]
Then 
\[
\jBo_\tb (\XX)
\ = \ 
\bigcup\limits_{K \in \jK} \ \jBo_\tb (K).
\]

[Show that the union in question is a $\delta$-ring.]
\\

(19) \quad 
True or False? \ Let $\XX$ be a Hausdorff topological space $-$then the bounded Borel sets in $\XX$ are precisely the 
relatively compact Borel sets in $\XX$.
\\

(20) \quad 
Let $\jS$ be a $\sigma$-ring in $\XX$; 
let $\jT$ be a $\sigma$-ring in $\YY$ $-$then
any $E \in \jS \ \ovs{\otimes} \ \jT$ has at most $\fc$ distinct horizontal or vertical sections.
\\[-.5cm]

[Fix $E \in \jS \ \ovs{\otimes} \ \jT$ $-$then there exist $\sigma$-rings 
$\jS_E \subset \jS$ and 
$\jT_E \subset \jT$
such that 
$E \in \jS_E \ \ovs{\otimes} \ \jT_E$
and such that both 
$\jS_E$ and 
$\jT_E$ 
are generated by no more than $\aleph_0$ elements (cf. Exer. 7).  
Owing to Lemma 2 (\S5), 
\[
\begin{cases}
\ E_x \in \jT_E \hspace{0.5cm} \forall  \ x \in \XX
\\[4pt]
\ E^y \in \jS_E \hspace{0.5cm} \forall  \ y \in \YY
\end{cases}
.
\]
On the other hand, 
\[
\begin{cases}
\ \card(\jS_E) \leq \fc
\\[4pt]
\ \card(\jT_E) \leq \fc
\end{cases}
\quad .]
\]
\\[-.5cm]

(21) \quad 
Let $\jS$ be a $\sigma$-ring in $\XX$.  Suppose that $\card(\XX) > \fc$ $-$then the diagonal $D$ in $\XX \times \XX$ 
does not belong to $\jS \ \ovs{\otimes} \ \jS$.
\\[-.5cm]

[This follows from Exer. 20.]
\\

(22) \quad 
Let $\XX$ and $\YY$ be Hausdorff topological spaces $-$then 
\[
\jBo (\XX) \ \ovs{\otimes} \ \jBo (\YY) 
\hsx \subset \  
\jBo (\XX \times \YY), 
\]
the containment being strict in general, but there being coincidence if the weight of $\XX$ and $\YY$ are both $\leq \aleph_0$.  
Does coincidence obtain if $\XX$ and $\YY$ are arbitrary Lindel\"of spaces?
\\[-.5cm]

[Note: \ 
Do Baire sets ``multiply''?  
While the answer is, of course, ``no'' in general, an important sufficient condition is this.  
Suppose that $\XX$ and $\YY$ are completely reguar, Hausdorff topological spaces for which $\XX \times \YY$ is 
$\jZ$-embedded in 
$\beta \XX \times \beta \YY$, 
the product of the Stone-\v Cech compactifications of $\XX$ and $\YY$ $-$then 
\[
\jBa (\XX) \ \ovs{\otimes} \ \jBa (\YY) 
\ = \  
\jBa (\XX \times \YY).
\]
For the details and further results, see 
R. Blair and 
A. Hager\footnote[2]{\vspace{.11 cm}\un{Set-Theoretic Topology}, Academic Press, New York, (1977), 47-72.}.]
\\

(23) \quad 
Let $\XX$ and $\YY$ be Hausdorff topological spaces $-$then 
\[
[\jBo_\tb (\XX)]_\sigma 
\ \ovs{\otimes} \
[\jBo_\tb (\YY)]_\sigma 
\hsx \subset \ 
[\jBo_\tb (\XX\times \YY)]_\sigma ,
\]
the containment being strict in general, but there being coincidence if the weights of $\XX$ and $\YY$ are both $\leq \aleph_0$.
Does coincidence obtain if $\XX$ and $\YY$ are arbitrary metric spaces?
\\

(24) \quad 
Take for $\XX$ the \unz{Sorgenfrey line $E$}, i.e., $E$ is the real line equipped with the topology generated by the $[a,b[$ $-$then 
\[
\jBo(E) 
\ = \ 
\jBo(\R)
\]
but
\[
\jBo(E \times E) 
\ \neq \ 
\jBo(\R\times \R).
\]

[To establish the second point, consider the line $L \hsy : \hsy x + y = 0$ $-$then, in the relative topology per 
$E \times E$, $L$ is discrete.  
Use now the fact that Borel sets relativize.]
Is
\[
\jBa(E) 
\ = \ 
\jBa(\R)?
\]
Is 
\[
\jBa(E \times E) 
\ = \ 
\jBa(\R\times \R)?
\]
\\[-.75cm]

(25)  \quad 
Given an example of an infinite class $\jS$ of subsets of $\R$ such that 
\[
\R \in \jS 
\quad \text{and} \quad 
\jS = \jS_\tB
\]
but such that $\jS$ is not a $\sigma$-algebra.
\\

(26) \quad 
Estimate the cardinality of $\jS_\tB$.  
Can the same be done of $\jS_{\tB_\td}$?
\\

(27) \quad 
True or False? \ Let $\jS$ be a ring.  
Suppose that for some limit ordinal $\lambda < \Omega$,
\[
\jS_\tB 
\ = \ 
\bigcup\limits_{\alpha < \lambda} \ \jB^{[\alpha]} (\jS).
\]
Then there is an $\alpha < \lambda$ such that 
\[
\jS_\tB 
\ = \ 
\jB^{[\alpha]} (\jS).
\]
\\[-1cm]

(28) \quad 
There exists a completely regular, nonnormal, Hausdorff topological space $\XX$ for which $\jF \subset \jG_\delta$.
\\[-.5cm]

[The classical example is the so-called \unz{Moore plane $\Gamma$}, 
i.e., $\Gamma$ is the closed upper half-plane $\{(x, y) \in \R^2 \hsy : \hsy y \geq 0\}$, topologized by specifying local open 
neighborhoods: 
The open neighborhoods of $(x, y)$ $(y > 0)$ are to be the usual open neighborhoods but the open 
neighborhoods of $(x, 0)$ are to be the sets $\{x\} \cup U$, 
where $U$ is an open disk in the upper half-plane tangent to the $x$-axis at $x$.]
\\

(29) \quad 
Let $\star$ be the property: \ 
$? = ?_{\sigma_\td}$ 
and 
$? = ?_\tc$.  
Verify that $\star$ is extensionally attainable.  
Given any nonempty $\jS$, we then write $\jS_{\tB_\tc}$ for $\star (\jS)$ and refer to $M_\star$ as \un{operation $\tB_\tc$}.  
Determine the properties of this operation.  
Show by example that $\jS_{\tB_\tc}$ need not coincide with $\sigma\dash\jAlg (\jS)$.  
Prove that 
\[
\jS_{\tB_\tc} 
\ = \ 
\sigma\dash\jAlg (\jS)
\]
iff
\[
\jS_\trx \hsx \subset \ \jS_{\tB_\tc} 
\quad \text{or} \quad 
\jS_\td \hsx \subset \ \jS_{\tB_\tc} .
\]

[So, in particular, if $\XX$ is a topological space, then
\[
\jBo (\XX) \ = \ 
\begin{cases}
\ \jF_{\tB_\tc}
\\[8pt]
\ \jG_{\tB_\tc}
\end{cases}
.]
\]

(30) \quad 
Let $\XX$ be a nonnormal, Hausdorff topological space $-$then $\jZ$ is properly contained in $\jF$.
\\

(31) \quad 
A compact Hausdorff space $\XX$ is 0-dimensional iff $\jZ_\sigma = \jZ_{\sigma_\td}$.
\\

(32) \quad 
Let $\XX = [0, \Omega]$, equipped with the order topology $-$then the Baire sets in $\XX$ consist of those subsets $S$ 
such that either
\[
\card(S) \leq \aleph_0 
\quad \text{or} \quad
\card(\XX - S) \leq \aleph_0.
\]
Thus, in this case, $\jBa (\XX)$ is strictly contained in $\jBo (\XX)$ (cf. Exer. 10). 
\\

(33) \quad 
Take for $\XX$ the real line topologized by specifying that the open sets are to be all sets of the form $U \cup V$, 
where $U$ is open in the usual topology and $V$ is any subset of $\PP = \XX - \Q$ $-$then
\[
\jBo (\XX) \ = \ 
\begin{cases}
\ \jF_{\tB_\td}
\\[8pt]
\ \jG_{\tB_\td}
\end{cases}
.
\]
However, $\XX$, while normal and Hausdorff, is not perfectly normal.  
Is \hsx $\jBa (\XX) = \jBo(\XX)$?
\\[-.25cm]

(34) \quad 
There exists a compact Hausdorff space $\XX$ for which $\jF_{\tB_\td} \neq \jF_\tB$.
\\[-.25cm]

[Let $A = D \cup \{\infty\}$ be the Alexandroff compactification of an uncountable discrete set $D$.  
Form the product $A \times \N$ and let $S$ be the set obtained by identifying 
$(\{\infty\}, n)$ $(n \in \N)$.  
Equip $S$ with the quotient topology $-$then $S$ is a completely regular, $\sigma$-compact, Hausdorff topological space.  
Let $\XX = \beta S$, the Stone-\v Cech compactification of $S$ $-$then 
$S \in \jF_\tB$ but $S \notin \jF_{\tB_\td}$.]
\\[-.25cm]

(35)  \quad 
Consider $\XX = [0,1]^{[0,1]}$ in the product topology.  
Is the subspace of all continuous $f : [0,1] \ra [0,1]$ a Borel (Baire) set in $\XX$?
\\[-.25cm]

(36) \quad 
Take $\XX = \R$ $-$then
\[
\PP \in \jF_{\sigma_d \delta \sigma_d \delta}.
\]
\\[-.75cm]

(37) \quad 
Let $\jS$ be nonempty $-$then we have: 

\qquad (i) \quad 
$\jS_\tB = \sigma\dash\jAlg (\jS) \quad \text{iff} \quad \jS_\tc \subset \jS_\tB$;
\\[-.25cm]

\qquad (ii) \quad 
$\jS_{\tB_\td} = \sigma\dash\jRin(\jS) \quad \text{iff} \quad \jS_\trx \subset \jS_{\tB_\td}$.
\\[-.25cm]

[Compare these statements with Lemmas 4 and 5.]
\\

(38) \quad 
Let $\jS$ be nonempty $-$then
\[
\delta\dash\jRin (\jS) 
\ = \ 
\bigcup\limits_{S \in \jRin (\jS)} \ \tr_S (\sigma\dash\jRin (\jS)).
\]
\\[-.75cm]

(39) \quad 
True or False? \ 
Let $\jS$ be a $\sigma$-ring in $\XX$; let $\jT$ be a $\sigma$-ring in $\YY$ $-$then
\[
\jS_\loc \ \ovs{\otimes}  \ \jT_\loc 
\ = \ 
(\jS \ \ovs{\otimes} \ \jT)_\loc .
\]
Retaining the given hypotheses, determine the validity of the relation
\[
\tr_{X_0 \times Y_0} (\jS \ \ovs{\otimes} \ \jT) 
\ = \ 
\tr_{X_0} (\jS) 
\ \ovs{\otimes}  \ 
\tr_{Y_0} (\jT).
\]
\\[-.75cm]

(40) \quad 
Let $\XX = [0, \Omega[$, equipped with the order topology $-$then
\[
[\jBo_\tb (\XX)]_\loc 
\ = \ 
\jP (\XX).
\]
Therefore, in this case, $\jBo (\XX)$ is strictly contained in $[\jBo_\tb (\XX)]_\loc $ 
(cf. Exer. 10 and 15).
\\[-.25cm]

[For a somewhat different example, discuss $\XX = \R \times \R$, where, in the first factor, 
$\R$ has the usual topology and, in the second factor, $\R$ has the discrete topology.]


\chapter{
$\boldsymbol{\S}$\textbf{6}.\quad  Problems}
\setlength\parindent{2em}
\setcounter{theoremn}{0}
\renewcommand{\thepage}{\S6-\text{P-}\arabic{page}}



\textbf{I} \ DYNKIN CLASSES
\\[-.25cm]

Let $\XX$ be a nonempty set; 
let $\jS$ be a nonempty subset of $\jP (\XX)$ $-$then $\jS$ is said to be a \un{Dynkin class} if 
$\jS = \jS_{\sigma_\td}$ and 
\[
S, \hsx T \in \jS, 
\ 
S \supset T 
\implies 
S - T \in \jS.
\]

Take $\star$ to be the property: ? is a Dynkin class.
It is clear that $\star$ is extensionally attainable.  
Given any nonempty $\jS$, we then call $\star (\jS)$ the Dynkin class generated by $\jS$ and denote it by 
$\jDyn (\jS)$.
\\[-.25cm]

Every $\sigma$-ring is a Dynkin class but a Dynkin class is a $\sigma$-ring iff it is closed under the formation of finite intersections. 
\\[-.5cm]

[For a simple example of a class which is a Dynkin class but is not a $\sigma$-ring, take 
$\XX = \{1, 2, 3, 4\}$ and consider 
\[
\jS
\ = \ 
\{\emptyset, \{1, 2\}, \{1, 3\}, \{2, 4\}, \{3, 4\}, \{1, 2, 3, 4\}\}.]
\]

If $\jS = \jS_\td$, then
\[
\text{$\sigma$-$\jRin (\jS)$} 
\ = \ 
\jDyn (\jS).
\]
\\[-.75cm]

\noindent
\un{Ref}
E. Dynkin\footnote[2]{\vspace{.11 cm}\textit{\un{Die Grundlagen der Theorie der Markoffschen Prozesse}}, 
Springer-Verlag, Berlin, (1961) pp. 1-2.}.
\\[-.25cm]

[Note: \ 
Results substantially the same as these were obtained many years earlier by 
W. Sierpi\'nski\footnote[3]{\vspace{.11 cm}\textit{\un{Fund. Math.}}, \textbf{12} (1928), pp. 206-210.}.]
\\[-.25cm]


There is a variant on the preceding theme which is sometimes useful.  
Consider the following properties of a nonempty $\jS$:

(1) \quad 
$\jS = \jS_{\ts_\td}$; 
\\[-.25cm]

(2) \quad 
$\forall \ S \in \jS$:  
$\forall \ S_i \in \jS$:  
\[
S_1, \hsx S_2, \ldots \subset S, 
\quad 
S_i \cap S_j = \emptyset 
\quad (i \neq j)
\]

\hspace{1.5cm} $\implies$
\[
\bigcup \ S_i \in \jS;
\]

(3) \quad
$S$, $T \in \jS$, $S \supset T \implies S - T \in \jS$ .
\\[-.25cm]

Let $\star$ be the conjunction of (1), (2), and (3) $-$then $\star$ is extensionally attainable 
and the above results on Dynkin classes can be carried over to this setting in the obvious way.  
In particular, observe that if $\jS = \jS_\td$, then $\star (\jS)$ is simply 
$\delta-\jRin (\jS)$.
\\[1cm]

\textbf{II} \ STABILITY OF SECTIONS
\\[-.25cm]

If
\[
\begin{cases}
\ \jS \subsetx \jP (\XX)\\
\ \jS \subsetx \jP (\YY)
\end{cases}
\quad 
\]
both contain $\emptyset$, then 
\[
\begin{cases}
\ \forall \ x : \jB^{(\alpha)} (\jS \hsx \boxtimes \hsx \jT)_x \subset \jB^{(\alpha)} (\jT) 
\quad \text{and} \quad
\jB_{(\alpha)} (\jS \hsx \boxtimes \hsx \jT)_x \subset \jB_{(\alpha)} (\jT) 
\\[8pt]
\ \forall \ y : \jB^{(\alpha)} (\jS \hsx \boxtimes \hsx \jT)^y \subset \jB^{(\alpha)} (\jS) 
\quad \text{and} \quad
\jB_{(\alpha)} (\jS \hsx \boxtimes \hsx \jT)^y \subset \jB_{(\alpha)} (\jS) 
\end{cases}
(\alpha < \Omega).
\]


[This follows by an easy transfinite induction on $\alpha$.]
\\[-.25cm]

Take now $\jS = \jP (\XX)$ and suppose that $\card(\jT) \leq \aleph_0$.  
Let $E$ be a nonempty subset of $\XX \times \YY$ $-$then given $\alpha$ $(0 < \alpha < \Omega)$, 
\[
E \in \jB^{[\alpha]} (\jP (\XX) \hsx \boxtimes \hsx \jT
\quad \text{iff} \quad
E_x \in \jB^{[\alpha]} (\jT) 
\qquad 
(\forall \ x \in \pi_\XX (E)).
\]
\\[-1.5cm]

[To discuss the nontrivial point, viz. that
\[
E_x \in \jB^{[\alpha]} (\jT) 
\ 
(\forall \ x \in \pi_\XX (E)) 
\ \implies \ 
E \in \jB^{[\alpha]} (\jP (\XX) \hsx \boxtimes \hsx \jT),
\]
one can argue by transfinite induction on $\alpha$, treating first the case when $\alpha = 1$ and then looking at the cases when 
$\alpha$ is odd or even separately.  
Here is the proof for $\alpha = 1$.  
Let $\jT = \{T_1, T_2, \ldots\}$ be an enumeration of $\jT$.  
Put
\[
S_i 
\ = \ 
\{x \in \pi_\XX (E) : T_i \subset E_x\}.
\]
Then
\[
E 
\ = \ 
\bigcup \ (S_i \times T_i) \in \jB^{[1]} \hsx (\jP (\XX) \hsx \boxtimes \hsx \jT.]
\]
\\[-.75cm]

\un{Ref}
R. Bing, W. Bledsoe, and R. Mauldin\footnote[2]{\vspace{.11 cm}\textit{\un{Pacific J. Math.}}, \textbf{51} (1974), pp. 27-36.}. 
\\[1cm]

\textbf{III} \ SETS GENERATED BY RECTANGLES

Let $\XX$ be a nonempty set $-$then, in \S5, we discussed the question: 
Is $\jP(\XX) \hsx \otimes \hsx \jP(\XX)$ dense in $\jP (\XX \times \XX)_\ts$?  
As has been seen there, the answer depends on the cardinality of $\XX$, the case of mystery being when 
$\aleph_1 < \card (\XX) \leq \fc$.  
\\[-.25cm]

If $\card (\XX) \leq \aleph_1$, then it is actually true that


\[
\jP (\XX \times \XX) 
\ = \ 
\jB^{[2]} (\jP(\XX) \hsx \boxtimes \hsx \jP(\XX)),
\]
i.e., each subset of $\XX \times \XX$ can be generated from the rectangles in just two steps.  
Assuming Martin's axiom, this conclusion remains in force if only $\card (\XX) \leq \fc$.  
\\[-.25cm]

On the other hand, the density of 
$\jP(\XX) \hsx \otimes \hsx \jP(\XX)$ 
in 
$\jP (\XX \times \XX) _\ts$
or still, the relation 
\[
\jP (\XX \times \XX) 
\ = \ 
\jP (\XX) \hsx \ovs{\otimes} \ \jP (\XX),
\]
is equivalent to the existence of a countable ordinal $\alpha \geq 2$ such that 
\[
\jP (\XX \times \XX) 
\ = \ 
\jB^{[\alpha]} (\jP(\XX) \hsx \boxtimes \hsx \jP(\XX)).
\]
\\[-.5cm]

\un{Ref}
Bing, W. Bledsoe, and R. Mauldin (op. cit.).
\\[-.25cm]

[Note: \ 
One could ask: Does
\[
\jP (\XX \times \XX) 
\ = \ 
\jP (\XX) \hsx \ovs{\otimes} \ \jP (\XX)
\]

\hspace{1.5cm} $\implies$
\[
\jP (\XX \times \XX) 
\ = \ 
\jB^{[2]} (\jP(\XX) \hsx \boxtimes \hsx \jP(\XX)) ?
\]
For a discussion of this question, see 
A. Miller\footnote[3]{\vspace{.11 cm}\textit{\un{Ann. Amer. Logic}}, \textbf{16} (1979), pp. 233-267.}. 
Consequences and implications may be found in 
R. Mauldin\footnote[4]{\vspace{.11 cm}\textit{\un{Fund. Math.}}, \textbf{95} (1977), pp. 129-139.}.]
\\[1cm]

\textbf{IV} \ POINT-FINITE CLASSES
\\[-.25cm]

Let $\XX$ be a nonempty set.  
Fix a subset $\jS$ of $\jP (\XX)$ containing $\emptyset$ and $\XX$.
\\[-.25cm]

A nonempty class $\jX \subset \jP (\XX)$  is said to be \un{point-finite} if each point of $\XX$ belongs to 
at most a finite number of elements of $\jX$.
\\[-.25cm]


(\un{\textbf{H}}) \ 
Suppose that $\jX$ is a point-finite class in $\XX$ such that 
$\jX_\Sigma \subset \jS_\tB$ $-$then, for some $\alpha < \Omega$, 
\[
\jX \subset \jB^{[\alpha]} (\jS).
\]

It will be simplest to examine first a special case.
\\

(\un{\textbf{P}}) \ 
Suppose that $\jX$ is a disjoint class in $\XX$ such that 
$\jX_\Sigma \subset \jS_\tB$ $-$then, for some $\alpha < \Omega$, 
\[
\jX \subset \jB^{[\alpha]} (\jS).
\]

[Proceed by contradiction $-$then there exist $\aleph_1$ disjoint subclasses  $\jX_\beta$ of $\jX$ 
such that for all $\alpha, \hsx \beta < \Omega$ $\jX_\beta \not\subset \jB^{[\alpha]} \hsx (\jS)$.  
Because 
$\jX_\Sigma \subset \jS_\tB$, 
there is a function 
$f : [0, \Omega[ \ra  [0, \Omega[$ such that 
$\bigcup \ \jX_\beta \in \jB^{[f(\beta)]}\hsx  (\jS)$ $(\beta \leq f(\beta))$.  
Choose $X_\beta \in \jX_\beta$, $X_\beta \notin \jB^{[f(\beta)]} (\jS)$.  
Put
$A = \bigcup\limits_{\beta < \Omega} \ X$ $-$then, 
for some $\alpha$, $A \in \jB^{[\alpha]} (\jS)$.  
But now 
$X_\alpha = A \cap \bigcup\jX_\alpha \in \jB^{[f(\alpha)]} (\jS)$, 
a contradiction.]
\\

\noindent
\un{Ref}
D. Preiss\footnote[2]{\vspace{.11 cm}\textit{\un{Comment. Math. Univ. Carolinae}}, \textbf{15} (1972), pp. 341-344.}. 
\\[-.5cm]

[The above proof is due to 
Fleissner\footnote[3]{\vspace{.11 cm}\textit{\un{Trans. Amer. Math. Soc.}}, \textbf{251} (1979), pp. 309-328.}.]
\\

In order to deduce (H) from (P), the following artifice will be needed.
\\[-.25cm]

{\small\bf Lemma} \ 
Let $\YY$ be a separable metric space with topology $\jT$.  
Suppose that 
$\{X(y) : y \in \YY\}$ is a point-finite class in $\XX$ such that 
\[
\{X(y) : y \in \YY\}_\Sigma \subset \jS_\tB.
\]
\un{Then}
\[
\{X(y) \times \{y\} : y \in \YY\}_\Sigma \ \subset \  \jS \hsx \boxtimes \hsx \jT)_\tB.
\]
\\[-1.5cm]

[Choose, as is possible, a basis $N_n$ $(n \in \N)$ for $\YY$ satisfying the diameter condition, 
i.e., $\diam(N_n) \ra 0$ and with the property that each point of $\YY$ belongs to $N_n$ for arbitrarily large values of $n$.  
Given a nonempty subset $Y_0$ of $\YY$, put
\[M_n 
\ = \ 
\bigcup \ \{\XX (y) : y \in N_n \cap Y_0\}.
\]
Then
\[
\bigcup \ \{\XX (y) \times \{y\} : y \in Y_0\} 
\ = \ 
\varlimsup \hsx (M_n \times N_n),
\]
hence is in $(\jS \hsx \boxtimes \hsx \jT)_\tB$.]
\\[-.25cm]

[\un{Proof} of (H) \ 
Proceed by contradiction $-$then
\[
\jX
\not\subset 
\jB^{[\alpha]} (\jS)
\qquad \forall \ \alpha < \Omega.
\]
Accordingly, one may select sets
\[
X_\alpha \in \jX - 
\left(
\jB^{[\alpha]} (\jS) \cup \{X_\beta : \beta < \alpha\}
\right)
\quad (\alpha < \Omega).
\]
Viewing $Y = \{\alpha : \alpha < \Omega\}$ as a subspace of $\R$, statement (P), 
in conjunction with the lemma supra, 
allows one to conclude that
\[
A 
\ = \ 
\bigcup \ \{X_\alpha \times \{\alpha\} : \alpha < \Omega\} 
\in 
\jB^{[\beta + 1]} (\jS \hsx \boxtimes \hsx \jT)
\]
for some $\beta > 1$.  Since


\[
\XX \times \{\alpha\} \in \jB^{[2]} (\jS \hsx \boxtimes \hsx \jT),
\]
it follows that 
\[
A 
\cap 
(\XX \times \{\alpha\})
\ = \ 
X_\alpha \times \{\alpha\} 
\in 
\jB^{[\beta + 1]} (\jS \hsx \boxtimes \hsx \jT).
\]
However (cf. Prob. II), this implies that $X_\alpha \in \jB^{[\beta + 1]} (\jS)$ $\forall \ \alpha < \Omega$, 
a contradiction.//
\\

It can be easily shown by example that statement (H) is no longer true if 
``point-finite'' is replaced by ``point-countable'' (defined in the obvious way).
\\[-.25cm]

\noindent
\un{Ref}
R. Hansell\footnote[2]{\vspace{.11 cm}\textit{\un{Proc. Amer. Math. Soc.}}, \textbf{83} (1981), pp. 375-378.}. 
\\[1cm]

\textbf{V} \ THEOREMS OF MILLER AND KUNEN 
\\[-.25cm]

Suppose that $\XX$ is a topological space for which 
$\jB\jo (\XX) = \jP (\XX)$.
Does $\exists$ an $\alpha < \Omega$ such that 
\[
\jB\jo (\XX)
\ = \ 
\Sigma_\alpha \ (\fS)?
\]
The answer, in general, is unknown.  
However, if $\XX$ is a metric space, then the response is positve.
\\[-.25cm]

{\small\bf Theorem} \ (Miller) \ 
Suppose that $\XX$ is a separable metric space for which 
$\jB\jo (\XX) = \jP (\XX)$ 
$-$then $\exists$ an $\alpha < \Omega$ such that 
\[
\jB\jo (\XX)
\ = \ 
\Sigma_\alpha \ (\fS).
\]

[First note that the cardinality of $\XX$ is necessarily $< \fc$.  
For otherwise, 
\[
\card (\jB\jo (\XX)) 
\ \geq \ 
2^\fc 
\ > \ 
\fc,
\]
which is impossible as there can be at most $\fc$ Borel sets in a separable metric space.  
If 
$\card(\XX) \leq \aleph_0$, then the assertion is clear.  
Let us consider the simplest nontrivial case, viz. when 
$\card(\XX) = \aleph_1$, referring the reader to the paper infra for the details when 
$\aleph_1 < \card(\XX) < \fc$.  
Write 
$\XX = \{x_\alpha : \alpha < \Omega\}$
and proceed by contradiction.  
For each $\alpha < \Omega$, let 
$A_\alpha \in \Sigma_{\alpha + 1} \ (\fS) - \Sigma_\alpha (\fS)$ 
and put 
$A = \{(x_\alpha, a) : a \in \A_\alpha\}$ 
$-$then it need only be shown that 
$A \in \Sigma_\beta \ (\fS \times \fS)$ 
for some $\beta < \Omega$ as this would entail
\[
A_{\beta + 1} 
\ = \ 
A \cap 
\left(
\{x_{\beta + 1}\} \times \XX
\right)
\in \Sigma_\beta (\fS).
\]
But, in view of the fact that $\XX$ is of cardinality $\aleph_1$ and of weight $\aleph_0$, we have

\[
\begin{matrix*}[l]
\jP (\XX) \ \ovs{\otimes} \ \jP (\XX)  \hspace{0.35cm} = \ \jP (\XX \times \XX) 
\\[4pt]
\hspace{0.5cm} \rotatebox{90}{=}  \hspace{1.25cm} \rotatebox{90}{=}  \hspace{2.4cm} \rotatebox{90}{=}
\\[4pt]
\jB\jo (\XX) \ \ovs{\otimes} \ \jB\jo (\XX)  \ = \ \jB\jo (\XX \times \XX) 
\end{matrix*}
,
\]
making the contention plain enough.]
\\

\noindent
\un{Ref}
A. Miller\footnote[2]{\vspace{.11 cm}\textit{\un{Ann. Math. Logic}}, \textbf{16} (1979), pp. 233-267.}. 
\\

[Note: \ 
Observe that the continuum hypothesis denies the existence of an uncountable separable metric space all of whose subsets are Borel.  
On the other hand, in the presence of Martin's axiom and the negation of the continuum hypothesis, 
it can be shown that there exists an uncountable set $X \subset \R$ in which every subset is an $F_\sigma$ 
(or, equivalently, $G_\delta$); 
cf. 
F. Tall\footnote[3]{\vspace{.11 cm}\textit{\un{Dissertationes Math.}}, \textbf{148} (1977), pp. 1-57.}.]
\\

{\small\bf Theorem} \ (Kunen) \ 
Suppose that $\XX$ is a metric space for which 
$\jB\jo (\XX) = \jP (\XX)$ 
$-$then $\exists$ an $\alpha < \Omega$ such that 
\[
\jB\jo (\XX)
\ = \ 
\Sigma_\alpha \ (\fS).
\]

[Kunen's proof is given in the paper of Miller cited above.
It runs as follows.
Because $\XX$ is a metric space, $\XX$ admits a $\sigma$-discrete basis 
$\fN =\bigcup  \{\fN_n  : n \in \N\}$.  
For each $N \in \fN$, let $\alpha (N)$ be the smallest ordinal $\alpha$ such that $\jP(N) = \Sigma_\alpha (\tr_N (\fS))$.  
Given $n \in \N$ and $\alpha < \Omega$, let
\[
C_{n, \alpha} 
\ = \ 
\{N \in \fN_n : \alpha (N) < \alpha\}.
\]
Claim: \ 
$\forall \ n$ $\exists \ \alpha(n)$ such that
\[
\card(\fN_n  - C_{n, \alpha(n)})
\ \leq \ 
\aleph_0.
\]
Indeed, if not, then for some $n$ it would be possible to find $A_\alpha$, $N_\alpha$ $(\alpha < \Omega)$ with:
\\[-.25cm]

\qquad (1) \quad 
$N_\alpha \in \fN_n$;
\\[-.25cm]

\qquad (2) \quad 
$N_\alpha \neq N_\beta$ $(\forall \ \alpha \neq \beta)$;
\\[-.25cm]

\qquad (3) \quad 
$A_\alpha \in \Sigma_{\alpha + 1} \ (\tr_{N_\alpha} (\fS)) - \Sigma_\alpha (\tr_{N_\alpha} (\fS))$.
\\[-.25cm]

\noindent
Since the union $\bigcup  A_\alpha$ cannot be Borel under these circumstances, we have a contradiction.  
The claim established, let $\alpha^* = \sup \{\alpha (n)\}$.  
Put
\[
X_0
\ = \ 
\XX 
\hsx - \hsx 
\bigcup \ 
\{N \in \fN : \alpha (N) < \alpha^*\}.
\]
Thanks to the claim, $X_0$ is a separable subspace of $\XX$, so, by Miller's theorem, 
$\exists \ \alpha_0 < \Omega$ such that 
$\jB\jo (\XX_0) = \Sigma_{\alpha_0} \ (\fS_0)$.  
If now 
$\alpha = \sup \{\alpha_0, \alpha^* + 1\}$, then 
$\jB\jo (\XX) = \Sigma_{\alpha} \ (\fS)$.] 
\\[1cm]

\textbf{VI} \ POINT-FINITE CLASSES (BIS)
\\[-.25cm]

As in Prob. IV, let $\XX$ be a nonempty set.  
Fix a subset $\jS$ of $\jP (\XX)$ containing $\emptyset$ and $\XX$.
\\[-.25cm]

Suppose that $\jX$ is a point-finite class in $\XX$ such that $\jX_\Sigma \subset \jS_\tB$ $-$then, 
as seen above, $\jX$ is contained in $\jB^{[\alpha]} \hsx (\jS)$ for some $\alpha < \Omega$.  
We now ask: 
Does there exist an $\alpha < \Omega$ such that $\jX_\Sigma \subset \jB^{[\alpha]} \hsx (\jS)$?
\\[-.25cm]

To give an answer, write $\jX = \{\XX_i : i \in I\}$ $-$then there will be an $\alpha$ with the stated property if 
$\exists$ an uncountable set $J$ such that 

\[
\jP (I \times J) 
\ = \ 
\jP (I) \hsx \ovs{\otimes} \ \jP (J).
\]
\\[-.75cm]

[The proof is similar to that of statement (H) in Prob. IV, modulo an appropriate variant of the lemma appearing there.]
\\[-.25cm]

The question of the equality
\[
\jP (I \times J) 
\ = \ 
\jP (I) \hsx \ovs{\otimes} \ \jP (J)
\]
has been considered in Exer. 5 (\S5).  
Recall that it will hold if both $\card(I)$ and $\card(J)$ are $\leq \aleph_1$ 
(or even $\leq \fc$ if Martin's axiom is assumed).  
Consequently, the answer to the question supra is affirmative if 
$\card(I) \leq \aleph_1$. 
\\[-.25cm]

There is another condition on J which will force the equality
\[
\jP (I \times J) 
\ = \ 
\jP (I) \hsx \ovs{\otimes} \ \jP (J),
\]
namely that $\jP (J)$ be generated as a $\sigma$-algebra by a set of cardinality $\leq \aleph_0$.
\\[-.25cm]

[Use Exer. 9, Prob., V, and Prob. II.] 
\\

\textbf{Example} \ 
Let $\XX$ be a metric space.  
Take, in this context, $\jS = \jG$.  
Suppose that $\jX = \{\{x\} : x \in \XX\}$.  
If $\jX_\Sigma \subset \jS_\tB$, then $\jB \jo (\XX) = \sP (\XX)$, so in this case we are back in the setting of Prob. V. 
Assume now that $\XX$ is, in addition, separable.  
Let $\jX$ be a point-finite class in $\XX$ such that $\jX_\Sigma \subset \jS_\tB$ $-$then there exists an $\alpha < \Omega$ 
such that $\jX_\Sigma \subset \jB^{[\alpha]} \hsx (\jS)$.  
This, of course, is obvious if $\card(I) \leq \aleph_0$. 
On the other hand, if $\card(I) > \aleph_0$, fix a point $x_i $ in each $\XX_i$ $-$then
$X_I = \{x_i : i \in I\}$ is an uncountable separable metric space all of whose subsets are Borel, hence
\[
\jP (I \times \XX_I) 
\ = \ 
\jP (I) \hsx \ovs{\otimes} \ \jP (\X_I).
\]
\\[-.75cm]

\noindent
\un{Ref}
R. Hansell\footnote[2]{\vspace{.11 cm}\textit{\un{General Topology and Modern Analysis}}, 
Springer-Verlag, Berlin (1981), pp. 405-416.}. 
\\[1cm]

\textbf{VII} \ ZERO SETS IN UNIFORM SPACES 
\\[-.25cm]

Let $\XX$ be a uniform space $-$then the class $\jZ$ of zero sets of the bounded uniformly continuous functions 
$f : \XX \ra \R$ has the following properties
\[
\begin{matrix*}[l]
&(1) \quad \emptyset, \hsx \XX \in \jZ, 
&&&(2) \quad \jZ = \jZ_\ts, 
\\[4pt]
&(3) \quad \jZ = \jZ_\delta, 
&&&(4) \quad \jZ_\tc = \jZ_\sigma.\hspace{2cm}{\ }
\end{matrix*}
\]
In addition, given disjoint $Z_1$, $Z_2 \in \jZ$, there exist disjoint $U_1$, $U_2 \in \jZ_\tc$ such that 
\[
Z_1 \subset U_1, 
\quad 
Z_2 \subset U_2.
\]

One has: 
\[
\jZ_{\tB_\td}
\ = \ 
\text{$\sigma$-$\jAlg (\jZ)$.}
\]

[This can be seen by repeating the argument for its topological analogue virtually word-for-word.]

\noindent
\un{Ref}
J. Jayne\footnote[3]{\vspace{.11 cm}\textit{\un{Proc. Prague Symp. General Topology}}, Part B
(1976), pp. 187-194.}. 
\\[-.25cm]

Let $\XX$ be a nonempty set; let $\jZ$ be a class of subsets of $\XX$ possessing the five properties supra $-$then 
$\XX$ can be equipped with the structure of a uniform space with respect to which $\jZ$ is precisely the class of zero sets of 
the bounded uniformly continuous functions $f : \XX \ra \R$.  
Consequently
\[
\jZ_{\tB_\td}
\ = \ 
\text{$\sigma$-$\jAlg (\jZ)$.}
\]
\\[-1.25cm]

[In this connection, recall that a topology $\jT$ on $\XX$ is the uniform topology for some uniformity on $\XX$ 
iff the topological space $(\XX, \jT)$ is completely regular.]
\\[-.25cm]

\noindent
\un{Ref}
H. Gordon\footnote[5]{\vspace{.11 cm}\textit{\un{Pacific J. Math.}}, \textbf{36} (1971), pp. 133-157.}. 
\\[1cm]

\textbf{VIII} \ DISJOINT GENERATION
\\[-.25cm]

Let $\XX$ be a nonempty set; let $\jS$ be a nonempty class of subsets of $\XX$ such that 
\[
\jS \ = \ \jS_\ts, 
\quad 
\jS_\tc \subset \jS_\sigma.
\]
Suppose in addition that given disjoint $S_1$, $S_2 \in \jS$, there exist disjoint $C_1$, $C_2 \in \jS_\tc$ such that 
\[
S_1 \subset C_1, 
\quad 
S_2 \subset C_2.
\]
Then 
\[
\jS_{\tB_\td}
\ = \ 
\text{$\sigma$-$\jAlg (\jS)$}.
\]


[According to Lemma 5, it suffices to prove that $\jS_\tc \subset \jS_{\tB_\td}$.  
For this purpose, show by a direct set-theoretic construction that
\[
\jS_\tc 
\subset
\jS_{\delta \hsy \sigma_\td \hsy \delta \hsy \sigma_\td}.]
\]
\\[-1.25cm]

\noindent
\un{Ref}
J. Jayne\footnote[2]{\vspace{.11 cm}\textit{\un{Mathematika}},\textbf{24} (1977), pp. 241-256.}. 
\\[1cm]

\textbf{IX} \ INCREASING AND DECREASING LIMITS
\\[-.25cm]

Let $\jS \subset \jP (\XX)$ be nonempty.  Write
\[
\begin{cases}
\ (\uparrow) \ (\jS)\\
\ (\downarrow) \ (\jS)
\end{cases}
\]
for the class of all subsets of $\XX$ which are the limit of an 
\[
\begin{cases}
\ \text{increasing}\\
\ \text{decreasing}
\end{cases}
\]
sequence of sets in $\jS$.
\\[-.25cm]

Suppose now that $\jS$ is a lattice.
Put
\[
\Xi^0 (\jS) 
\ = \ 
(\uparrow) \ (\jS), 
\quad 
\Xi_0 (\jS) 
\ = \ 
(\downarrow) \ (\jS), 
\]
and define via transfinite recursion the classes $\Xi^\alpha (\jS) $, $\Xi_\alpha (\jS) $ by writing
\[
\begin{cases}
\ \Xi^\alpha (\jS)  
\ = \ 
(\uparrow) \ \big(\bigcup\limits_{\beta < \alpha} \ \Xi_\beta (\jS)\big)
\\[18pt]
\ \Xi_\alpha (\jS)  
\ = \ 
(\downarrow) \ \big(\bigcup\limits_{\beta < \alpha} \ \Xi^\beta (\jS)\big)
\end{cases}
\qquad 
(\alpha < \Omega).
\]


Investigate these classes.
\\[-.25cm]

\noindent
\un{Ref}
W. Sierpi\'nski\footnote[2]{\vspace{.11 cm}\textit{\un{Fund. Math.}}, \textbf{18} (1932), pp. 1-22.}. 
\\[1cm]

\textbf{X} \ $\aleph$-OPERATIONS
\\[-.25cm]

Let $\aleph$ be an infinite cardinal.  
Consider a map
\[
M \hsx : \hsx \jP (\jP(\XX)) \ra \jP (\jP(\XX))
\]
with the following properties:
\\[-.25cm]

\qquad (1) \quad 
If $f : \XX \ra \XX$ is a function and if $\jS \subset \jP(\XX)$ is a class, then 
\[
f^{-1} (M(\jS)) 
\subset 
M (f^{-1}(\jS));
\]
\\[-1cm]

\qquad (2) \quad 
If $\jS^\prime$, $\jS^{\prime\prime} \subset \sP(\XX)$ are classes, 
if $M(\jS^\prime) \subset M(\jS^{\prime\prime})$, and if $\jS^{\prime\prime} \in M(\jS^{\prime\prime})$,
then 
\[
M (\jS^\prime \cup \{\jS^{\prime\prime}\}) 
\subset 
M (\jS^{\prime\prime}).
\]
\\[-.75cm]

Under these circumstances, $M$ is said to be an \un{$\aleph$-operation} if for every initial ordinal $\zeta$ with 
$\card (\zeta) \leq \aleph$ and if for any increasing transfinite $\zeta$-sequence 
$\{\jS_\alpha : \alpha < \zeta\}$, the inclusions
\[
M(\jS_\alpha) 
\subset 
M(\jS)
\quad (\alpha < \zeta) 
\quad \implies \quad 
M
\bigg(
\bigcup\limits_{\alpha < \zeta} \ \jS_\alpha 
\bigg)
\hsx \subset \hsx 
M(\jS).
\]
\\[-.75cm]

Illustrate this concept by examining the various set-theoretic operations which have been discussed in this \S.
\\[-.25cm]

If $\star$ is extensionally attainable, then is it necessarily true that $M_\star$ is an $\aleph$-operation?

\noindent
\un{Ref}
M. Ershov\footnote[2]{\vspace{.11 cm}\textit{\un{SLN}}, \textbf{794} (1980), pp. 105-111.}. 
\\[-.25cm]

[Here also may be found a number of selection theorems of substantial generality.]

%
\printindex
\end{document}